\documentclass[11pt]{article}
\usepackage{fullpage}

\date{}

\usepackage{tikz}
\usetikzlibrary{calc}
\usetikzlibrary{matrix}
\usepackage{amsmath,amssymb,amsfonts,amsthm}
\usepackage{color}                    
\usepackage[pdftex]{hyperref}

\usepackage{graphicx,subcaption}
\allowdisplaybreaks

\providecommand{\V}[1]{\boldsymbol{#1}}
\providecommand{\Th}{\mathcal{T}_h}
\providecommand{\Tc}{\mathcal{T}_{c}}
\providecommand{\hTc}{\hat{\mathcal{T}}_{c}}
\providecommand{\M}{\mathbb{M}}

\providecommand{\J}{\mathbb{J}}

\DeclareMathOperator{\tr}{tr}

\newcommand{\bey}{\begin{eqnarray}}
\newcommand{\eey}{\end{eqnarray}}

\newcommand{\beq}{\begin{equation}}
\newcommand{\eeq}{\end{equation}}
\theoremstyle{plain}

\theoremstyle{definition}

\theoremstyle{remark}
\newtheorem{exam}{\hspace{6mm}Example}[section]
\newtheorem{rem}{\hspace{6mm}Remark}[section]

\title{An adaptive moving mesh finite element solution of the Regularized Long Wave equation}

\author{Changna Lu\thanks{College of Mathematics and Statistics,
Nanjing University of Information Science and Technology,
Nanjing, Jiangsu 210044, China. ({\tt luchangna@nuist.edu.cn})},
\and Weizhang Huang\thanks{
Department of Mathematics, the University of Kansas, Lawrence, KS 66045,
U.S.A. ({\tt whuang@ku.edu})},
\and
Jianxian~Qiu\thanks{School of Mathematical Sciences and Fujian Provincial
Key Laboratory of Mathematical Modeling \& High-Performance
Scientific Computing, Xiamen University,
Xiamen, Fujian 361005, China. ({\tt jxqiu@xmu.edu.cn})}
}
\begin{document}
\vskip 1cm
\maketitle

\begin{abstract}
An adaptive moving mesh finite element method is proposed for the numerical solution of
the regularized long wave (RLW) equation. A moving mesh strategy based on the so-called
moving mesh PDE is used to adaptively move the mesh to improve computational accuracy
and efficiency. The RLW equation represents a class of partial differential equations
containing spatial-time mixed derivatives. For the numerical solution of those equations, a $C^0$
finite element method cannot apply directly on a moving mesh since the mixed derivatives of the finite
element approximation may not be defined. To avoid this difficulty, a new variable is
introduced and the RLW equation is rewritten into a system of two coupled equations.
The system is then discretized using linear finite elements in space and the fifth-order Radau IIA scheme
in time.  A range of numerical examples in one and two dimensions, including
the RLW equation with one or two solitary waves and special initial conditions that lead to
the undular bore and solitary train solutions, are presented. Numerical results demonstrate that
the method has a second order convergence and is able to move and adapt
the mesh to the evolving features in the solution.
\end{abstract}



\noindent{\bf Abbreviated title.}
An adaptive moving mesh FE solution of RLW

\section{Introduction}

We consider the adaptive moving mesh finite element (FE) solution of
the regularized long wave (RLW) equation
(which is also called the Benjamin-Bona-Mahony or BBM equation)
 in one and two dimensions.
The initial-boundary value problem of the 2D RLW equation \cite{Goldstein-1985,Calvert-1976,Goldstein-1980}
reads as
\begin{equation}
\begin{cases}
u_t + \alpha u_x + \beta u_y + \gamma u u_x + \delta u u_y - \mu u_{xxt} - \mu u_{yyt}  = 0,
&\quad (x,y)\in \Omega,\; t \in (0,T]
\\
u(x,y,t) = g(x,y,t), &\quad (x,y)\in \partial \Omega,\; t \in (0,T]
\\
u(x,y,0) = u_0(x,y),& \quad (x,y)\in \Omega
\end{cases}
\label{RLW-1}
\end{equation}
where $\Omega$ is a bounded polygonal domain and $\alpha$,  $\beta$,  $\gamma$, $\delta$, $\mu$
are constants with $|\gamma| + |\delta| > 0$, and $\mu>0$,
and $u_0$ and $g$ are given functions.
The RLW equation has been used to model ion acoustic waves and
magnetohydrodynamics waves in plasmas,
longitudinal dispersive waves in elastic rods,
pressure waves in liquid gas bubbles, 
and nonlinear transverse waves in shallow water; for example see
\cite{Benjamin1972,Bona2000,Peregrine-1966}.
The RLW equation was proposed first by Peregrine \cite{Peregrine-1966} and later by
Benjamin et al. \cite{Benjamin1972} as a model for small amplitude long waves on the surface of water in a channel.
Generalizations such as the generalized regularized long wave equation (gRLW)
or the modified regularized long wave equation (MRLW) \cite{Goldstein-1985,Calvert-1976,Goldstein-1980}
and generalized Rosenau-Kawhara-RLW equation \cite{Rosenaua1986}
also arise from various applications.

The RLW equation is related to the Korteweg-de Vries (KdV) equation but has distinct features.
For example, Medeiros and Miranda \cite{Medeiros1977} discuss the problem of periodic solution and show that
RLW can almost cover all the application of KdV.
On the other hand,  Olver \cite{Olver1979} proves that RLW can have only three non-trivial
independent conservation laws. This is very different from KdV which is known to
have an infinite number of conservation laws.
Moreover, KdV is known to possess single and multiple solitons that maintain their shapes
and velocities after their interactions and can have inelastic collision.
RLW does not appear to admit an inverse-scattering theory which would lead to
an analytical representation for solitary wave solutions.
Nevertheless, the initial-value problem of RLW posed on the whole real line still has the property
that initial disturbances resolve into a train of solitary waves and a dispersive tail (e.g., see \cite{Bona1981};
also see Examples~\ref{exam3.4} and \ref{exam3.7} in section~\ref{SEC:numerics}).
Much effort has been made to understand whether or not RLW has the characteristics
of solitons. For example, Abdulloev \cite{Abdulloev1976} shows that
two solitons of RLW can have inelastic collision. Analytical solutions
have also been obtained by various researchers; e.g., see \cite{Ma1996,Matveev1991,Tian2003,Yang2012}.

The numerical solution of the RLW equation and its variants and generalizations have been
considered extensively in recent decades.  Among many existing works, we mention
Eilbeck and Mc{G}uire \cite{Eilbeck1975,Eilbeck1977} (finite difference methods),
Guo and Cao \cite{Guo1988} (a Fourier pseudospectral method with a restrain operator),
Luo and Liu \cite{Luo-1999} (a mixed Galerkin),
Zaki \cite{Zaki2001} (combined splitting with cubic B-spline FEM),
Dogan \cite{Dogan2002} (linear FEM),
Da{\v{g}} et al. \cite{Dag2004} (cubic B-spline collation),
Gu and Chen \cite{Gu-2008} (a least squares mixed Galerkin),
Gao et al. \cite{Qiu2009} (local Discontinuous Galerkin),
Mei et al. \cite{Gao2015,Mei-2012b,Mei2012,Mei-2015} (mixed Galerkin),
and Siraj-ul-Islam et al. \cite{SirajulIslam2009}  (meshfree method).
These works are for RLW, gRLW, or MRLW in 1D, and much less work has been done in 2D.
Dehghan and Salehi \cite{Dehghan2011} consider the numerical solution of 2D RLW
in fluids and plasmas using the boundary knot method (a meshless boundary-type
radial basis function collocation technique).

The objective of this paper is twofold. The first is to study the numerical solution of RLW using an adaptive
moving mesh method. The method works for a general spatial dimension but we focus
only on 1D and 2D in this work. As will be seen in section~\ref{SEC:numerics},
a large spatial domain often has to be used in the numerical solution to
reduce the boundary effects and to cover the evolving
features for the whole time period under consideration. This requires a large
number of mesh elements for a reasonable level of computational accuracy
especially in multi-dimensions. To improve computational efficiency, it is natural
to employ an adaptive moving mesh technique which dynamically adapts
the mesh to the local, evolving features in the solution of RLW.
In this work, we will employ the so-called moving mesh PDE (MMPDE) method \cite{HRR94b,HRR94a,HR11}
that moves the mesh continuously in time and orderly in space
using a PDE formulated as the gradient flow equation of a meshing
functional. We will use a newly developed discretization of the MMPDE \cite{HK2014}
that makes the implementation of the MMPDE method
not only significantly simpler in multi-dimensions but also much more reliable in the sense
that there is a theoretical guarantee for mesh nonsingularity.

The second objective of the paper is to study how to discretize space-time mixed derivatives
using finite elements on moving meshes. RLW (\ref{RLW-1}) represents a class of
PDEs containing space-time mixed derivatives. In addition to RLW, this class includes 
Boussinesq \cite{Calogero1982}, modified Buckley-Leverett \cite{Spayd2011},
and Sobolev \cite{Showalter1972} equations. A feature of these PDEs is that space-time mixed derivatives
are involved in their both strong and weak formulations. When the mesh is moving,
these derivatives of a $C^0$ finite element approximation are not defined (cf. section~\ref{SEC:fem}).
There are various ways to overcome this difficulty. We utilize a new variable (see (\ref{v-1}) below)
and demonstrate numerically that the resulting linear finite element discretization gives a second
order convergence on moving meshes. Since (\ref{v-1}) is not tailored to the special structure
of RLW, we may expect that this idea of treating space-time mixed derivatives
can also be used for the moving mesh solution of Boussinesq,
modified Buckley-Leverett, and Sobolev equations.

It is worth mentioning that a number of moving mesh methods have been developed in the past
and there is a considerable literature in the area. Instead of going over the literature, we refer
the interested reader to the books/review articles \cite{Bai94a,Baines-2011,BHR09,HR11,Tan05}
and references therein.

An outline of the paper is as follows. The adaptive moving mesh finite element method is described
in Section~\ref{SEC:method}. The transformation of RLW into a system
of two coupled PDEs, the discretization of the PDE system on moving meshes via linear finite elements,
and the conservation laws possessed by RLW are discussed in the section.
The generation of adaptive moving meshes using a new implementation of the MMPDE method
is discussed in section~\ref{SEC:mmpde}.
1D and 2D numerical examples of RLW (and MRLW) are presented in section~\ref{SEC:numerics}.
Finally, section~\ref{SEC:conclusion} contains conclusions and further comments.

\section{An adaptive moving mesh finite element method}
\label{SEC:method}

In this section we describe the adaptive moving mesh FE method for the numerical solution of
the RLW equation. We first describe the basic procedure of the method and then elaborate on
the linear FE discretization of the RLW equation on moving meshes, followed by
a discussion on the conservation laws possessed by the RLW equation.
An MMPDE-based moving mesh strategy will be discussed in the next section.
To be specific, we describe the method in two dimensions.
The one dimensional formulation is similar.

We start with introducing a new variable
\begin{equation}
v = u - \mu u_{xx} - \mu u_{yy}
\label{v-1}
\end{equation}
and rewriting (\ref{RLW-1}) into 
\begin{equation}
\begin{cases}
v_t + \alpha u_x + \beta u_y
+ \gamma u u_x + \delta u u_y  = 0,
&\quad (x,y)\in \Omega, \quad t \in (0, T]
\\
v = u - \mu u_{xx} - \mu u_{yy}, &\quad (x,y)\in \Omega, \quad t \in (0, T]
\\
u = g, &\quad (x,y)\in \partial \Omega, \quad t \in (0, T] .
\end{cases}
\label{RLW-2}
\end{equation}
The weak formulation is to find $u(\cdot, t) \in H^1(\Omega)\cap \{ u |_{\partial \Omega} = g \}$
and $v(\cdot, t) \in H^1(\Omega)$ for $ 0 < t \le T$ such that
\begin{equation}
\begin{cases}
 \int_\Omega \left (v_t + \alpha u_x + \beta u_y
+ \gamma u u_x + \delta u u_y\right ) \phi d x d y  = 0,
&\quad \forall \phi \in H^1(\Omega), \quad t \in (0, T]
\\
 \int_\Omega \left ( (v - u) \psi  - \mu u_x\psi_x
 - \mu u_y \psi_y \right ) d x d y = 0,
&\quad \forall \psi \in H_0^1(\Omega) , \quad t \in (0, T].
\end{cases}
\label{RLW-2-w}
\end{equation}

The basic procedure of the the adaptive moving mesh FE method for solving (\ref{RLW-2-w}) is as follows.
\begin{enumerate}
\item Given an initial mesh $\Th^0$ and an initial time step $\Delta t_0$.
\item For $n = 0, 1, ...$
	\begin{enumerate}
	\item An MMPDE-based moving mesh strategy (cf. section~\ref{SEC:mmpde}) is used to generate
		the new mesh $\Th^{n+1}$ based on the current mesh $\Th^n$ and the numerical solution
		$u_h^n \approx u(\cdot, t_n)$ defined thereon. Note that $\Th^{n+1}$ and $\Th^{n}$
		have the same number of the elements ($N$), the same number of the vertices ($N_v$), and the same
		connectivity. They differ only in the location of the vertices, $(x_i, y_i)$, $i = 1, ..., N_v$.
		
	\item For $t \in [t_n, t_{n+1}]$ with $t_{n+1} =t_n + \Delta t_n$,
		the coordinates and velocities of the vertices are defined as
		\begin{align*}
                 & x_i(t) = \frac{t^{n+1}-t}{\Delta t_n}x_i^n + \frac{t-t^n}{\Delta t_n}x_i^{n+1},\quad 
			y_i(t) = \frac{t^{n+1}-t}{\Delta t_n}y_i^n + \frac{t-t^n}{\Delta t_n}y_i^{n+1},\quad i = 1, ..., N_v
		\\
		& \dot{x}_i(t) = \frac{x_i^{n+1}-x_i^n}{\Delta t_n},\quad 
		\dot{y}_i(t) = \frac{y_i^{n+1}-y_i^n}{\Delta t_n},\quad i = 1, ..., N_v .
		\end{align*}
		The corresponding mesh is denoted by $\Th(t)$ ($ t_n \le t \le t_{n+1}$).
	\item The RLW equation (\ref{RLW-2-w}) is discretized in space using linear finite elements and then integrated
		in time for one step using the fifth-order Radau IIA method (e.g., see Hairer and Wanner \cite{HW96}).
		A standard procedure is used for the selection of the time step size, together with a two-step error estimator 
		of Gonz\'{a}lez-Pinto et al. \cite{Montijano2004}. If the actual step size
		(denoted by $\widetilde{\Delta t_n}$)
		is smaller than $\Delta t_n$, the time and mesh are updated as
		\[
		t_{n+1} \leftarrow t_n + \widetilde{\Delta t_n},\quad
		x_i^{n+1} \leftarrow x_i^n + \widetilde{\Delta t_n} \dot{x}_i,\quad 
		y_i^{n+1} \leftarrow y_i^n + \widetilde{\Delta t_n} \dot{y}_i ,\quad i =1, ..., N_v.
		\]
		The predicted time step size will be used as $\Delta t_{n+1}$.
	\end{enumerate}
\end{enumerate}

The FEM discretization of the RLW equation on $\Th(t)$ is discussed in the next subsection
while the generation of $\Th^{n+1}$ using the MMPDE-based moving mesh
strategy will be given in section~\ref{SEC:mmpde}.

\subsection{Linear finite element discretization on $\Th(t)$}
\label{SEC:fem}

For notational simplicity, we assume that the vertices of $\Th(t)$ are ordered in a way that the first $N_{vi}$
vertices are interior vertices. Let $\phi_i = \phi_i(x,y,t)$ be the linear basis function associated with
the $i$-th vertex $(x_i, y_i)$. Define
\begin{align}
& V^h(t) = \text{span}\{ \phi_1, ..., \phi_{N_{v}}\} ,
\label{V-1}
\\
& V_0^h(t) = V^h(t) \cap \{ v|_{\partial \Omega} = 0\} \equiv \text{span}\{ \phi_1, ..., \phi_{N_{vi}}\},
\label{V-2}
\\
& V_g^h(t) = V^h(t) \cap \{ v(x_i, y_i, t) = g(x_i, y_i, t), \, i = N_{vi}+1, ..., N_v\} .
\label{V-3}
\end{align}
The linear finite element approximation of (\ref{RLW-2-w}) is
to find $u_h(\cdot, t) \in V_g^h(t)$ and $v_h(\cdot, t)\in V^h(t)$, $t \in (0, T]$ such that
\begin{equation}
\begin{cases}
\int_\Omega\left ( \frac{\partial v_h}{\partial t}
 + \alpha \frac{\partial u_h}{\partial x}
+ \beta \frac{\partial u_h}{\partial y} + \gamma u_h \frac{\partial u_h}{\partial x}
+ \delta u_h  \frac{\partial u_h}{\partial y}\right ) \phi \; d x d y = 0, \quad \forall \phi \in V^h(t),\quad
t \in (0, T]
\\
\int_\Omega \left ((v_h - u_h) \psi  - \mu \frac{\partial u_h}{\partial x}
\frac{\partial \psi}{\partial x} - \mu  \frac{\partial u_h}{\partial y} \frac{\partial \psi}{\partial y}\right )
 d x d y = 0,\quad \forall \psi \in V_0^h(t),\quad t \in (0, T] .
\end{cases}
\label{fem-1}
\end{equation}

To cast (\ref{fem-1}) in a matrix form, we express $u_h$ and $v_h$ as
\begin{equation}
u_h = \sum_{i=1}^{N_v} u_i(t) \phi_i(x, y, t),
\quad v_h = \sum_{i=1}^{N_v} v_i(t) \phi_i(x, y, t),
\label{uv-1}
\end{equation}
subject to the boundary condition
\begin{equation}
u_i = g(x_i, y_i, t),\quad i = N_{vi}+1, ..., N_{v}.
\label{bc-2}
\end{equation}
Notice that
\[
\frac{\partial v_h}{\partial t} = \sum_{i=1}^{N_v} \frac{d v_i}{d t}  \phi_i
+ \sum_{i=1}^{N_v} v_i \frac{\partial \phi_i}{\partial t} .
\]
It is not difficult to show (e.g., see Jimack and Wathen \cite{Jimack-1991}) that
\beq
\frac{\partial \phi_i}{\partial t} = - \nabla \phi_i \cdot \dot{\V{X}}, \quad \text{a.e. in } \Omega
\eeq
where $\dot{\V{X}}$ is a piecewise linear mesh velocity defined by
\beq
\dot{\V{X}} = \sum_{i=1}^{N_v} \begin{bmatrix} \dot{x}_i \\ \dot{y}_i \end{bmatrix} \phi_i .
\label{Xdot}
\eeq
Using this we can rewrite $\partial v_h/\partial t$ as
\begin{equation}
\frac{\partial v_h}{\partial t} = \sum_{i=1}^{N_v} \frac{d v_i}{d t}  \phi_i
- \nabla  v_h \cdot \dot{\V{X}} .
\label{vt-1}
\end{equation}
Inserting (\ref{uv-1}), (\ref{bc-2}), and (\ref{vt-1}) into (\ref{fem-1}) and taking $\phi = \phi_i$ ($i = 1, ..., N_{v}$)
and $\psi = \phi_i$ ($i = 1, ..., N_{vi}$) successively, we get
\begin{equation}
\begin{cases}
\begin{bmatrix} M_{II} & M_{IB} \\ M_{BI} & M_{BB} \end{bmatrix} \frac{d }{d t}\begin{bmatrix} \V{v}_I \\ \V{v}_B \end{bmatrix}
+ \begin{bmatrix} \V{f}_I \\ \V{f}_B \end{bmatrix} = 0,
\\
\begin{bmatrix} M_{II} & M_{IB} \end{bmatrix} \left (\begin{bmatrix} \V{v}_I \\ \V{v}_B \end{bmatrix}
- \begin{bmatrix} \V{u}_I \\ \V{u}_B \end{bmatrix} \right )
- \begin{bmatrix} A_{II} & A_{IB} \end{bmatrix} \begin{bmatrix} \V{u}_I \\ \V{u}_B \end{bmatrix} = 0 ,
\\
\V{u}_B = \V{g}_B,
\end{cases}
\label{fem-2}
\end{equation}
where the vectors and matrices in (\ref{fem-2})  
are partitioned according to the entries associated with the interior vertices
(with symbol ``$I$'') and those associated with the boundary vertices (with symbol ``$B$''),
$\V{u} = (u_1, ..., u_{N_{vi}}, ..., u_{N_v})^T$ and $\V{v} = (v_1, ..., v_{N_{vi}}, ..., v_{N_v})^T$ are the unknown
vectors, $M$ and $A$ are the mass and stiffness matrices, respectively, and the entries of
$M$, $A$, $\V{f}$, and $\V{g}$ are given by
\begin{equation}
\begin{cases}
& M_{i,j} = \int_\Omega \phi_i \phi_j d x d y,\quad
A_{i,j} = \int_\Omega \left (\mu \frac{\partial \phi_i}{\partial x} \frac{\partial \phi_j}{\partial x} 
+ \mu \frac{\partial \phi_i}{\partial y} \frac{\partial \phi_j}{\partial y} \right ) d x d y,
\\
& f_{i} = \int_\Omega \left (\alpha \frac{\partial u_h}{\partial x}
+ \beta \frac{\partial u_h}{\partial y} + \gamma u_h \frac{\partial u_h}{\partial x}
+ \delta u_h  \frac{\partial u_h}{\partial y} - \nabla u_h \cdot \dot{{\V{X}}} \right ) \phi_i   d x d y ,
\\
& g_i = g(x_i, y_i, t).
\end{cases}
\label{fem-2b}
\end{equation}

When the mesh is fixed, both the mass and stiffness matrices are time independent. In this case,
by differentiating the second equation of (\ref{fem-2}) with respect to time and
subtracting it from the first equation we get
\begin{equation}
(M_{II} + A_{II}) \frac{d \V{u}_I}{d t} + \V{f}_I = 0 .
\label{fem-3}
\end{equation}
Since both $M_{II}$ and $A_{II}$ are symmetric and positive definite, $M_{II} + A_{II}$ is invertible
and (\ref{fem-3}) forms an ODE system. As a result, the solution existence and uniqueness of (\ref{fem-2})
can be derived from that of the ODE system (\ref{fem-3}). Moreover, it is not difficult to show
that (\ref{fem-3}) can also be obtained by applying the linear finite element discretization directly
to the original equation (\ref{RLW-1}).

When the mesh is varying with time, both $M$ and $A$ depend on time too. In this case,
(\ref{fem-2}) cannot reduce to (\ref{fem-3}) in general. Nevertheless, from the second and third equations
of (\ref{fem-2}) we get
\begin{equation}
\V{u}_I = (M_{II}+A_{II})^{-1}  \left ( M_{II} \V{v}_I + M_{IB} \V{v}_B  - (M_{IB} + A_{IB}) \V{g}_B\right ) .
\label{fem-4}
\end{equation}
Notice that $\V{f}$ is a function of $\V{u} = (\V{u}_I^T, \V{u}_B^T)^T$ and can be written as
$\V{f} = \V{f}(\V{u}_I, \V{u}_B)$. Inserting (\ref{fem-4}) into the first equation of (\ref{fem-2}) we obtain
\begin{equation}
M \frac{d \V{v}}{d t}
+ \V{f}((M_{II}+A_{II})^{-1}  \left ( [M_{II}\; M_{IB}] \V{v}  - (M_{IB} + A_{IB}) \V{g}_B\right ) , \V{g}_B) = 0,
\label{fem-5}
\end{equation}
which is also an ODE system. Then, the solution existence and uniqueness of (\ref{fem-2}) can be derived
from that of the ODE system (\ref{fem-5}). Once $\V{v}$ has been obtained, we can find $\V{u}_I$
from (\ref{fem-4}).

In our computation, (\ref{fem-2}) is solved directly, which is a DAE (differential-algebraic equation) system.
It is integrated using the fifth-order Radau IIA method with a variable step size determined by a two-step error estimator 
of Gonzalez-Pinto et al. \cite{Montijano2004}.

\begin{figure}[tb]
\centering
\begin{tikzpicture}[scale = 1]
\draw [thick, ->] (0,0) -- (8,0);
\draw [thick, ->] (0,3) -- (8,3);
\draw[right] (8.1,0) node {$t^n$};
\draw[right] (8.1,3) node {$t^{n+1}$};

\draw[thick] (1, 0) -- (2,3);
\draw[thick] (3, 0) -- (4.2,3);
\draw[thick] (4.5, 0) -- (6.3,3);

\draw[right] (3.5,4) node {$t$};
\draw[thick] (3.5, -1) -- (3.5,4);

\fill (canvas cs:x=1cm,y=0cm) circle (1.5pt);
\draw[below] (1,0) node {$x_{j-1}^n$};
\fill (canvas cs:x=3cm,y=0cm) circle (1.5pt);
\draw[below] (3,0) node {$x_j^n$};
\fill (canvas cs:x=4.5cm,y=0cm) circle (1.5pt);
\draw[below] (4.5,0) node {$x_{j+1}^n$};

\fill (canvas cs:x=2cm,y=3cm) circle (1.5pt);
\draw[above] (2,3) node {$x_{j-1}^{n+1}$};
\fill (canvas cs:x=4.2cm,y=3cm) circle (1.5pt);
\draw[above] (4.2,3) node {$x_j^{n+1}$};
\fill (canvas cs:x=6.3cm,y=3cm) circle (1.5pt);
\draw[above] (6.3,3) node {$x_{j+1}^{n+1}$};
\end{tikzpicture}
\caption{An illustration of the movement of element boundaries.}
\label{f-1}
\end{figure}
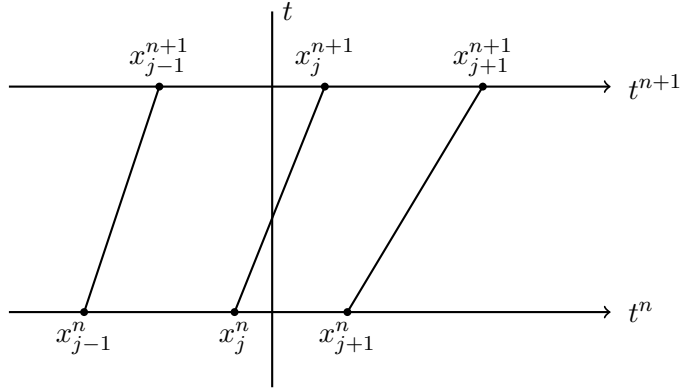

\begin{rem}
\label{remark:mixed-der-1}
On a moving mesh, a $C^0$ finite element method cannot be applied
to the original equation (\ref{RLW-1}) directly. Indeed, the weak formulation of (\ref{RLW-1}) takes the form
\[
\int_\Omega \left ( (u_t + \alpha u_x + \beta u_y + \gamma u u _x + \delta u u_y) \phi + \mu u_{xt} \phi_x
+ \mu u_{yt} \phi_y \right ) d x d y = 0, \quad \forall \phi \in H_0^1(\Omega) .
\]
Then a finite element approximation will contain space-time mixed derivatives
\begin{equation}
\frac{\partial^2 u_h}{\partial x \partial t}, \quad \frac{\partial^2 u_h}{\partial y \partial t},
\label{mixed-der}
\end{equation}
where $u_h$ is a finite element approximation to $u$. Notice that $u_h$ is piecewise continuous
and $\nabla u_h$ is discontinuous across element boundaries. Since these boundaries
vary with time for a moving mesh, $\nabla u_h$ has jumps in the time direction for spatial points
where the element boundaries sweep through (see Fig.~\ref{f-1}) and cannot be differentiated
with respect to time at these points (even in weak sense). Thus, the terms
in (\ref{mixed-der}) are not defined on $\Omega$, and a moving mesh finite element method does not apply to
(\ref{RLW-1}) directly.
\qed
\end{rem}

\begin{rem}
Several other choices of new variables have been used in the numerical solution
of the RLW equation. For example, Luo and Liu \cite{Luo-1999}
(also see Gu and Chen \cite{Gu-2008} for a least squares mixed FEM)
use the new variable $ p = a u^2/2 - \delta u_{x t}$ for a mixed finite element approximation of the 1D RLW equation
\[
u_t + a u u_x - \delta u_{xxt} = 0,
\]
subject to a homogeneous Dirichlet boundary condition. 
They use the weak formulation
\[
\begin{cases}
\int_{\Omega} (u_t \phi - p \phi_x) d x = 0,\quad \forall \phi \in H_0^1(\Omega)
\\
\int_{\Omega} \left (p - \frac{a}{2} u^2 + \delta u_{x t}\right ) \psi d x  = 0, \quad \forall \psi \in L^2(\Omega) .
\end{cases}
\]
It does not work with a moving mesh finite element method since it contains space-time mixed derivatives.
More recently, Gao and Mei \cite{Gao2015} define $p = u_x$ for
the 1D RLW equation
\[
u_t + u_x + 6 u^2 u_x - \mu u_{xxt} = 0
\]
with a homogeneous Dirichlet boundary condition. They use the weak formulation
\[
\begin{cases}
\int_{\Omega} (p - u_x ) \phi_x d x = 0,\quad \forall \phi \in H_0^1(\Omega)
\\
\int_{\Omega} \left ( p_t \psi - p \psi_x - 6 u^2 p \psi_x + \mu p_{xt} \psi_x  \right ) d x  = 0,
\quad \forall \psi \in H^1(\Omega)
\end{cases}
\]
which once again does not work with a moving mesh finite element method since it contains
a space-time mixed derivative.
\qed
\end{rem}

\subsection{Conservation laws}

Olver \cite{Olver1979} shows that the RLW equation possesses three non-trivial
independent conservation laws. Each of such laws corresponds to an invariant quantity
if the solution vanishes on the boundary (i.e., $g \equiv 0$).
The first two for (\ref{RLW-1}) are  
\begin{equation}
E_1(t) = \int_\Omega u d x d y, \quad E_2(t) = \int_{\Omega}  (u^2 + \mu u_x^2 + \mu u_y^2 ) d x d y,
\label{invariant-1}
\end{equation}
which can readily be verified by multiplying (\ref{RLW-1}) with $1$ and $u$, respectively,
integrating the resulting equation over $\Omega$, and performing integration by parts.

We first consider if these quantities are conserved by the finite element approximation on a fixed mesh.
For this case, $\dot{\V{X}} \equiv 0$ and both $A$ and $M$ are independent of time.
Summing the rows of (\ref{fem-3}) and using (\ref{fem-2b}) and $u_h = \sum_{j=1}^{N_{vi}} u_j \phi_j$,
we have
\[
\sum_{i=1}^{N_{vi}} \int_\Omega \left [ \phi_i \frac{\partial u_h}{\partial t} + \nabla \phi_i \cdot \nabla
\frac{\partial u_h}{\partial t} + \left ( \alpha \frac{\partial u_h}{\partial x}
+ \beta \frac{\partial u_h}{\partial y} + \gamma u_h \frac{\partial u_h}{\partial x}
+ \delta u_h  \frac{\partial u_h}{\partial y} \right )\phi_i \right ] d x d y = 0.
\]
This can be rewritten as
\begin{align*}
& \sum_{i=1}^{N_{v}} \int_\Omega \left [ \phi_i \frac{\partial u_h}{\partial t} + \nabla \phi_i \cdot \nabla
\frac{\partial u_h}{\partial t} + \left ( \alpha \frac{\partial u_h}{\partial x}
+ \beta \frac{\partial u_h}{\partial y} + \gamma u_h \frac{\partial u_h}{\partial x}
+ \delta u_h  \frac{\partial u_h}{\partial y} \right )\phi_i \right ] d x d y
\\
&\quad = \sum_{i=N_{vi}+1}^{N_{v}} \int_\Omega \left [ \phi_i \frac{\partial u_h}{\partial t} + \nabla \phi_i \cdot \nabla
\frac{\partial u_h}{\partial t} + \left ( \alpha \frac{\partial u_h}{\partial x}
+ \beta \frac{\partial u_h}{\partial y} + \gamma u_h \frac{\partial u_h}{\partial x}
+ \delta u_h  \frac{\partial u_h}{\partial y} \right )\phi_i \right ] d x d y .
\end{align*}
Noticing that $\sum_{i=1}^{N_{v}} \phi_i \equiv 1$ and using the divergence theorem
and the fact that $u_h|_{\partial \Omega} = 0$, we obtain
\begin{align}
& \frac{d }{d t} \int_\Omega u_h d x d y
\notag
\\
&\quad =  \sum_{i=N_{vi}+1}^{N_{v}} \int_\Omega \left [ \phi_i \frac{\partial u_h}{\partial t} + \nabla \phi_i \cdot \nabla
\frac{\partial u_h}{\partial t} + \left ( \alpha \frac{\partial u_h}{\partial x}
+ \beta \frac{\partial u_h}{\partial y} + \gamma u_h \frac{\partial u_h}{\partial x}
+ \delta u_h  \frac{\partial u_h}{\partial y} \right )\phi_i \right ] d x d y .
\label{invariant-2}
\end{align}
Thus, $E_1$ is not conserved by (\ref{fem-3}) since the right-hand side does not vanish in general.
An estimate of the derivation can be obtained as follows. Noticing that $E_1(t) = E_1(0)$ and
using Schwarz's inequality, we have
\begin{align*}
& \left | \int_\Omega u_h(x, y, t) d x  d y - \int_\Omega u_h(x, y, 0) d x  d y\right |
\\
\qquad 
& \le \int_\Omega | u_h(x, y, t) - u(x, y, t)| d x  d y + \int_\Omega | u_h(x, y, 0) - u(x, y, 0)| d x  d y
\\
\qquad & \le \left ( \| e_h(\cdot, t)\|_{L^2(\Omega)} + \| e_h(\cdot, 0)\|_{L^2(\Omega)}\right ) |\Omega|^{\frac 1 2}  .
\end{align*}
Assuming that the finite element error is second order in $L^2$ norm, we have
\begin{equation}
\Delta E_1(t) \equiv \int_\Omega u_h(x, y, t) d x  d y - \int_\Omega u_h(x, y, 0) d x  d y = \mathcal{O} (h^2),
\label{invariant-4}
\end{equation} 
where $h$ is the maximal diameter of the elements. It is interesting to point out that the numerical examples
in section~\ref{SEC:numerics} show that the difference decreases much faster than what shown in (\ref{invariant-4})
as $N \to \infty$. This may be attributed to the cancellation between terms on the right-hand
side of (\ref{invariant-2}) and the fact that $u_h$ and its derivatives are getting smaller on the boundary elements
which are getting closer to the boundary as $N$ increases.

Similarly, multiplying the $i$-th row of (\ref{fem-3}) with $u_i$ and summing all of the resulting rows
we can get
\begin{equation}
\frac{d }{d t} \int_\Omega \left [ u_h^2 + \mu \left (\frac{\partial u_h}{\partial x}\right )^2
+ \mu \left (\frac{\partial u_h}{\partial y}\right )^2 \right ] d x d y = 0,
\label{invariant-3}
\end{equation}
which implies that $E_2$ is conserved by (\ref{fem-3}). It is noted that this conservation holds only for
the semi-discrete scheme (\ref{fem-3}). It may not necessarily hold for the fully discrete scheme.
Nevertheless, (\ref{invariant-3}) implies that the difference will be small when the time step is small.

We now consider the moving mesh situation. Generally speaking, $\dot{\V{X}} \not \equiv 0$ and both $A$ and $M$
are time dependent for this case. In principle, we can perform a similar analysis as for the fixed mesh case.
Since the derivation is very tedious and the results are not that useful, we choose to not give the analysis
here. Instead, we simply state that the finite element method with a moving mesh does not conserve
either quantity. This will be verified by the numerical examples.
Moreover, assuming that the finite element error is second order in $L^2$ norm and first order in semi-$H^1$
norm, we can show that the FE approximation on a moving mesh possesses the property (\ref{invariant-4}) and
\begin{align}
 \Delta E_2(t) \equiv & \int_\Omega \left [ u_h^2 + \mu \left (\frac{\partial u_h}{\partial x}\right )^2
+ \mu \left (\frac{\partial u_h}{\partial y}\right )^2 \right ](x,y,t) d x  d y 
\notag
\\
& \qquad - \int_\Omega \left [ u_h^2 + \mu \left (\frac{\partial u_h}{\partial x}\right )^2
+ \mu \left (\frac{\partial u_h}{\partial y}\right )^2 \right ](x,y,0) d x  d y = \mathcal{O} (h) .
\label{invariant-5}
\end{align} 
Moreover, the numerical examples show that both $\Delta E_1(t)$ and $\Delta E_2(t)$ decreases
much faster than what indicated in (\ref{invariant-4}) and (\ref{invariant-5}).
Particularly, $\Delta E_1(t)$ behaves similarly for both fixed and moving meshes.

\section{An MMPDE-based moving meshes strategy}
\label{SEC:mmpde}

In this section we describe the generation of $\Th^{n+1}$ based on $\Th^n$ and $\V{u}^n$
using an MMPDE-based moving mesh strategy. The strategy uses a metric tensor (a symmetric
and uniformly positive definite matrix-valued function) to specify the information of the size
shape, and orientation of the elements throughout the domain. We take a Hessian-based metric tensor as
\begin{equation}
\M = \det( \alpha_h I + |H(u_h^n)| )^{-\frac{1}{6}} ( \alpha_h I + |H(u_h^n)| ),
\label{M-1}
\end{equation}
where $I$ is the identity matrix, $\det(\cdot)$ denotes the determinant of a matrix,
$H(u_h^n)$ is a recovered Hessian from the finite element
solution $u_h^n$, $|H(u_h^n)| = Q \text{diag}(|\lambda_1|, |\lambda_2|) Q^T$
with $Q \text{diag}(\lambda_1, \lambda_2) Q^T$ being the eigen-decomposition
of $H(u_h^n)$, and $\alpha_h$ is a regularization parameter defined through the equation
\[
\sum_{K \in \Th} |K| \det (\M)^{\frac{1}{2}} \equiv \sum_{K \in \Th} |K| \det(\alpha_h I + |H(u_h^n)|)^{\frac 2 3}
= 2 \sum_{K \in \Th} |K| \det(|H(u_h^n)|)^{\frac 2 3} .
\]
It is noted that the above equation equation uniquely defines $\alpha_h$ and can be solved using, for instance,
the bisection method. Moreover, 
the metric tensor (\ref{M-1}) is optimal for the $L^2$ norm of linear interpolation
error \cite{HS03}. Furthermore, in our computation $H(u_h^n)$ at any vertex is recovered by differentiating
a quadratic polynomial that fits the values of $u_h^n$ at the neighboring vertices in the least square sense
(e.g., see \cite{KaHu2013}).

A key of the MMPDE-based moving mesh strategy is to view any nonuniform mesh
as a uniform one in the metric $\M$. To explain this, we consider a physical mesh $\Th$
and a computational mesh $\Tc$, either of which can be viewed as a deformation of the other.
Then, $T_h$ is said to be an $\M$-uniform mesh in the metric $\M$ (e.g., see \cite{Hua06,HR11})
if it satisfies
\begin{align}
|K| \det(\M_K)^{\frac{1}{2}} = \frac{|K_c| \sigma_h}{|\Omega_c|}, \quad \forall K \in \Th
\label {equidistribution}
\\
\frac{1}{2} \tr\left ( (F'_K)^{-1}\M_K^{-1}(F'_K)^{-T} \right ) = 
\det \left ( (F'_K)^{-1}\M_K^{-1}(F'_K)^{-T} \right )^ {\frac{1}{2}}, \quad \forall K \in \Th
\label {alignment}
\end{align}
where $K$ is an element of $\Th$, $K_c$ is the element of $\Tc$ corresponding to $K$,
$|K|$ and $|K_c|$ denote the volumes of $K$ and $K_c$, respectively, $|\Omega_c| = \sum_{K_c \in \Tc} |K_c|$,
 $\sigma_h = \sum_{K\in \Th} |K| \det(\M_K)^{\frac{1}{2}}$,
$F'_K$ is the Jacobian matrix of the affine mapping $F_K: K_c \rightarrow K$,
$\M_K$ is the average of $\M$ over $K$, and $\tr(\cdot)$ denotes the trace of a matrix.
The condition (\ref{equidistribution}) is referred to as the equidistribution condition which determines 
the size of elements through the metric tensor $\M$. The bigger $\det(\M_K)^{\frac{1}{2}}$ is, the smaller
the element $K$ is. On the other hand, (\ref{alignment}) is called the alignment condition, which
requires $K$, when measured in the metric $\M_K$, to be similar to $K_c$ and in this way, determines
the shape and orientation of $K$ though $\M_K$ and $K_c$.

The meshing strategy we use is to generate a mesh satisfying (\ref{equidistribution}) and (\ref{alignment})
as closely as possible. This is done by minimizing the energy
\begin{align}
I_h(\Th, \Tc) & = \frac{1}{3} \sum_K |K| \det(\M_K)^{\frac{1}{2}} (  \tr ( (F'_K)^{-1}\M_K^{-1}(F'_K)^{-T} ) )^{2}
+ \frac{4}{3} \sum_K |K| \det(\M_K)^{-\frac{1}{2}}  \det(F'_K)^{-2},
\label{energy}
\end{align}
which is a Riemann sum of a continuous functional developed in \cite{Hua01b} based on
equidistribution and alignment for variational mesh generation and adaptation.
Instead of minimizing $I_h(\Th, \Tc)$ directly, we define the mesh equation as a gradient system of $I_h(\Th, \Tc)$
(the MMPDE approach). For example, assume that we have chosen a quasi-uniform reference computational
mesh $\hTc$. Then $I_h(\Th, \hTc)$ is a function of $\Th$ or the coordinates of its vertices, $\V{x}_i$, $i =1, ..., N_v$.
The mesh equation is
\begin{equation}
\frac{d \V{x}_i}{d t} = -\frac{\det( \M (\V{x}_i)  )^{\frac{1}{2}}}{\tau}
\left ( \frac{\partial I_h}{\partial  \V{x}_i} \right )^T,\quad i = 1,...,N_v
\label{mmpde-x}
\end{equation}
where ${\partial I_h}/{\partial  \V{x}_i}$ is considered as a row vector and $\tau$ is a parameter used
for adjusting the time scale for the mesh movement to respond the changes in $\M$.
This $\V{x}$-formulation of the mesh equation, under suitable modifications for the boundary vertices
(to keep them on the boundary),
can be integrated from $t_n$ to $t_{n+1}$ (starting with $\Th^n$) to obtain the new mesh $\Th^{n+1}$.
Moreover, it is shown in \cite{HK2015} that $\Th^{n+1}$ is nonsingular and its minimal volume
and minimal height of the elements have positive lower bounds that depend only on the number of elements,
the initial mesh, and the metric tensor.

A major disadvantage of the above $\V{x}$-formulation is that we need to consider the dependence
of $\M$ on $\V{x}$ when computing the derivatives ${\partial I_h}/{\partial  \V{x}_i}$. The metric tensor $\M$
needs to be updated constantly (through interpolation) during the integration of (\ref{mmpde-x}) since
$\M$ is typically available only at the vertices of $\Th^n$. To avoid this difficulty, we use
the so-called $\V{\xi}$-formulation where we take $\Th = \Th^n$ and consider $I_h(\Th^n, \Tc)$ as a function
of the coordinates of the computational vertices, $\V{\xi}_i$, $i =1, ..., N_v$. The mesh equation is defined as
\begin{equation}
\frac{d \V{\xi}_i}{d t} = -\frac{\det( \M (\V{x}_i)  )^{\frac{1}{2}}}{\tau}
\left ( \frac{\partial I_h}{\partial  \V{\xi}_i} \right )^T,\quad i = 1,...,N_v .
\label{mmpde-xi}
\end{equation}
This equation, under suitable modifications for the boundary vertices (to keep them on the boundary),
can be integrated from $t_n$ to $t_{n+1}$ (starting with $\hTc$) to obtain the new mesh $\Tc^{n+1}$.
Note that $\Th^n$ is kept fixed and there is no need to update $\M$ during the integration.
Denote the correspondence between $\Tc^{n+1}$ and $\Th^n$ by $\Phi_h$, i.e.,
$\Th^n = \Phi_h(\Tc^{n+1})$. The new physical mesh is defined as $\Th^{n+1} = \Phi_h(\hTc)$,
which can be computed using linear interpolation.

Numerical experiment has shown that both $\V{x}$- and $\V{\xi}$-formulations are effective
in producing adaptive meshes. However, the latter will lead to simpler formulas since
there is no need to consider the dependence on $\M$ when calculating 
${\partial I_h}/{\partial  \V{\xi}_i}$. Using the notion of scalar-by-matrix differentiation,
we can find the analytical expressions for these derivatives; the interested reader is referred to
\cite{HK2014} for the detailed derivation. With those formulas, we can rewrite (\ref{mmpde-xi}) into
\begin{equation}
\frac{d \V{\xi}_i} {dt} = \frac{\det( \M (\V{x}_i)  )^{\frac{1}{2}}}{\tau}  \sum_{K \in \omega _i} |K| \V{v}_{i_K}^K,
\label{MMPDE-xi-2}
\end{equation}
where $\omega_i$ is the element patch associated with the vertex $\V{x}_i$, 
$i_K$ is the local index of $\V{x}_i$ in $K$ and the local velocities $\V{v}_{i_K}^K$ are given by
\begin{align}
\begin{bmatrix} (\V{v}_1^K)^T \\ (\V{v}_2^K)^T \end{bmatrix} =
 -E_K^{-1} \frac{\partial G}{\partial \J } - 
 \frac{\partial G}{\partial \det(\J )} \frac{ \det (E_{K_c})} {\det(E_K)} E_{K_c}^{-1},
 \quad 
\V{v}_0^K = - \V{v}_1^K- \V{v}_2^K .
\label{local-v}
\end{align}
Here, $E_K = [\V{x}_1^K - \V{x}_0^K, \V{x}_2^K - \V{x}_0^K]$  and
$E_{K_c} = [\V{\xi}_1^K - \V{\xi}_0^K, \V{\xi}_2^K - \V{\xi}_0^K]$ are the edge matrices of $K$
and $K_c$, respectively,
and the function $G = G(\J, \det(\J))$ (with $\J = (F_K')^{-1} = E_{K_c} E_K^{-1}$)
is associated with the energy (\ref{energy}).
It and its derivatives are given by
\begin{align*}
& G(\J, \det(\J)) =  \frac{1}{3} \det (\M)^{\frac{1}{2}} ( \tr (\J \M_K^{-1}) \J^T )^{2} 
+ \frac{4}{3}  \det (\M_K)^{-\frac{1}{2}} \det(\J)^2,
\\
&\frac{\partial G}{\partial \J} =  \frac{4}{3} \det(\M)^{\frac{1}{2}} \tr(\J \M_K^{-1} \J^T ) \M_K^{-1} \J^T,
\\
&\frac{\partial G}{\partial \det(\J)} = \frac{8}{3}  \det(\M_K)^{-\frac{1}{2}} \det ( \J ).
\end{align*}
In actual computation, the edge matrices and local velocities are first computed for all elements.
Then the nodal mesh velocities are assembled according to (\ref{MMPDE-xi-2}).

\section{Numerical results}
\label{SEC:numerics}

In this section we present numerical results obtained with the moving mesh finite element method
described in the previous sections for a number of 1D and 2D examples for the RLW and MRLW equations.
We shall demonstrate the second order convergence of the method in space and its ability to concentrate
mesh points in needed regions.
The error in the numerical solution is measured in the (global) $L^2$ and $L^\infty$ norm, i.e.,
\[
\int_0^T \| e^h(\cdot, t) \|_{L^2(\Omega)}d t, \quad \int_0^T \| e^h(\cdot, t) \|_{L^\infty(\Omega)}d t .
\]
The parameter $\tau$ for mesh movement is taken as $\tau = 10^{-4}$ for 1D examples and $\tau = 10^{-2}$
for 2D examples.

\begin{figure}[htb]
\centering
\hbox{
\begin{minipage}[b]{3in}
\centerline{(a): Computed solution}
\centering
\includegraphics[width=3in]{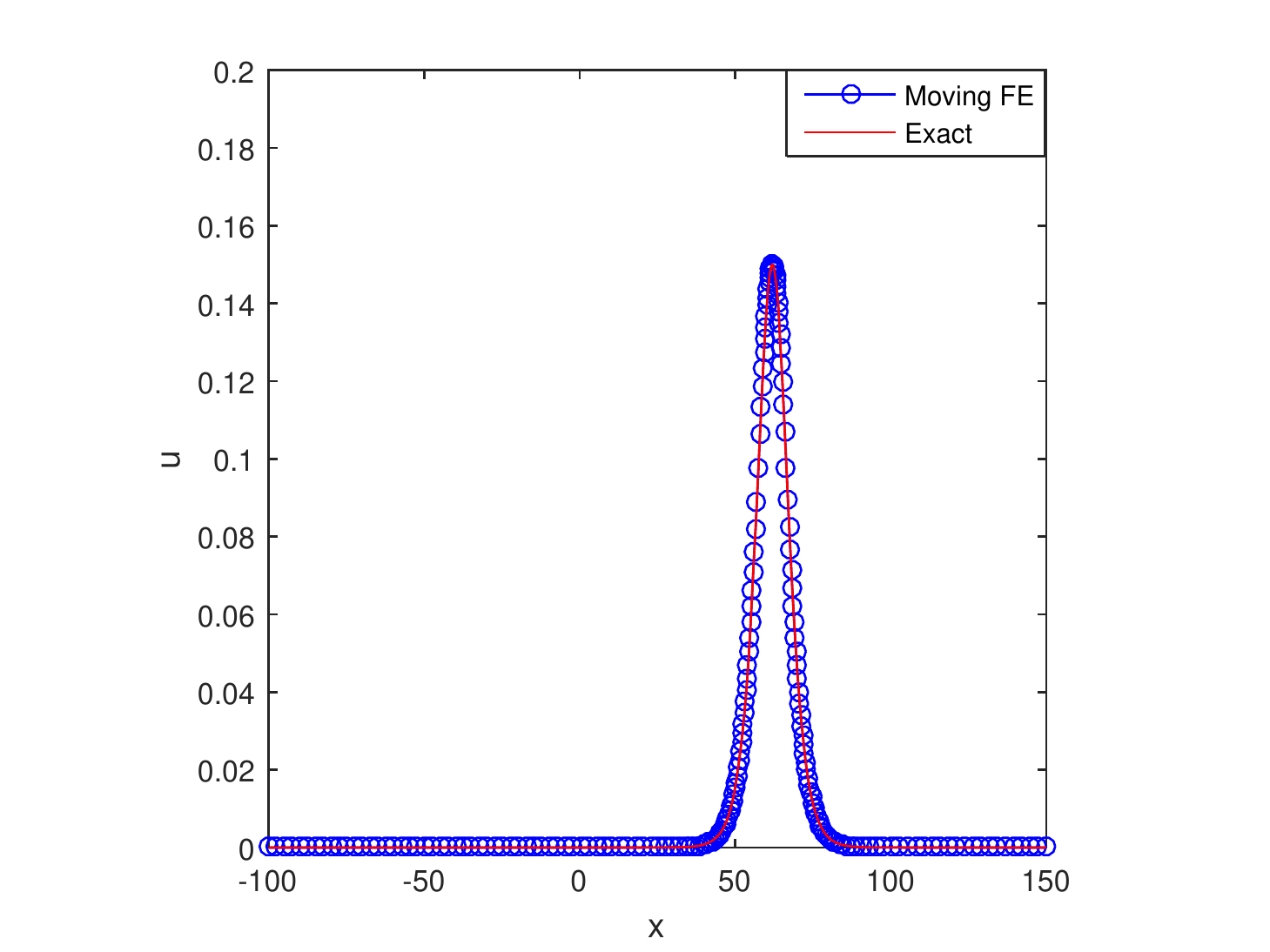}
\end{minipage}
\hspace{5mm}
\begin{minipage}[b]{3in}
\centerline{(b): Mesh trajectories}
\centering
\includegraphics[width=3in]{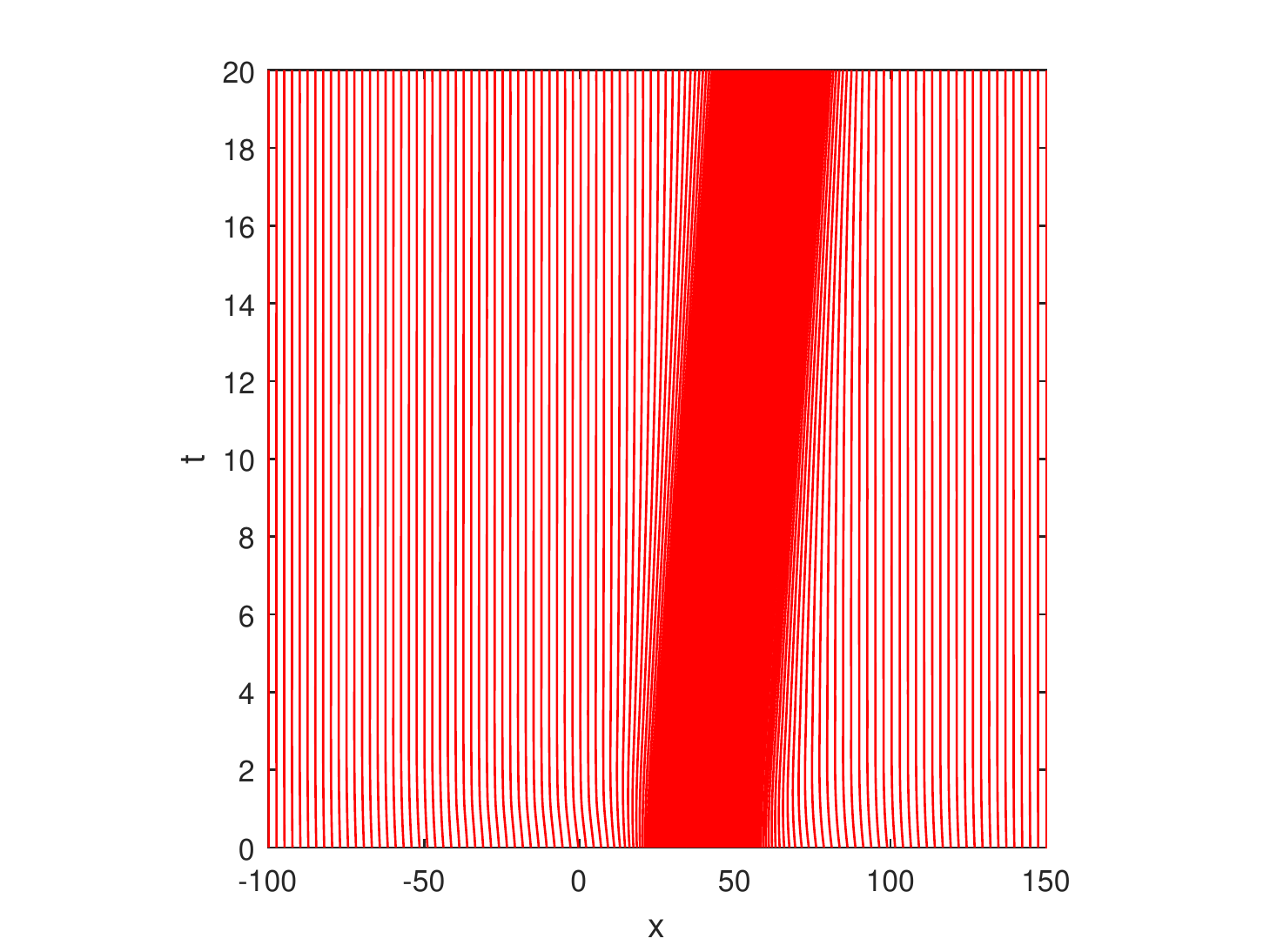}
\end{minipage}
}
\caption{Example~\ref{exam3.1}. The numerical solution and mesh trajectories are obtained with
the moving mesh finite element method ($N=200$) for the 1D RLW equation with a single soliton.}
\label{exam-2.1-1}
\end{figure}

\begin{figure}[thb]
\centering
\hbox{
\begin{minipage}[t]{3in}
\centerline{(a)}
\centering
\includegraphics[scale=0.4]{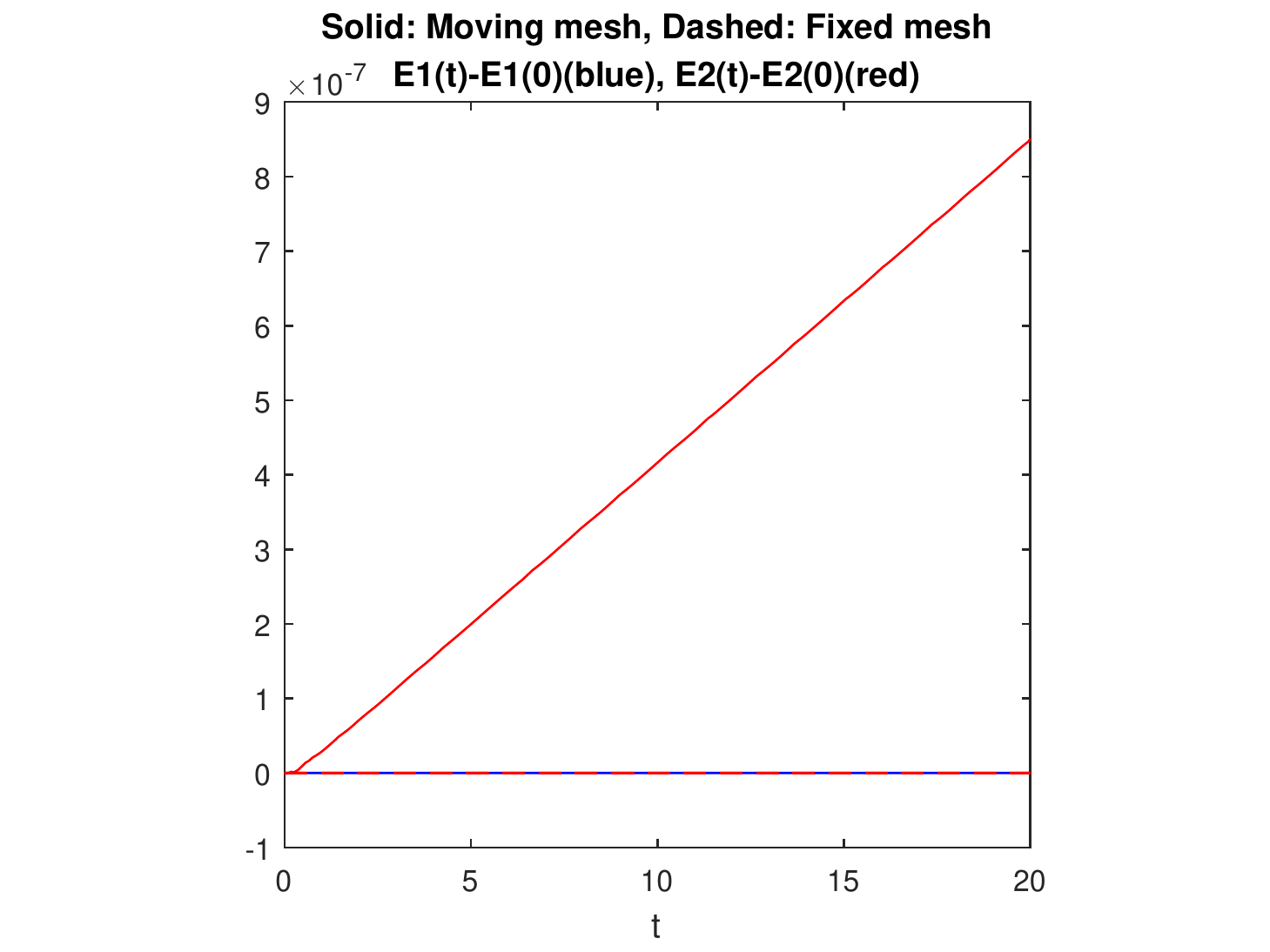}
\end{minipage}
\hspace{5mm}
\begin{minipage}[t]{3in}
\centerline{(b)}
\centering
\includegraphics[scale=0.4]{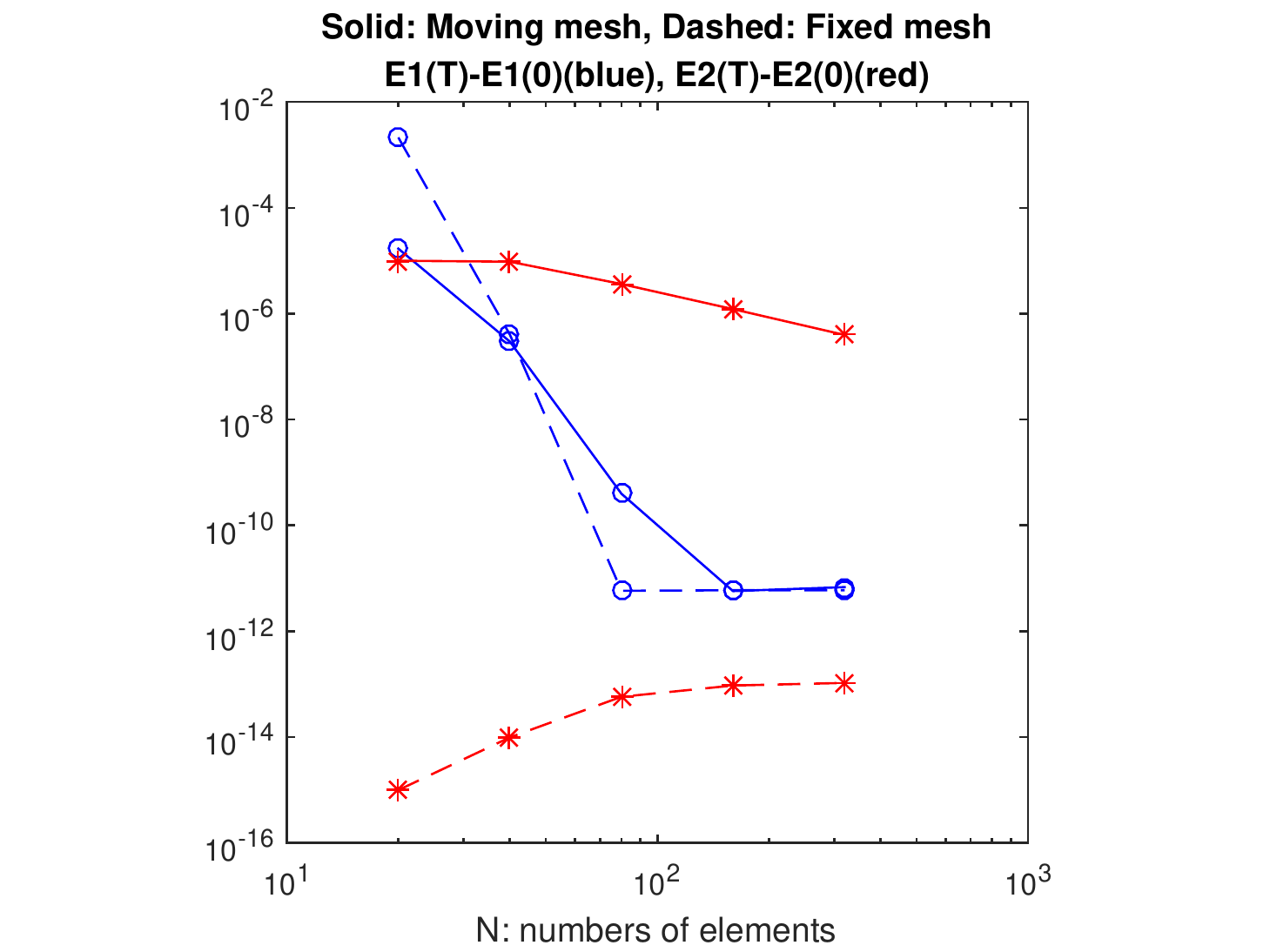}
\end{minipage}
\hspace{5mm}
}
\caption{Example~\ref{exam3.1}. The solid and dashed blue curves are
for $E_1(T)-E_1(0)$ with moving and fixed meshes, respectively,
while the solid and dashed red curves for $E_2(T)-E_2(0)$ with moving and fixed meshes, respectively.
(a) The difference of the conserved quantities for a mesh of $N=200$ is plotted as a function of time.
The solid blue, dashed blue, and dashed red curves are almost indistinguishable for this relativelys fine mesh.
(b) $E_1(T)-E_1(0)$ and $E_2(T)-E_2(0)$ are plotted as functions of $N$.}

\label{exam-2.1-2}
\end{figure}

\begin{exam}
\label{exam3.1} (1D RLW with a single soliton)
We consider the 1D RLW equation
\begin{equation}
u_t + u_x + \gamma u  u_x - \mu u_{xxt} = 0,
\label {RLW-1d}
\end{equation}
with $\gamma = 2$, $\mu = 1$, and $\Omega = (-100,150)$.
The Dirichlet boundary condition is chosen such that
the exact solution is a solitary wave
\[
u(x,t) = \frac{3 c}{2} \; \text{sech}^2\left (k (x- v t-x_0)\right ),
\]
where
$k = \frac{1}{2}\sqrt{\frac{v}{\mu( v+1)}}$, $v=c+1$, $x_0 = 40$, and $c = 0.1$.
The soliton has an amplitude $\frac{3c}{2}$ and a propagation velocity $v$.
A large spatial domain is chosen so that the solution is almost zero at the boundary and the example
can be used to check the conservation of $E_1$ and $E_2$. The computation is performed with $T = 20$.

The error and convergence order are listed in Table~\ref{exam-2.1-order} for both fixed and moving
meshes. It can be seen that while both types of mesh lead to the same second order of convergence, moving meshes
produce more accurate solutions (with the error being an order of magnitude smaller)
than fixed meshes. A typical numerical solution and the corresponding mesh trajectories are shown
in Fig.~\ref{exam-2.1-1}. It can be seen that the mesh points are concentrated in the peak area of the soliton
for the whole time, demonstrating the mesh adaptation ability of the method.

In Fig.~\ref{exam-2.1-2}(a), the difference of the conserved quantities is plotted as a function of $t$
for $N = 200$.  Notice that $\Delta E_1(t)$ for fixed and moving meshes (blue solid and dashed lines)
and $\Delta E_2(t)$ for the fixed mesh are indistinguishable. (In fact, they are almost at the level of roundoff error.)
The difference of the conserved quantities is plotted as a function of $N$ in Fig.~\ref{exam-2.1-2}(b).
We can see that $\Delta E_1(T)$ for both fixed and moving meshes
and $\Delta E_2(T)$ for moving meshes are quite significant for relatively small $N$.
However, $\Delta E_1(T)$ decreases quickly to the level of roundoff error as $N$ increases
for both fixed and moving meshes.
On the other hand, with fixed meshes $\Delta E_2(T)$ remains very small for the considered
range of $N$, consistent with the fact that $E_2$ is conserved on a fixed mesh by the semi-discrete
system of the method. With moving meshes, $\Delta E_2(T)$ is much bigger, reflecting the fact that
$E_2$ is not conserved by the method on moving meshes. Nevertheless, it decreases at a rate
$\mathcal{O}(N^{-1.6})$, much faster than the first order predicted in (\ref{invariant-5}).
Thus far we have seen that this example the fixed mesh method has better conservation properties
than the moving mesh method but gives less accurate solutions. It could be interesting to explore
what advantages the conservation of the quantities gives to the scheme for the RLW equation.

\end{exam}

\begin{table}[htb]
\begin{center}
\caption{Example~\ref{exam3.1}. $L^2$ and $L^\infty$ error and convergence order on moving and fixed meshes. }
\begin{tabular}{c|l|l|l|l|l|l|l|l}\hline \hline
                       & \multicolumn{4}{|c}{Moving Mesh}& \multicolumn{4}{|c}{Fixed Mesh}  \\ \hline
 $N$ & {$L^2$ error} &  order & $L^{\infty}$ error & order & {$L^2$ error} &  order & $L^{\infty}$ error & order\\ \hline\hline
20    &  2.62E-1  &             &   1.09E-1   &           &  4.86E-0  &             &   1.28E-0   & \\
40    &  5.44E-2  &   2.27   &   2.04E-2   &  2.42  &  1.57E-0  &   1.63   &   5.94E-1   &  1.11\\
80    &  1.29E-2  &   2.08   &   4.51E-3   &  2.17  &  3.34E-1  &   2.17   &   1.82E-1   &  1.71\\
160  &  3.15E-3  &   2.03   &   1.08E-3   &  2.07  &  7.86E-2  &   2.14   &   4.75E-2   &  1.94\\
320  &  7.84E-4  &   2.01   &   2.66E-4   &  2.02  &  1.91E-2  &   2.04   &   1.93E-2   &  1.99\\
640  &  1.96E-4  &   2.00   &   6.23E-5  &   2.00   &  4.76E-3  &   2.01   &   2.98E-3   &  2.00\\\hline\hline
\end{tabular}
\label{exam-2.1-order}
\end{center}
\end{table}

\begin{figure}[thb]
\centering
\hbox{
\begin{minipage}[t]{3in}
\centerline{(a) Computed solution}
\centering
\includegraphics[width=3.5in]{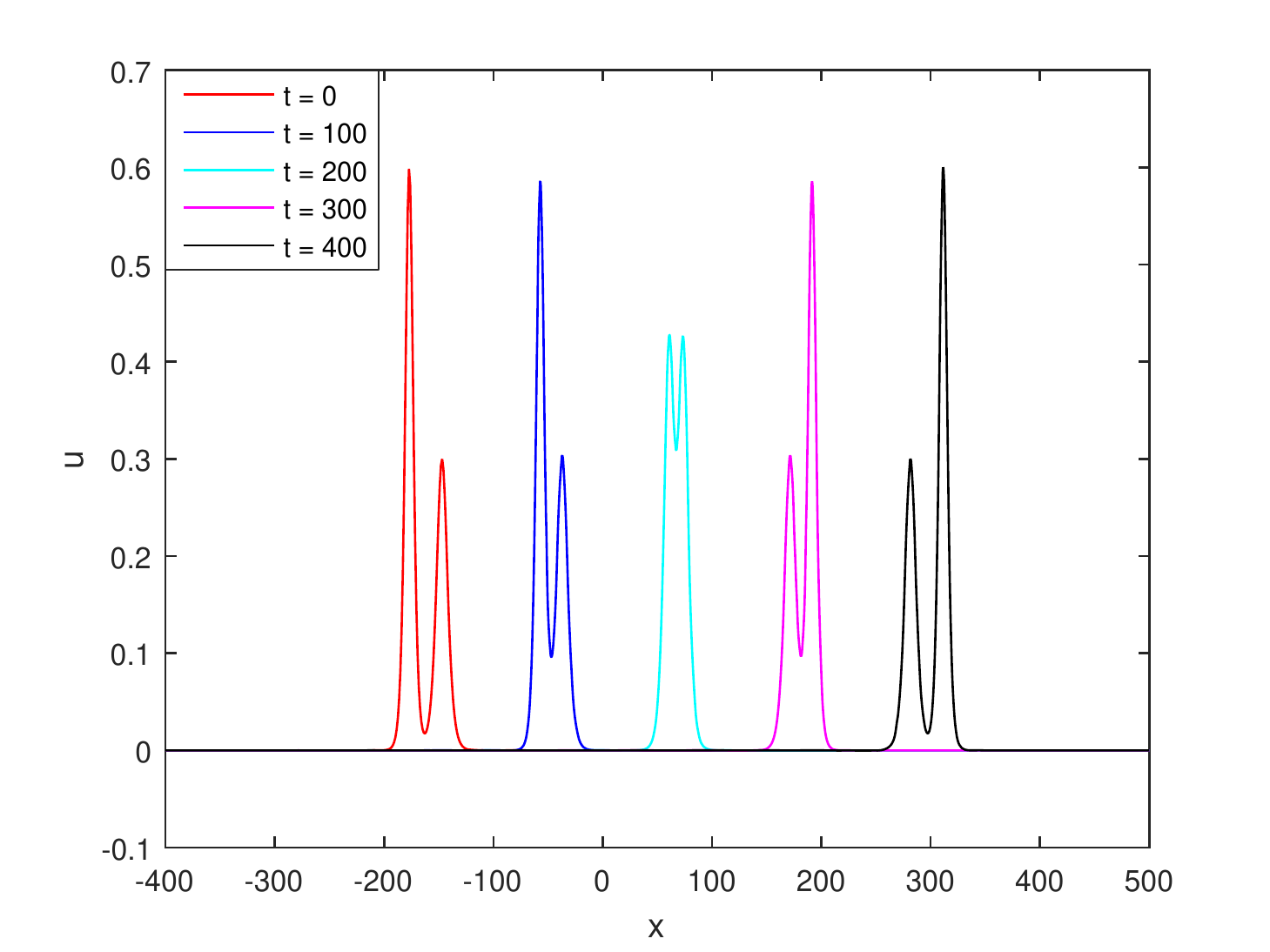}
\end{minipage}
\hspace{5mm}
\begin{minipage}[t]{3in}
\centerline{(b): Mesh trajectories}
\centering
\includegraphics[width=3.5in]{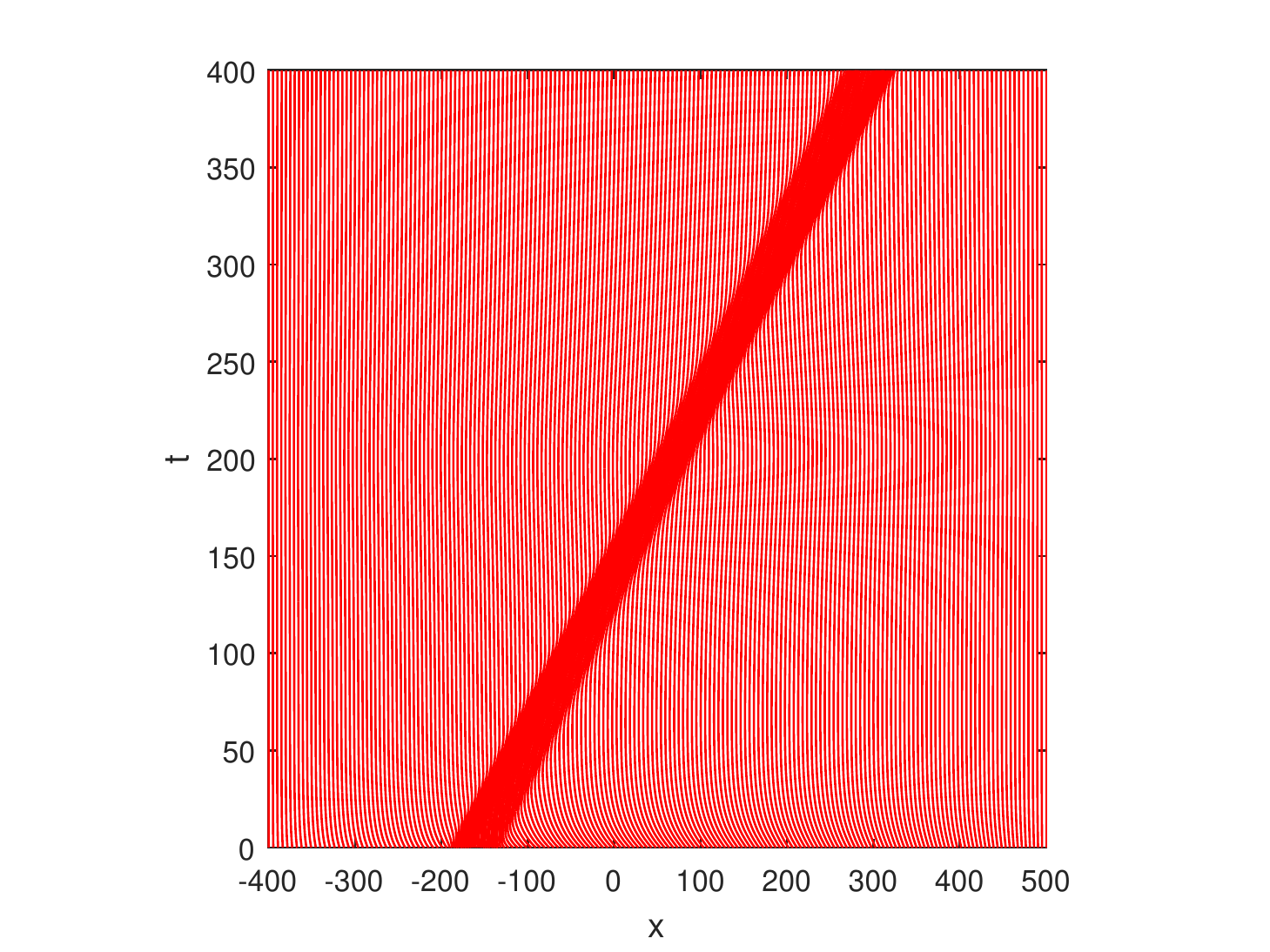}
\end{minipage}
\hspace{5mm}
}
\caption{Example~\ref{exam3.2}. A numerical solution at $t = 0,100,200,300,400$
and the mesh trajectories are obtained with the moving mesh finite element method with $N=800$.
As the value of $N$ is large, we only plot mesh trajectories every 4 nodes.
}
\label{exam-3.2-1}
\end{figure}

\begin{figure}[thb]
\centering
\hbox{
\begin{minipage}[t]{3in}
\centerline{(a) Computed solution with $N=800$}
\centering
\includegraphics[width=3.5in]{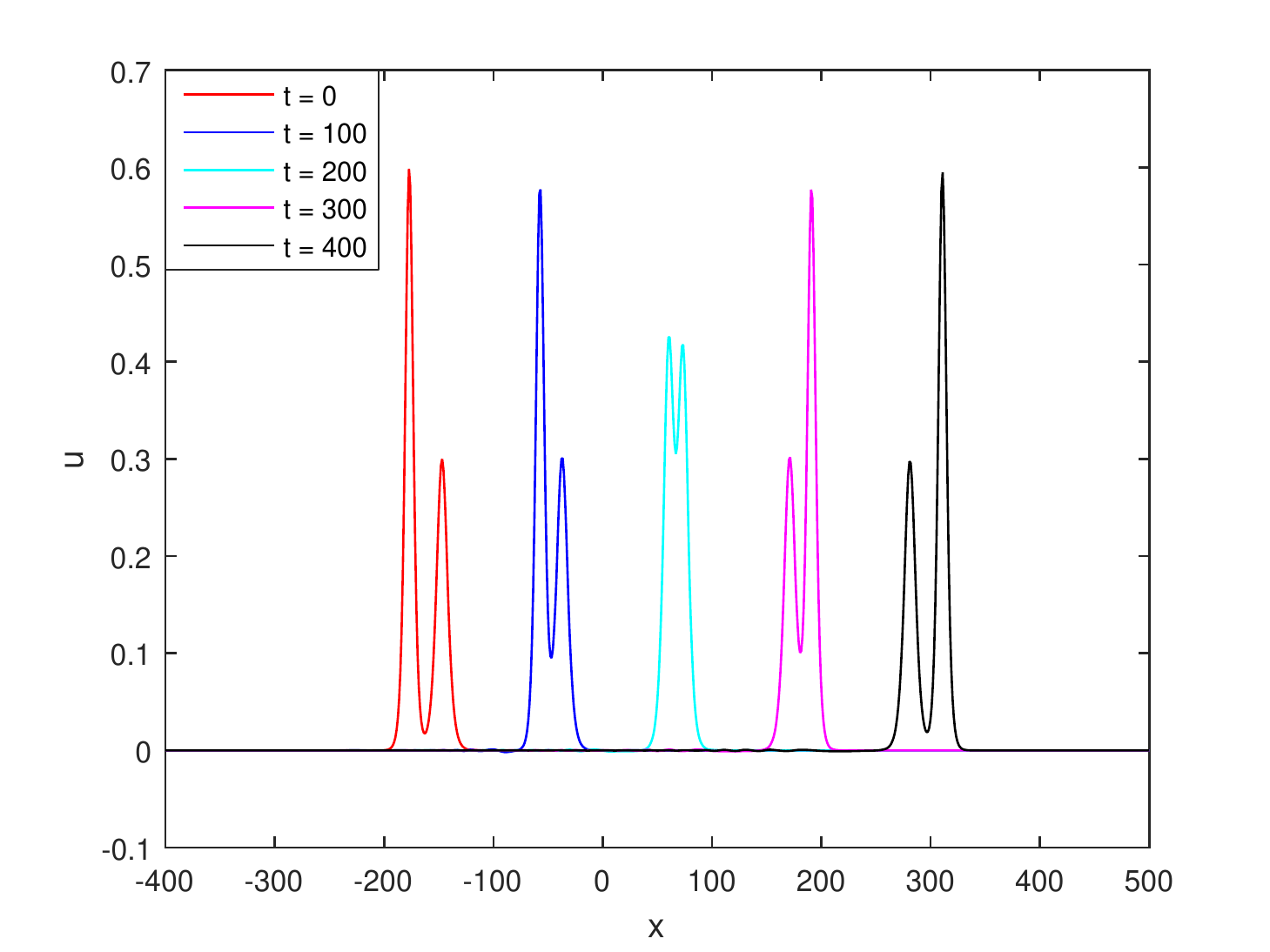}
\end{minipage}
\hspace{5mm}
\begin{minipage}[t]{3in}
\centerline{(b) Computed solution with $N=8000$}
\centering
\includegraphics[width=3.5in]{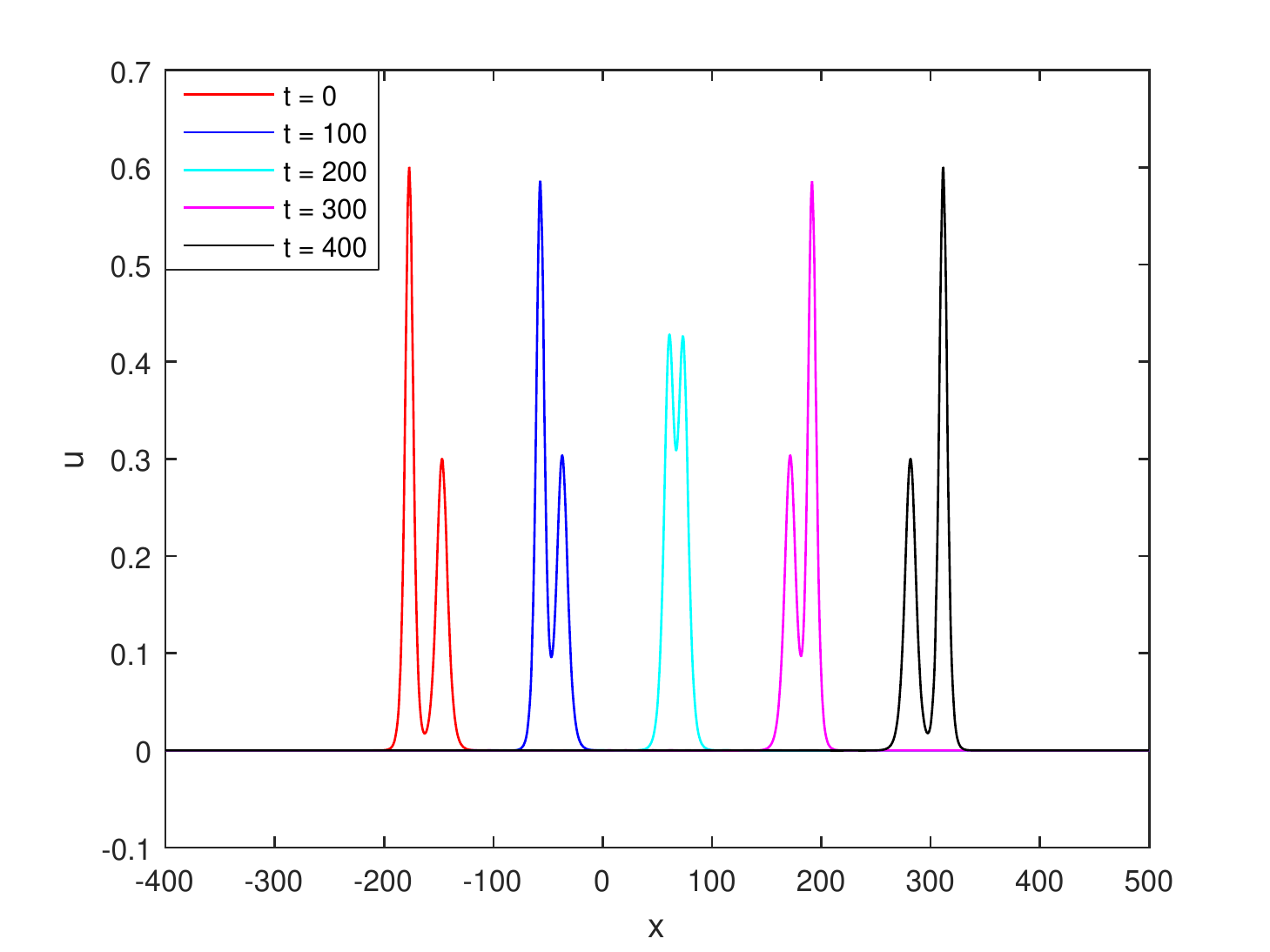}
\end{minipage}
\hspace{5mm}
}
\caption{Example~\ref{exam3.2}. Numerical solutions at $t = 0,100,200,300,400$
are obtained with fixed meshes of $N = 800$ and $8000$.
}
\label{exam-3.2-1F}
\end{figure}

\begin{figure}[thb]
\centering
\begin{minipage}[t]{3in}
\includegraphics[width=3.5in]{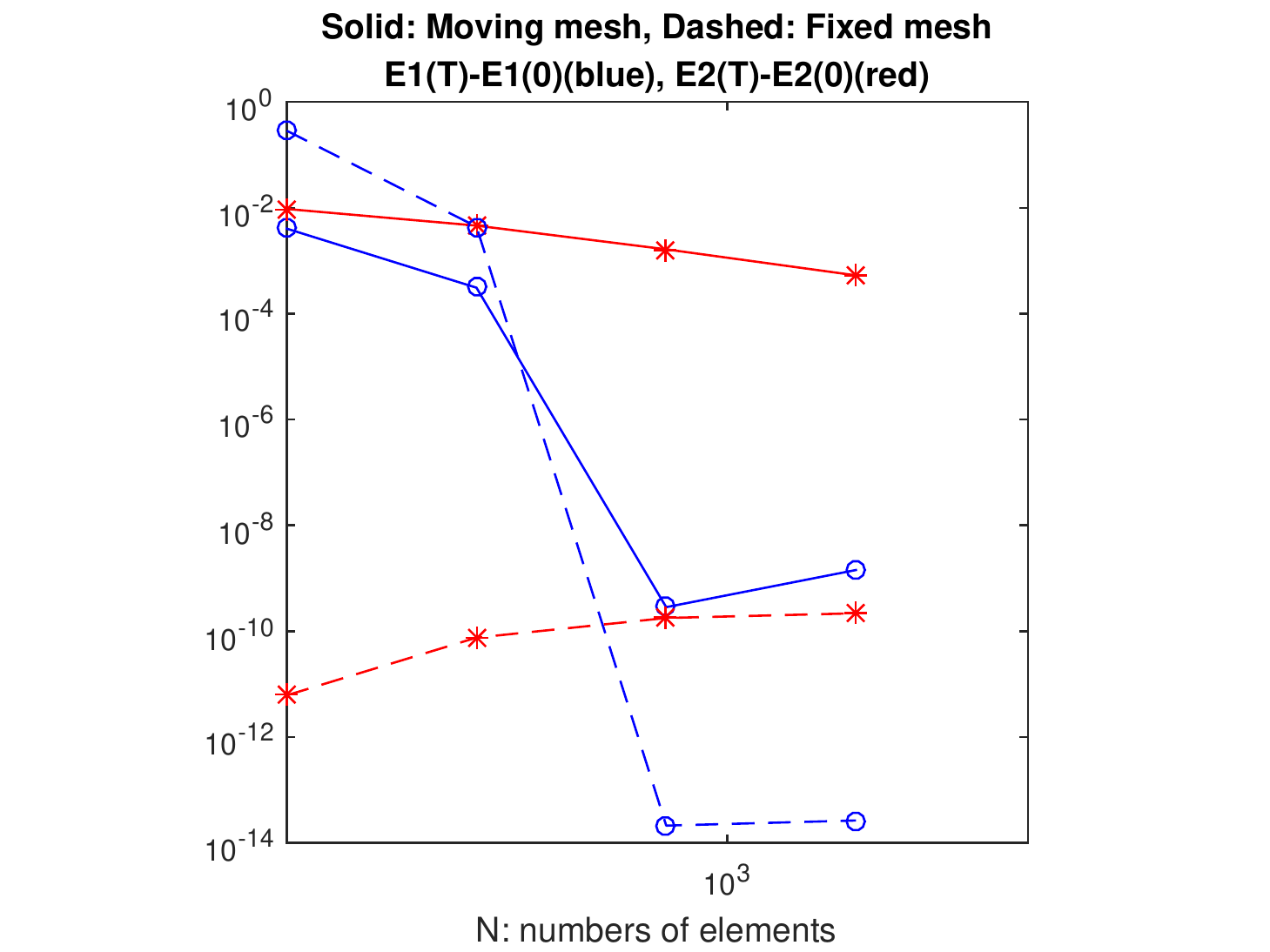}
\end{minipage}
\caption{Example~\ref{exam3.2}. 
The solid and dashed blue curves are
for $E_1(T)-E_1(0)$ with moving and fixed meshes, respectively,
while the solid and dashed red curves for $E_2(T)-E_2(0)$ with moving and fixed meshes, respectively.}
$E_1(T)-E_1(0)$ and $E_2(T)-E_2(0)$ are plotted as functions of $N$.
\label{exam-3.2-2}
\end{figure}

\begin{exam}
\label{exam3.2}
(1D RLW with interaction of two solitary waves)
In this example, we study the interaction of two solitary waves for 
the 1D RLW equation (\ref{RLW-1d})
with a homogeneous Dirichlet boundary condition and the initial condition
\[
u(x,0) = \sum_{j=1}^2 3 c_j \text{sech}^2\left (k_j (x-x_j)\right ),
\]
where $\gamma = \mu = 1$, $k_j = \frac{1}{2}\sqrt{\frac{\gamma v_j}{\mu(\gamma v_j+1)}}$, $v_j = 1+\gamma c_j$,
$x_1 = -177$, $x_2 = -147$, $c_1= 0.2$, and $c_2 = 0.1$.
Initially, the solitons have the amplitude $3 c_j$ and location $x_j$ ($j=1,2$)
and the larger soliton is placed on the left of the smaller one.
An interaction occurs as the larger one is catching up with and eventually passes
the smaller one. The simulation is performed on a domain $\Omega = (-400,500)$ until $t = 400$.
The exact analytical solution is unavailable for this example.

A numerical solution at $t = 0, 100, 200, 300, 400$ and the mesh trajectories are shown in Fig.~\ref{exam-3.2-1}.
The interaction of the two solitons can be clearly seen from the figure. Moreover, the width of the mesh concentration
also changes with time, becoming narrower during the interaction.
For comparison purpose, the solutions obtained with fixed meshes of $N=800$ and $8000$ are
plotted in Fig.~\ref{exam-3.2-1F}. Oscillations are visible along the $x$-axis
in the solution with the fixed mesh of $N=800$. The differences, $\Delta E_1(T)$ and $\Delta E_2(T)$,
are plotted as functions of $N$ in Fig.~\ref{exam-3.2-2}. 
Once again, $\Delta E_1(T)$ for fixed and moving meshes and $\Delta E_2(T)$ for moving meshes
are significant for small $N$.
$\Delta E_1(T)$ drops quickly
as $N$ increases for both fixed and moving meshes. On the other hand, $\Delta E_2(T)$ stays very small
for fixed meshes. It is relatively large for moving meshes although it decreases at a rate of
about $\mathcal{O}(N^{-1.6})$, which is faster than what indicated
by (\ref{invariant-5}).
\end{exam}

\begin{figure}[thb]
\centering
\hbox{
\begin{minipage}[t]{2in}
\centerline{(a) $d = 2$ (fixed mesh)}
\includegraphics[width=2.5in]{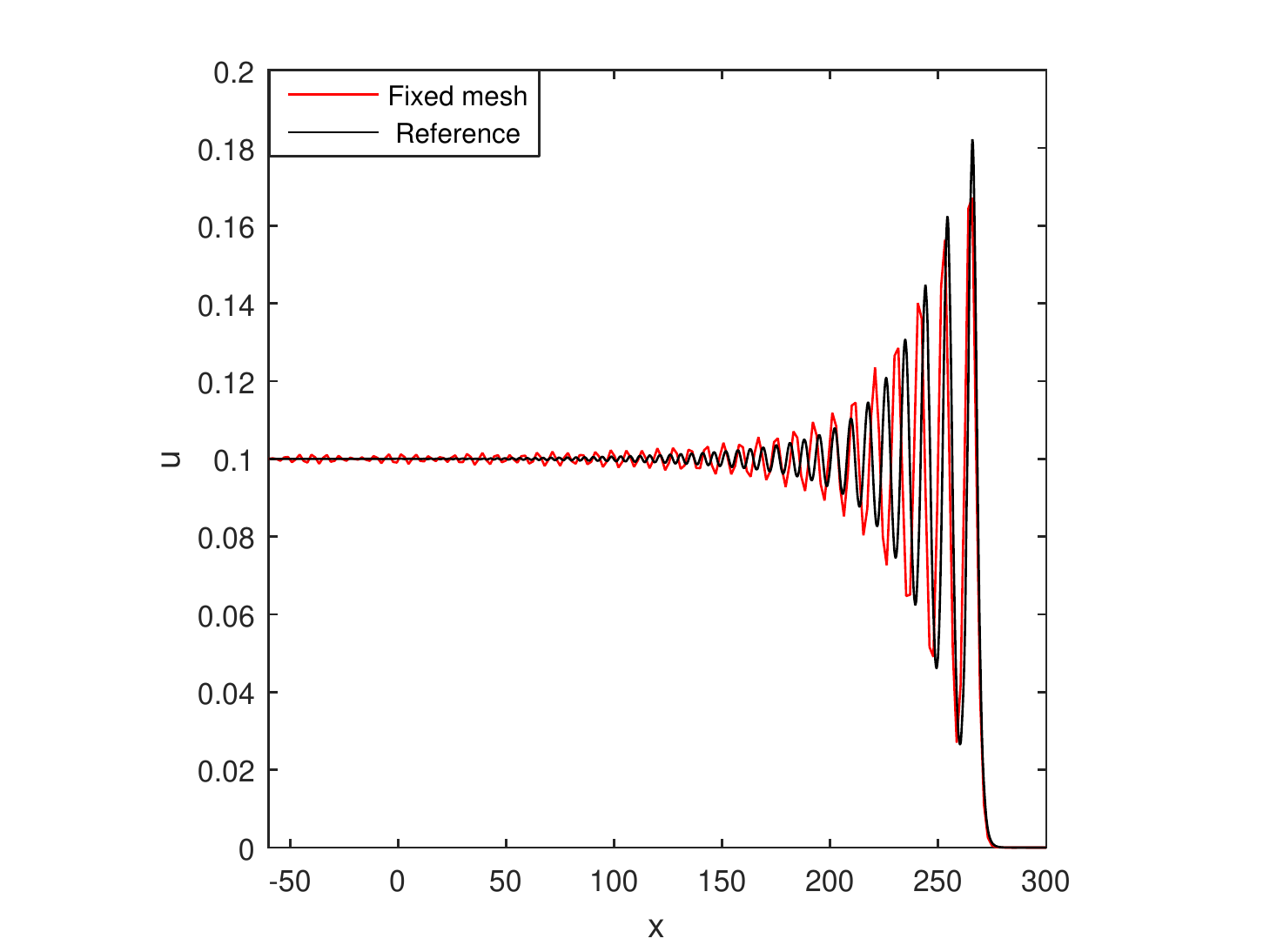}
\end{minipage}
\hspace{2mm}
\begin{minipage}[t]{2in}
\centerline{(b) $d = 2$ (moving mesh)}
\includegraphics[width=2.5in]{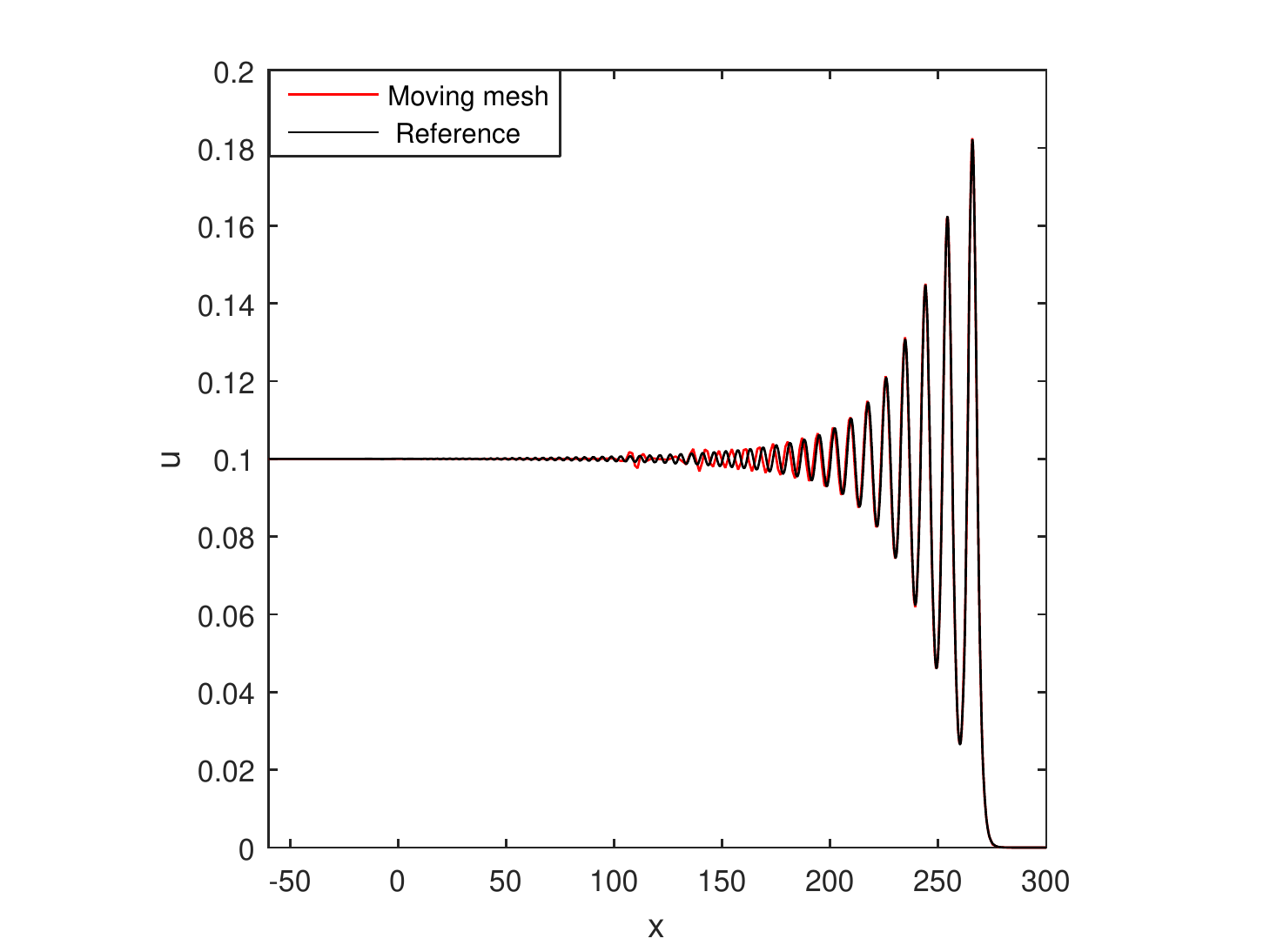}
\end{minipage}
\hspace{2mm}
\begin{minipage}[t]{2in}
\centerline{(c) $d = 2$}
\includegraphics[width=2.5in]{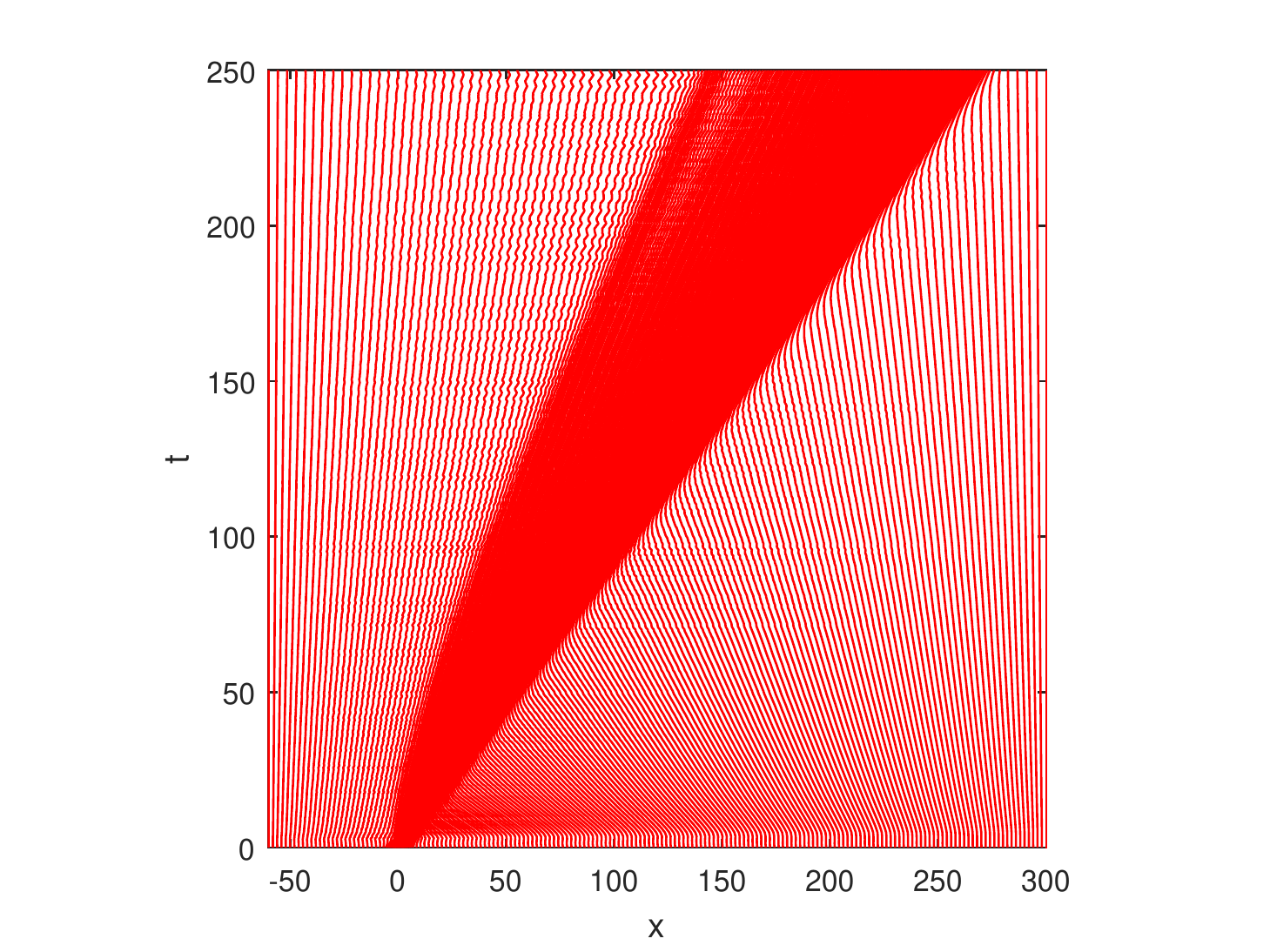}
\end{minipage}
}
\hspace{2mm}
\centering
\hbox{
\begin{minipage}[t]{2in}
\centerline{(d) $d = 5$ (fixed mesh)}
\includegraphics[width=2.5in]{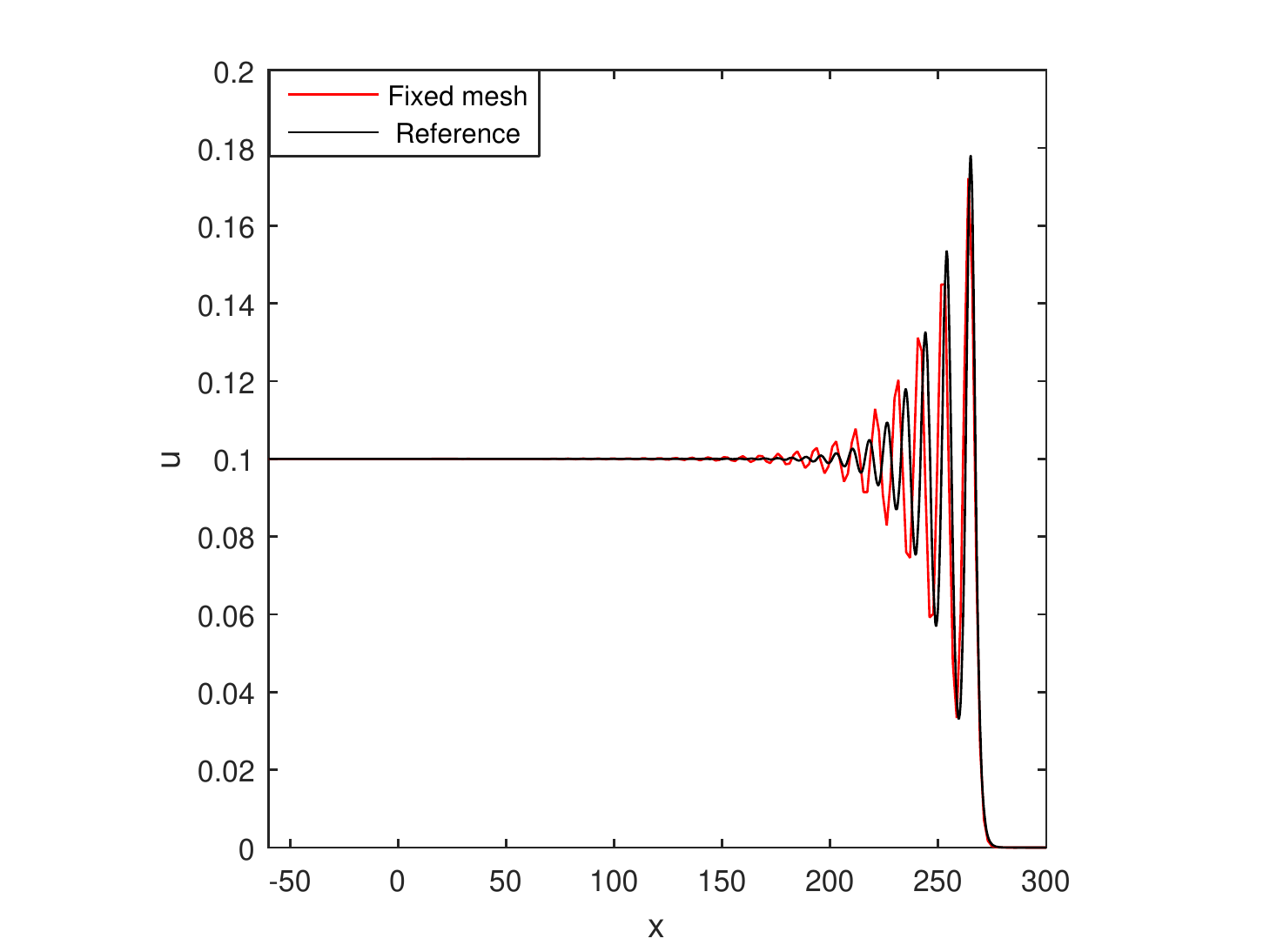}
\end{minipage}
\hspace{2mm}
\begin{minipage}[t]{2in}
\centerline{(e) $d = 5$ (moving mesh)}
\includegraphics[width=2.5in]{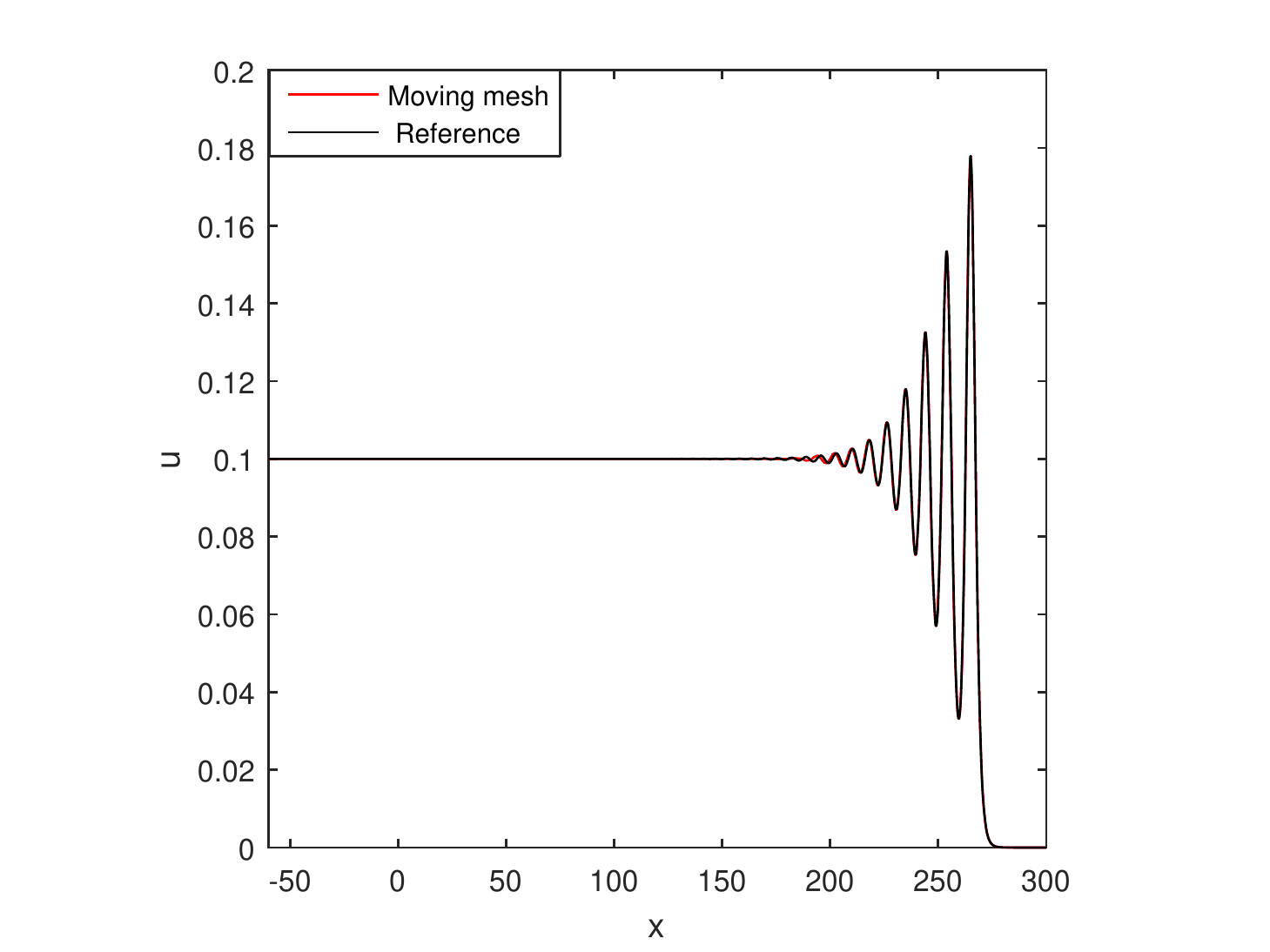}
\end{minipage}
\hspace{2mm}
\begin{minipage}[t]{2in}
\centerline{(f) $d = 5$}
\includegraphics[width=2.5in]{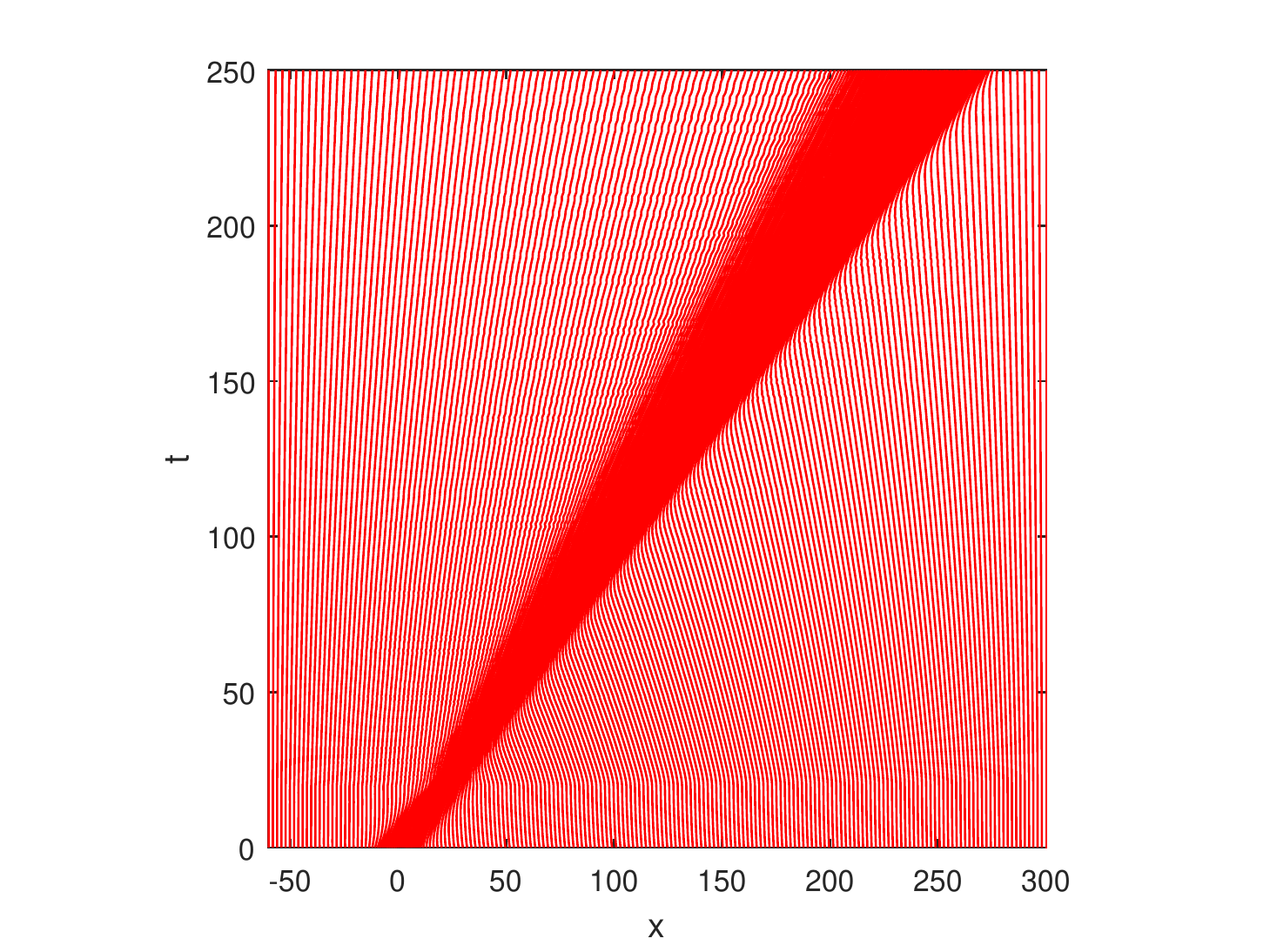}
\end{minipage}
\hspace{2mm}
}
\caption{Example~\ref{exam3.3}. The numerical solutions at $t=250$ obtained with fixed and moving meshes
for the 1D RLW equation with undular bore ($N=200$). The reference solution is obtained with a fixed mesh of
$N = 6000$.}
\label{exam-2.3-1}
\end{figure}

\begin{figure}[thb]
\centering
\hbox{
\begin{minipage}[t]{3in}
\centerline{(a) $d = 2$}
\centering
\includegraphics[width=2.5in]{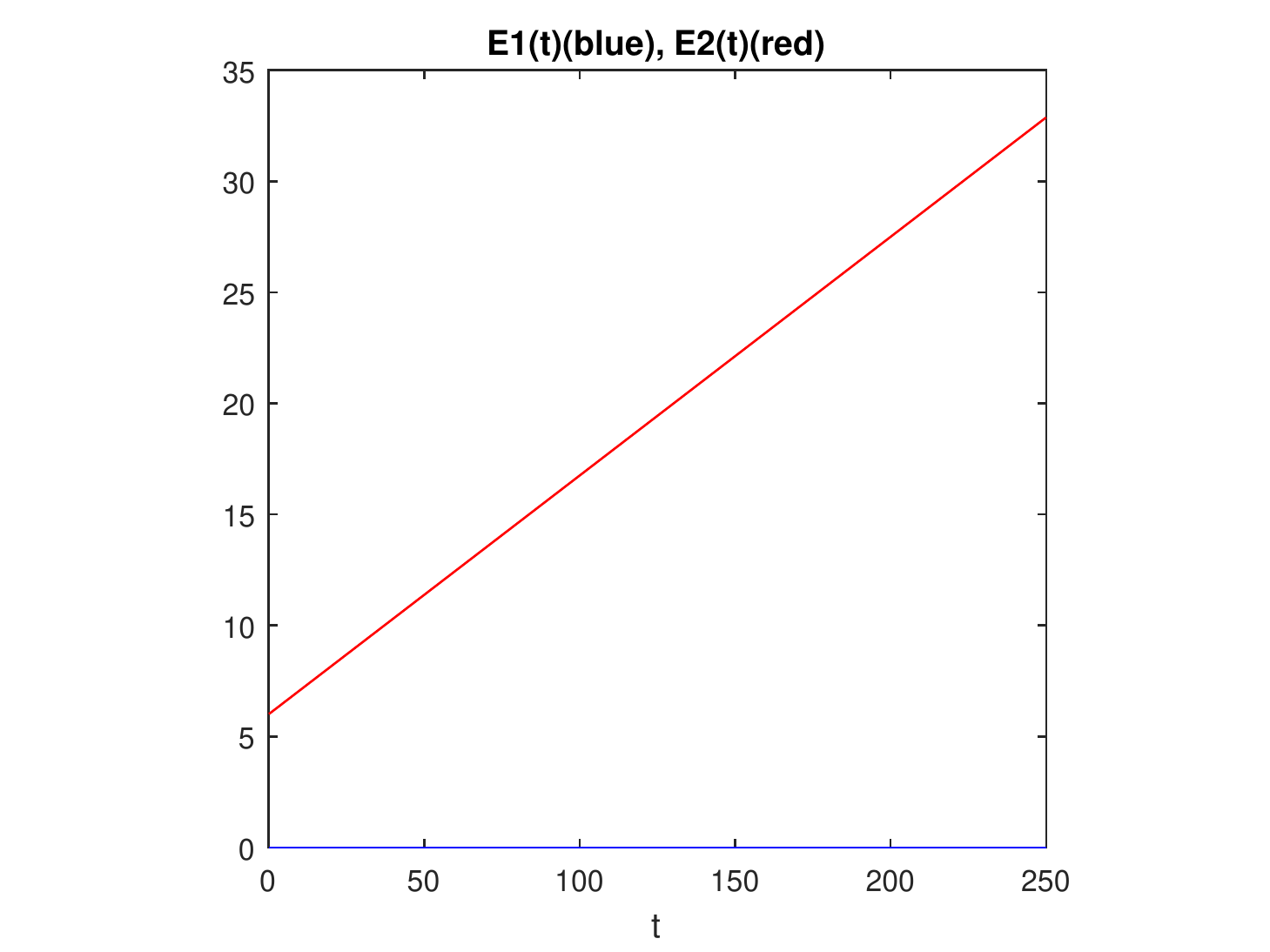}
\end{minipage}
\hspace{5mm}
\begin{minipage}[t]{3in}
\centerline{(b) $d = 5$}
\centering
\includegraphics[width=2.5in]{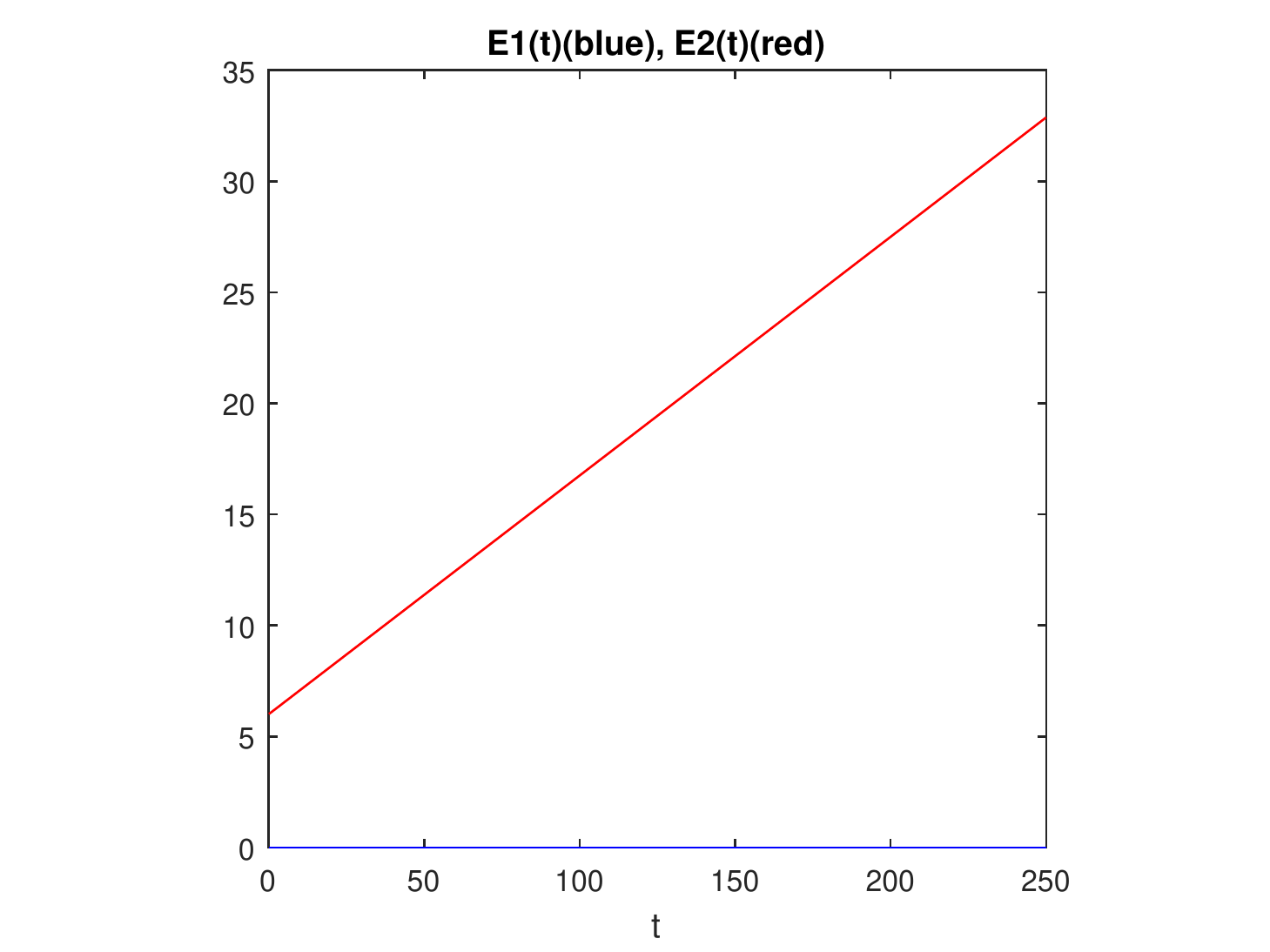}
\end{minipage}
\hspace{5mm}
}
\caption{Example~\ref{exam3.3}. The quantities $E_1$ and $E_2$
are plotted as functions of time. $E_1$ is almost zero and its graph is indistinguishable from
the $x$-axis.
}
\label{exam-2.3-2}
\end{figure}

\begin{exam}
\label{exam3.3}
(1D RLW with undular bore)
We consider the development of an undular bore (e.g., see \cite{Mei2012}) for
the 1D RLW equation (\ref{RLW-1d}) with the initial condition
\[
u(x,0) = \frac{u_0}{2} \left (  1-\tanh\left (\frac{x-x_0}{d}\right ) \right ),
\]
where  $\gamma = 1.5$, $\mu = 1/6$, $u_0 = 0.1$, $x_0 = 0$, and $d = 2$ or $5$.
The boundary condition is $u = u_0$ at $x = -60$ (upstream) and $u = 0$ at $x = 300$ (downstream).
In this example,  $u$ can be thought as the water depth above the equilibrium level
and $d$ as the slope between the still water and deeper water.
The computation is done until $t = 250$.
Due to the continuous injection at the left boundary and the finite propagation velocity, 
the undular bore forms and then is expanding its range as time evolves.

Numerical solutions at $t= 250$ and mesh trajectories with $N=200$ are shown in Fig.~\ref{exam-2.3-1}
for fixed and moving meshes. A solution obtained with the fixed mesh of $N=6000$ is used as the reference solution.
It can be seen that the solution obtained with a moving mesh is more accurate than that with a fixed mesh
of the same number of elements and the mesh concentration reflects correctly
the development of the undular bore.
The quantities $E_1$ and $E_2$ are plotted in Fig.~\ref{exam-2.3-2}.
As the water coming from the left boundary at a constant rate, these quantities grow linearly with time.
Nevertheless, $E_1$ remains very small, almost indistinguishable from
the $x$-axis. (Recall that the error in preserving $E_1$ on a moving mesh is at the level of roundoff error
for a sufficiently fine mesh.)
Finally, numerical results show that the undular bore is very stable.
\end{exam}

\begin{figure}[thb]
\centering
\hbox{
\begin{minipage}[t]{2in}
\centerline{(a): $\mu = 1$ (fixed mesh)}
\includegraphics[width=2.5in]{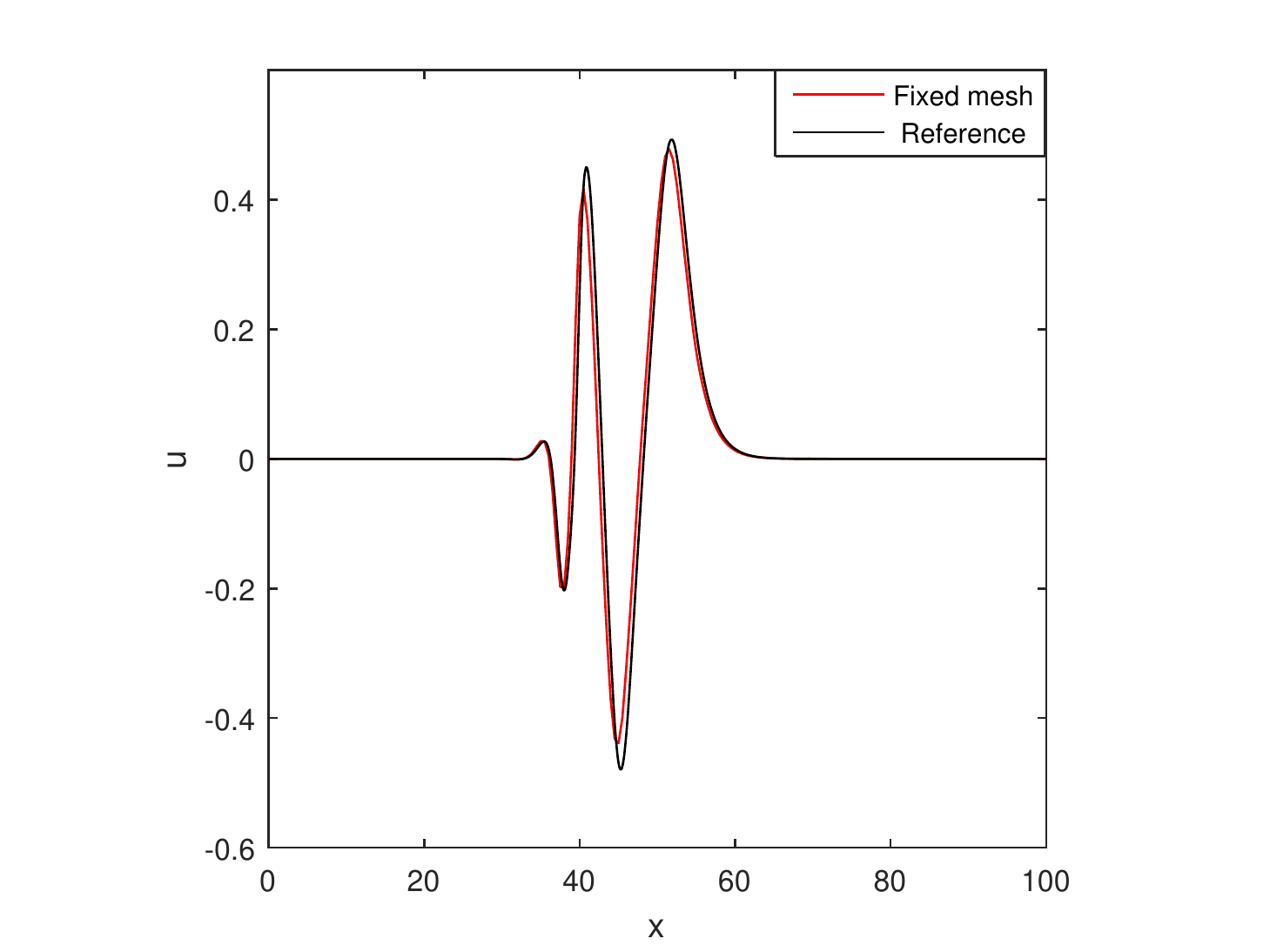}
\end{minipage}
\hspace{2mm}
\begin{minipage}[t]{2in}
\centerline{(b): $\mu = 1$ (moving mesh)}
\includegraphics[width=2.5in]{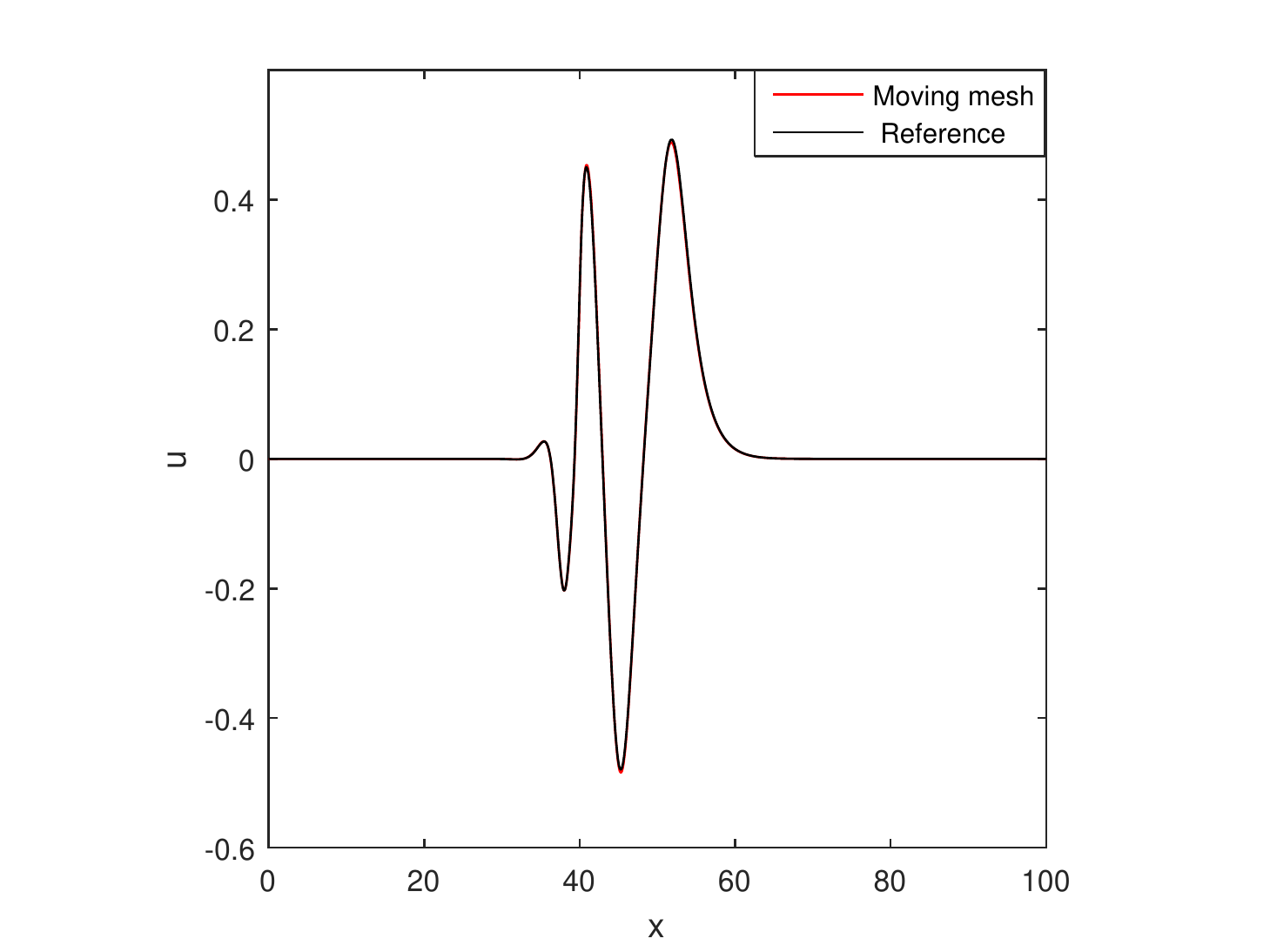}
\end{minipage}
\hspace{2mm}
\begin{minipage}[t]{2in}
\centerline{(c): $\mu = 1$}
\includegraphics[width=2.5in]{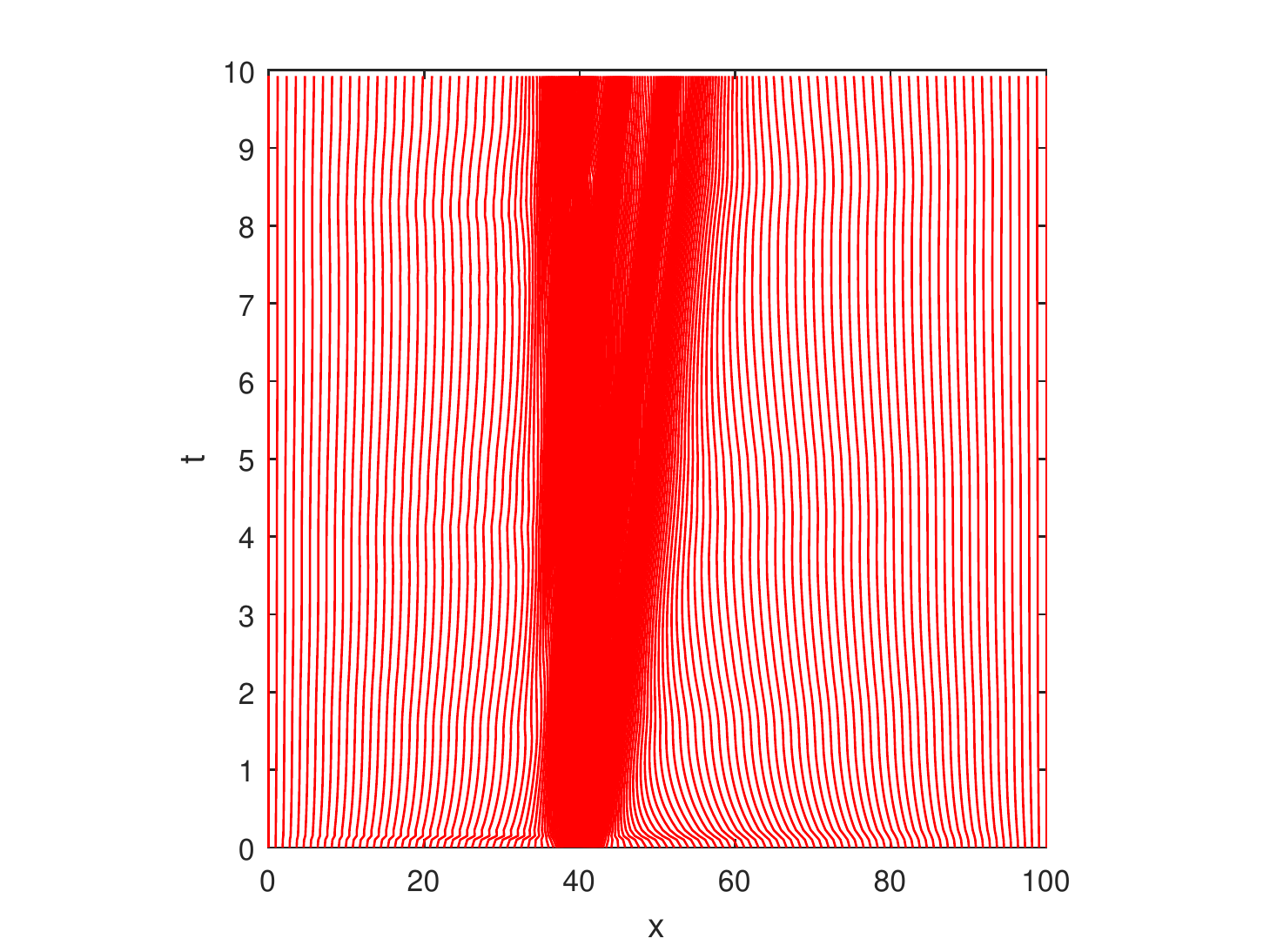}
\end{minipage}
}
\hspace{2mm}
\centering
\hbox{
\begin{minipage}[t]{2in}
\centerline{(d): $\mu = 0.5$ (fixed mesh)}
\includegraphics[width=2.5in]{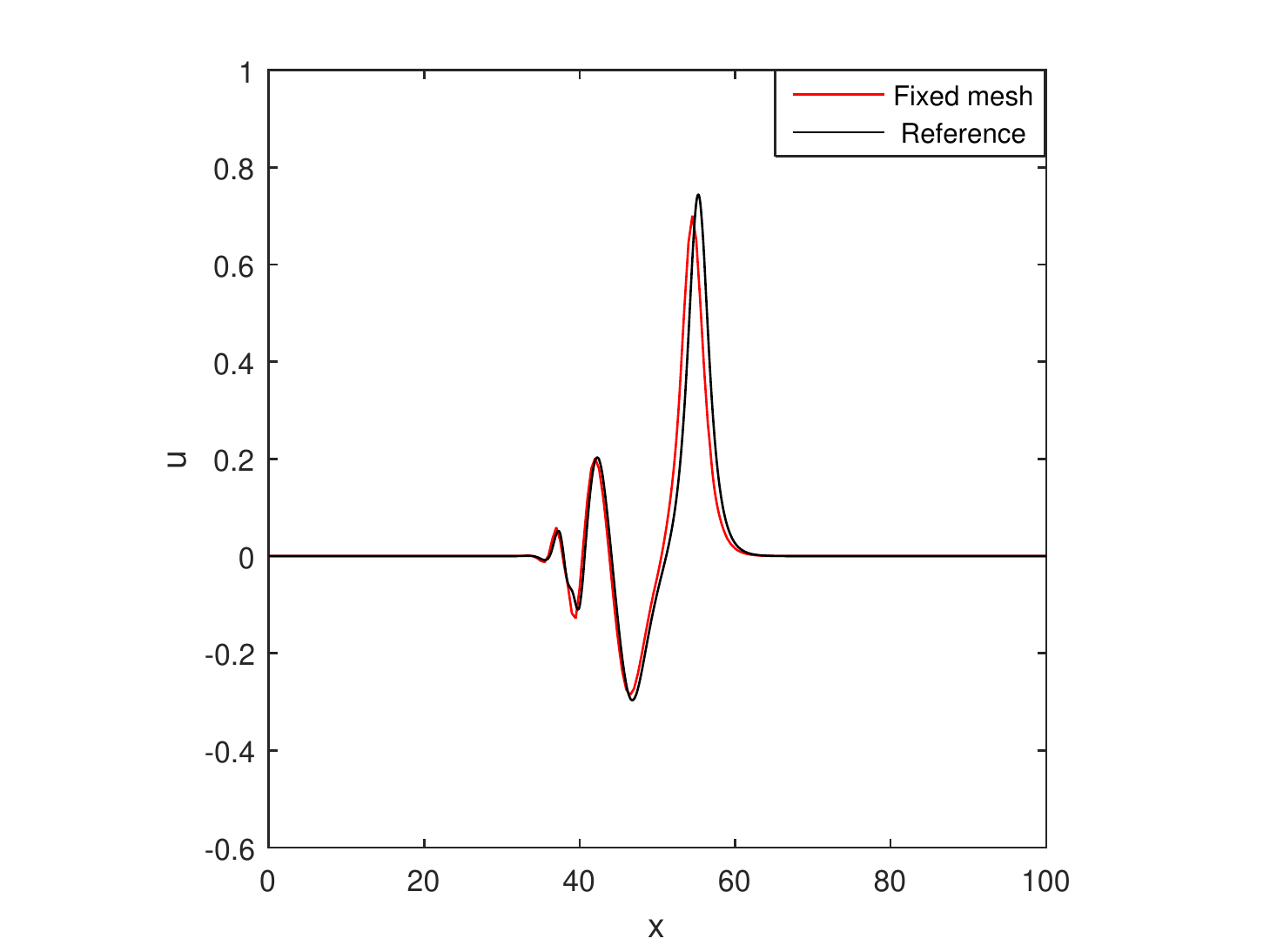}
\end{minipage}
\hspace{2mm}
\begin{minipage}[t]{2in}
\centerline{(e): $\mu = 0.5$ (moving mesh)}
\includegraphics[width=2.5in]{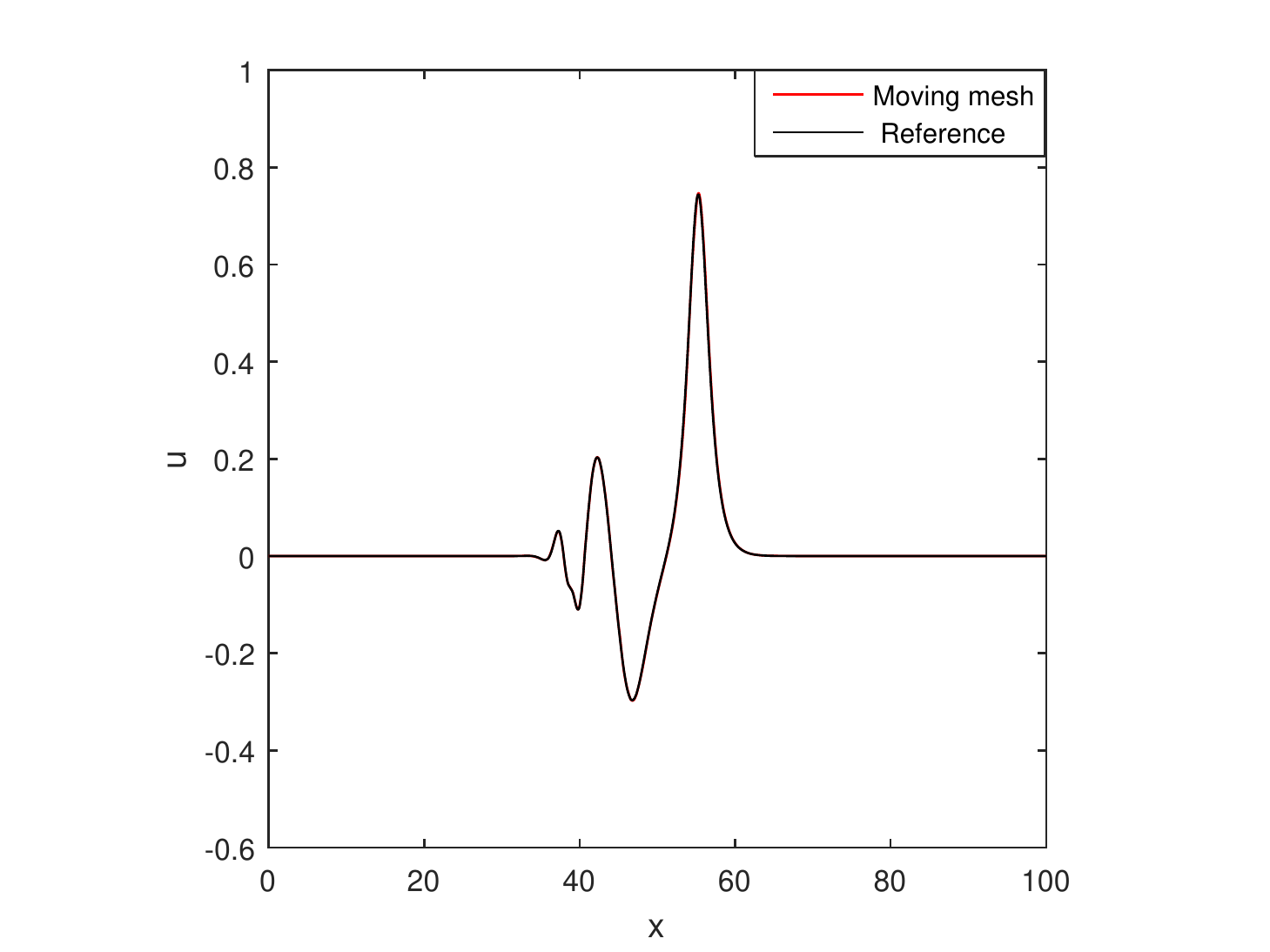}
\end{minipage}
\hspace{2mm}
\begin{minipage}[t]{2in}
\centerline{(f): $\mu = 0.5$}
\includegraphics[width=2.5in]{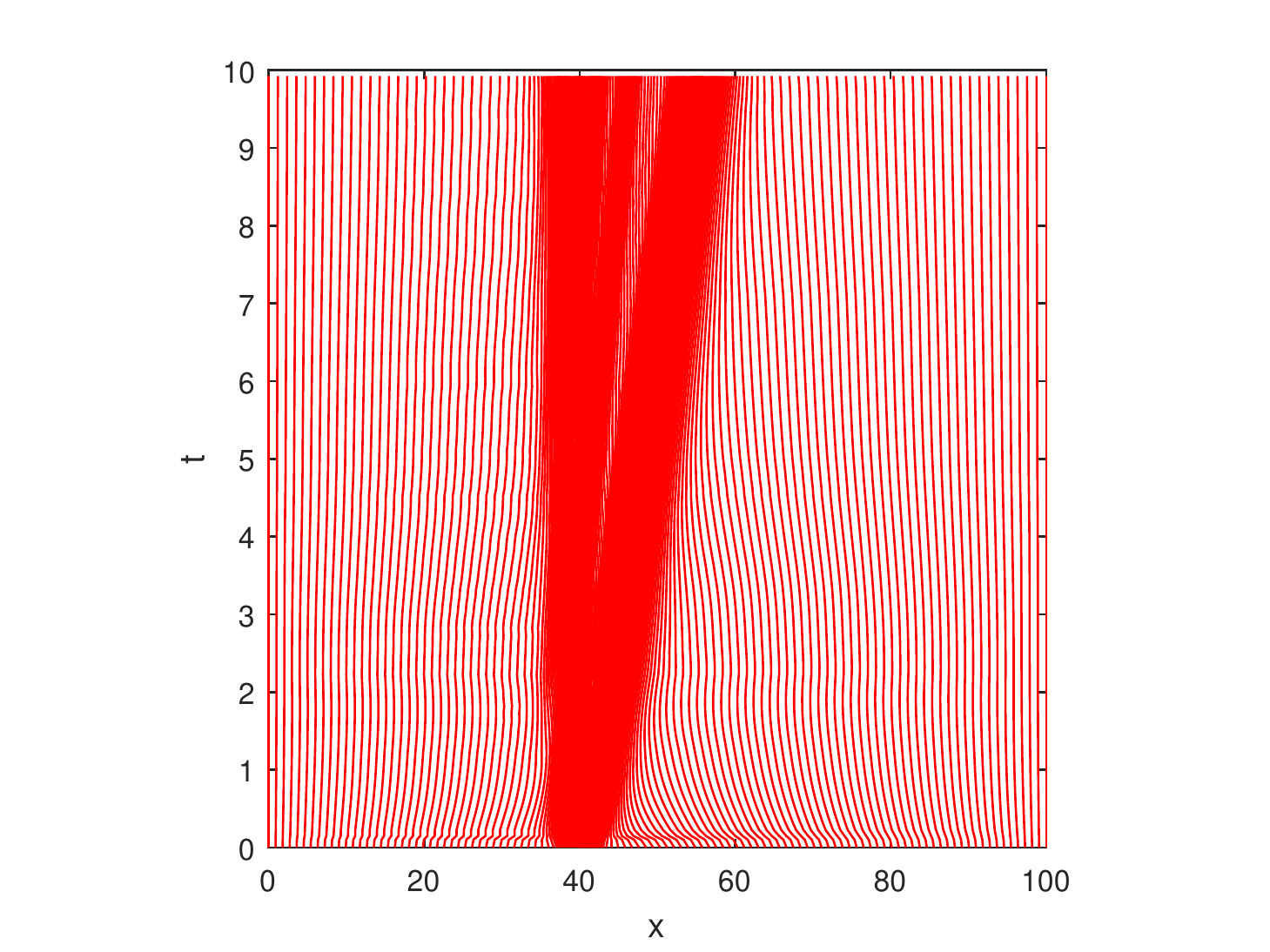}
\end{minipage}
\hspace{2mm}
}
\caption{Example~\ref{exam3.4}. Numerical solutions and mesh trajectories
are obtained with fixed and moving meshes of $N=200$ for the 1D modified RLW equation
with the Maxwellian initial condition. The reference solution is obtained with a fixed mesh
of $N=6000$.}
\label{exam-3.4}
\end{figure}

\begin{exam}
\label{exam3.4}
(1D modified RLW with the Maxwellian initial condition)
In this test, we consider the 1D modified RLW (MRLW) equation
\[
u_t + u_x + \gamma u^2  u_x - \mu u_{xxt} = 0
\]
subject to a homogeneous Dirichlet boundary condition and the Maxwellian initial condition \cite{Gao2015}
\[
u(x,0) = e^{-(x-40)^2}.
\]
We take $\gamma = 6$ and $\mu = 1$ or $\mu = 0.5$.
For the time being, the Maxwellian initial condition develops into a train of solitary waves,
with the wave number and amplitude depending on the value of $\mu$.
The smaller $\mu$ is, the more solitary waves will form.


The computation is performed with $T = 10$ and $\Omega = (0,100)$.
Numerical results obtained with fixed and moving meshes are shown
in Fig.~\ref{exam-3.4}. It can be seen that the solution with a moving mesh
is more accurate than that with a fixed mesh and, indeed, the former is almost indistinguishable
from the reference solution which is obtained with a fixed mesh of $N = 6000$.
Numerical experiment also shows that the train of the solitons are stable.
\end{exam}

\begin{figure}[thb]
\centering
\hbox{
\begin{minipage}[t]{3in}
\centerline{(a): On moving mesh}
\includegraphics[width=2.5in]{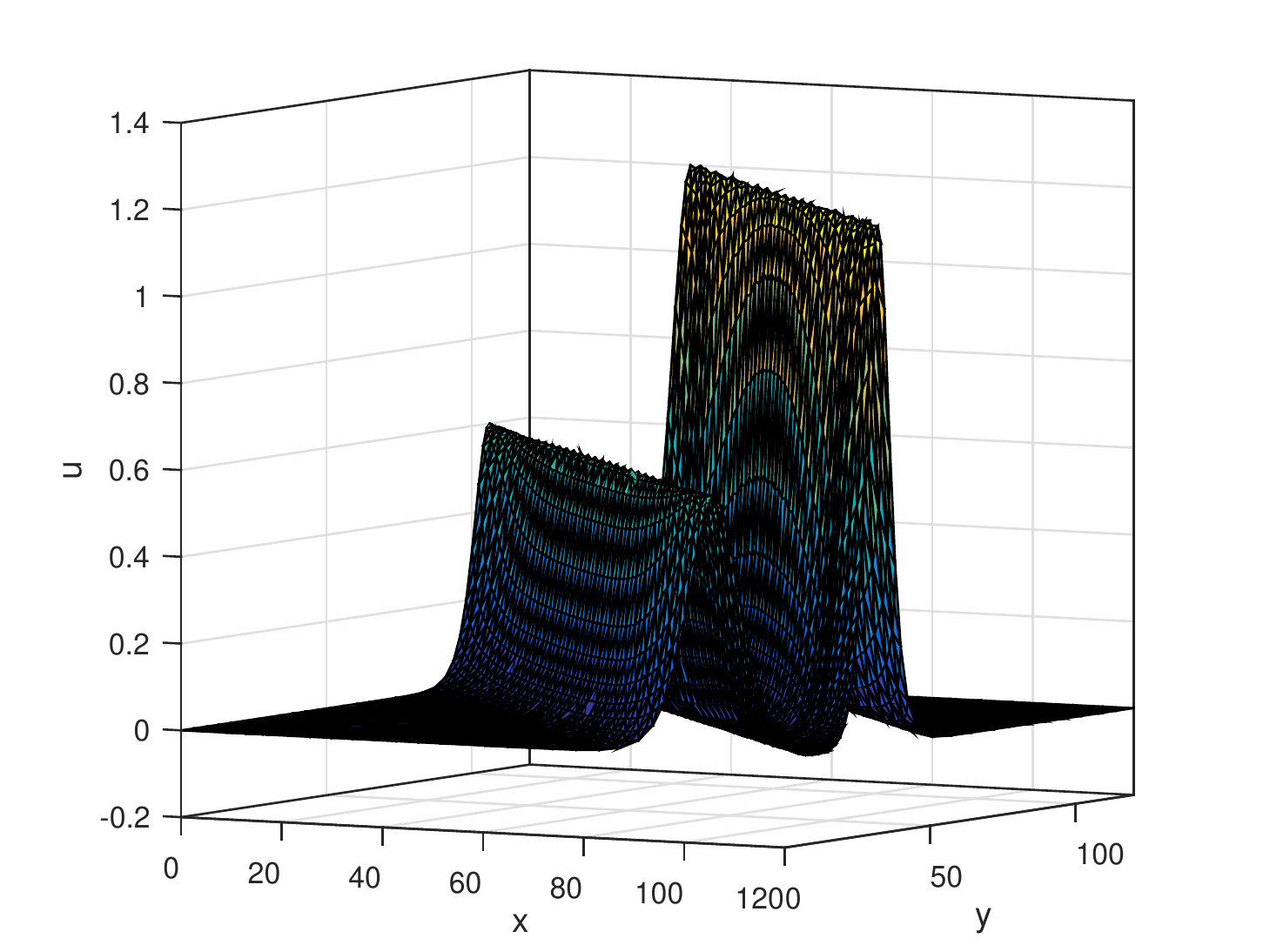}
\end{minipage}
\hspace{5mm}
\begin{minipage}[t]{3in}
\centerline{(b): Moving mesh}
\includegraphics[width=2.5in]{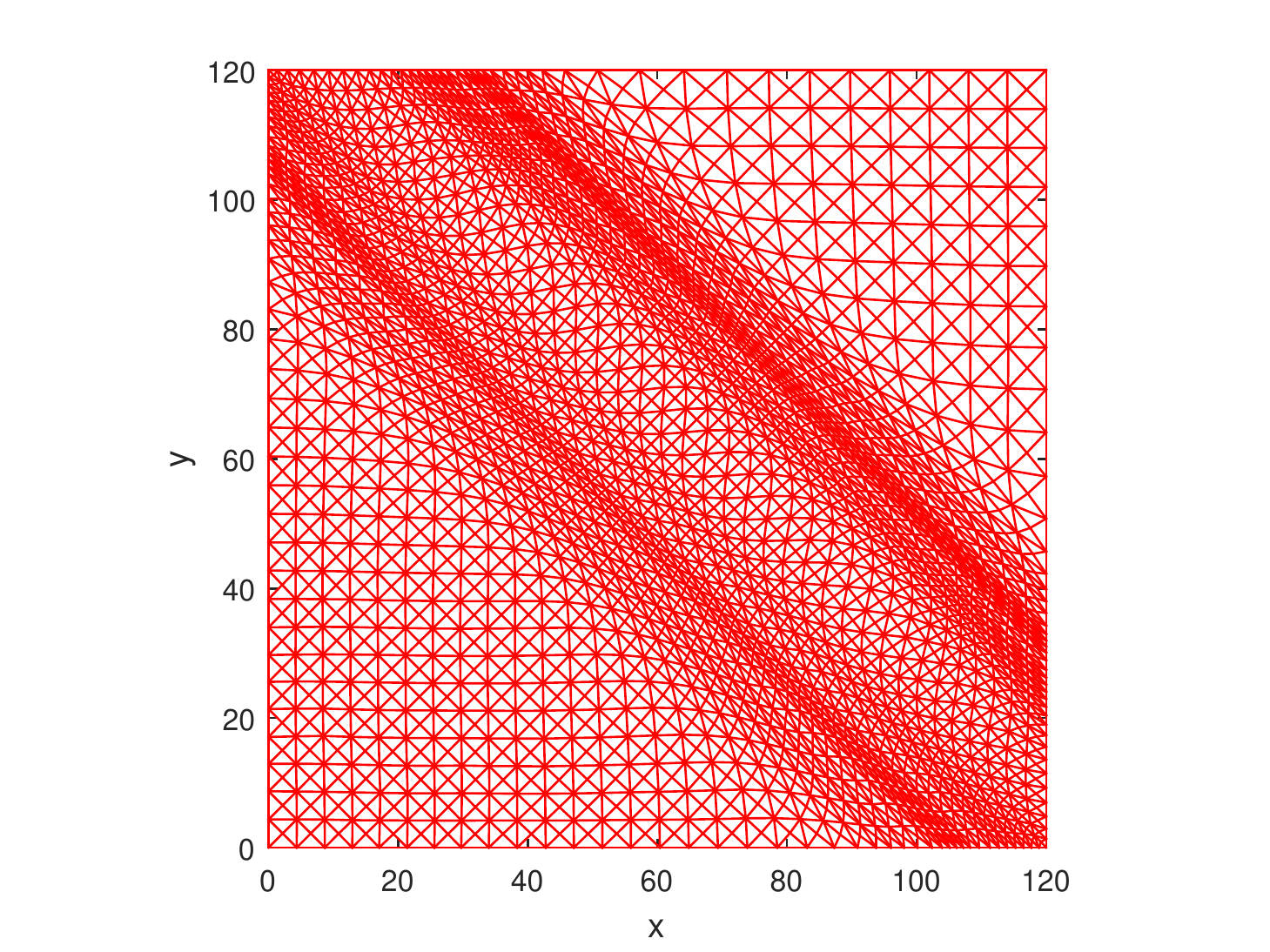}
\end{minipage}
\hspace{5mm}
}
\centering
\hbox{
\begin{minipage}[t]{3in}
\centerline{(c): On fixed mesh}
\includegraphics[width=2.5in]{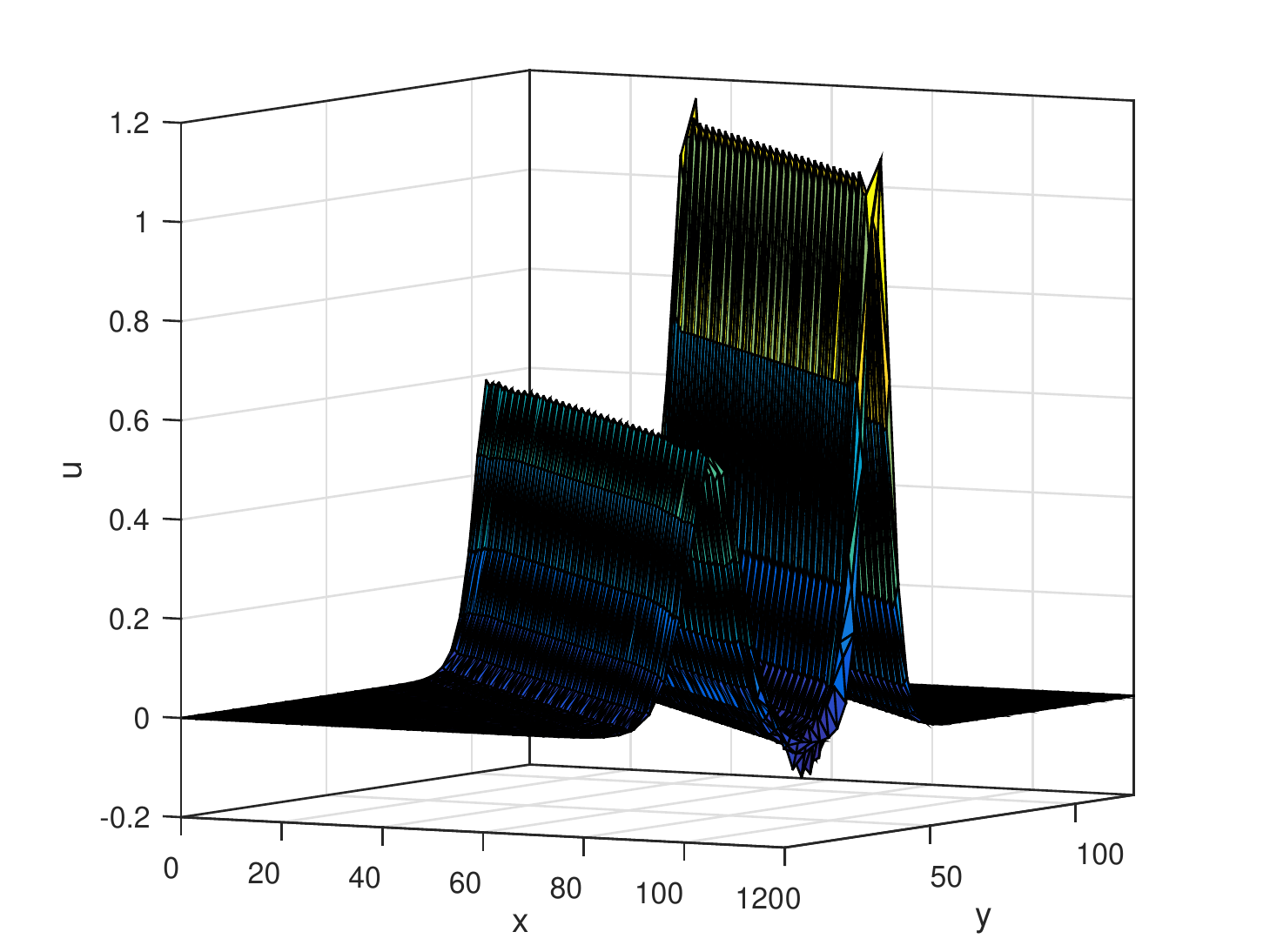}
\end{minipage}
\hspace{5mm}
\begin{minipage}[t]{3in}
\centerline{(d): Fixed mesh}
\includegraphics[width=2.5in]{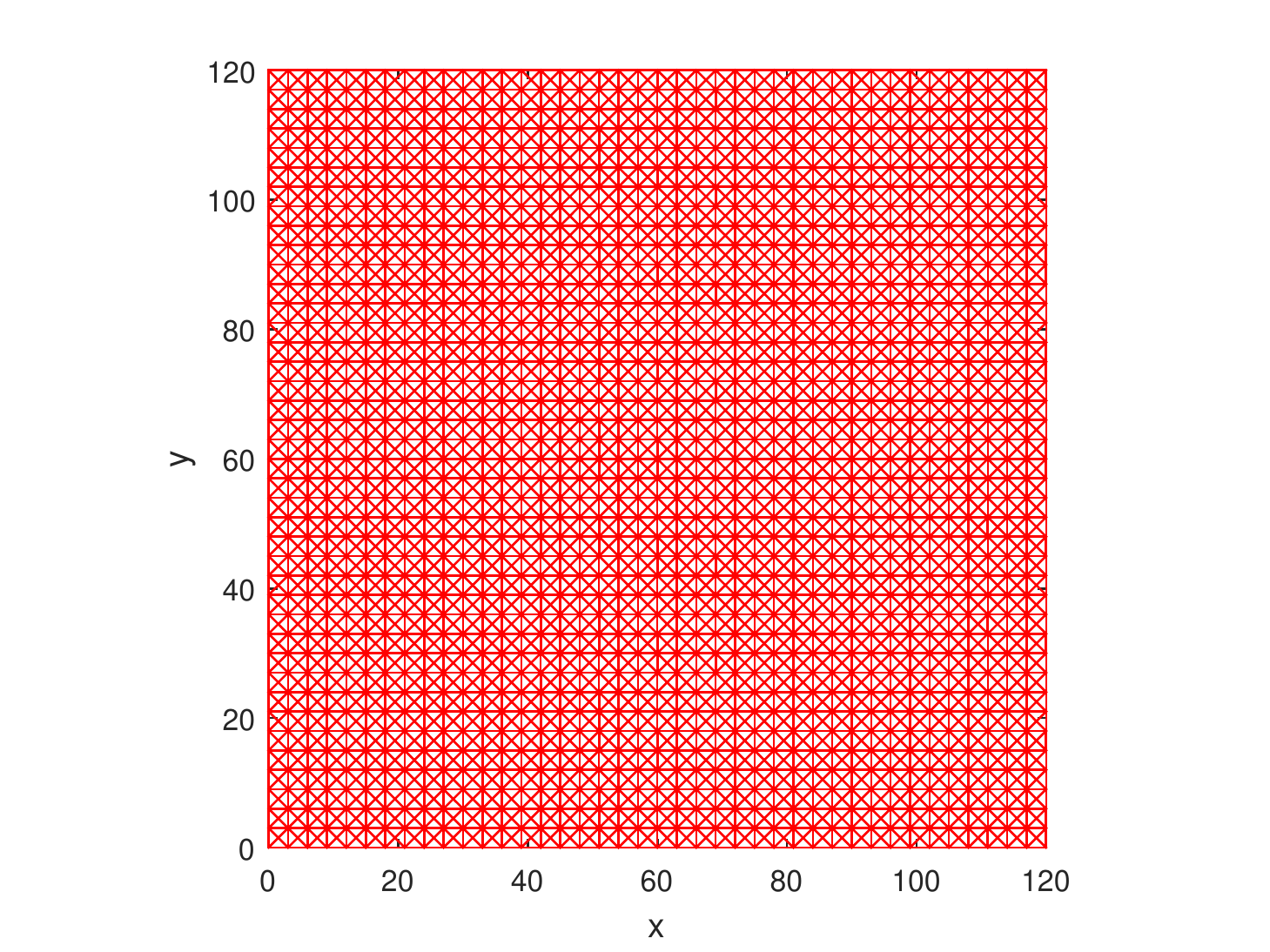}
\end{minipage}
\hspace{5mm}
}
\caption{Example~\ref{exam3.5}. Numerical solutions and meshes at $t = 15$ for fixed and
moving meshes of $N=6400$.
}
\label{exam-3.5-1}
\end{figure}

\begin{exam}
\label{exam3.5}
(2D RLW with two solitary waves)
In this test we consider the 2D RLW equation (\ref{RLW-1}) with
$\alpha = \beta = \gamma = \delta = \mu = 1$.
The Dirichlet and initial conditions are chosen such that
the exact solution is given by
\[
u(x,y,t) = \sum_{j=1}^2
3c_j \text{sech}^2\left (k_j (x+y-v_j t-x_j-y_j)\right ),
\]
where $k_j = \frac{1}{2}\sqrt\frac{c_j}{2(1+c_j)}$, $v_j = 2(1+c_j)$,
$c_1= 0.2$, $c_2 = 0.4$, $ v_1= 2.4$, $v_2 = 2.8$, $x_1=y_1 = 35$,
and $x_2 = y_2 = 55$.
Notice that $3c_j$ is the maximum amplitude and
$v_j$ is the circular frequency.
The computation is performed on $\Omega =  (0,120)\times (0,120)$ with $T = 15$.

Numerical results are shown in Table~\ref{exam-3.5-order} and Fig.~\ref{exam-3.5-1}. They
indicate that the finite element method is second order for both fixed and moving meshes.
Moreover, a moving mesh leads to more accurate solutions, with roughly an order of magnitude smaller error,
than a fixed mesh of the same number of elements.

\end{exam}

\begin{table}[htb]
\begin{center}
\caption{Example~\ref{exam3.5}. $L^2$ and $L^{\infty}$ error and convergence order for the 2D RLW equation.}
\begin{tabular}{c|l|l|l|l|l|l|l|l}\hline \hline
                       & \multicolumn{4}{|c}{Moving Mesh}& \multicolumn{4}{|c}{Fixed Mesh}  \\ \hline
 $N$ & {$L^2$ error} &  order & $L^{\infty}$ error & order & {$L^2$ error} &  order & $L^{\infty}$ error & order\\ \hline\hline
100    &  3.59E-1  &             &   1.73E1   &           &  3.84E-1  &             &  1.77E1   & \\
400    &  1.02E-1  &   1.82   &   4.87E-0   &  1.82  &  2.77E-1  &   0.47   &   1.19E1   &  0.56\\
1600  &  1.45E-2  &   2.81   &   1.02E-0   &  2.25  &  1.32E-1  &   1.07   &   6.72E-0   &  0.83\\
6400  &  2.82E-3  &   2.36   &   1.97E-1   &  2.38  &  3.45E-2  &   1.93   &   2.20E-0   &  1.61\\
25600&  6.24E-4  &   2.18   &   4.18E-2   &  2.24  &  8.34E-3  &   2.05   &   6.11E-1   &  1.85\\ \hline\hline
\end{tabular}
\label{exam-3.5-order}
\end{center}
\end{table}

\begin{figure}[thb]
\centering
\hbox{
\begin{minipage}[t]{2in}
\includegraphics[width=2.3in]{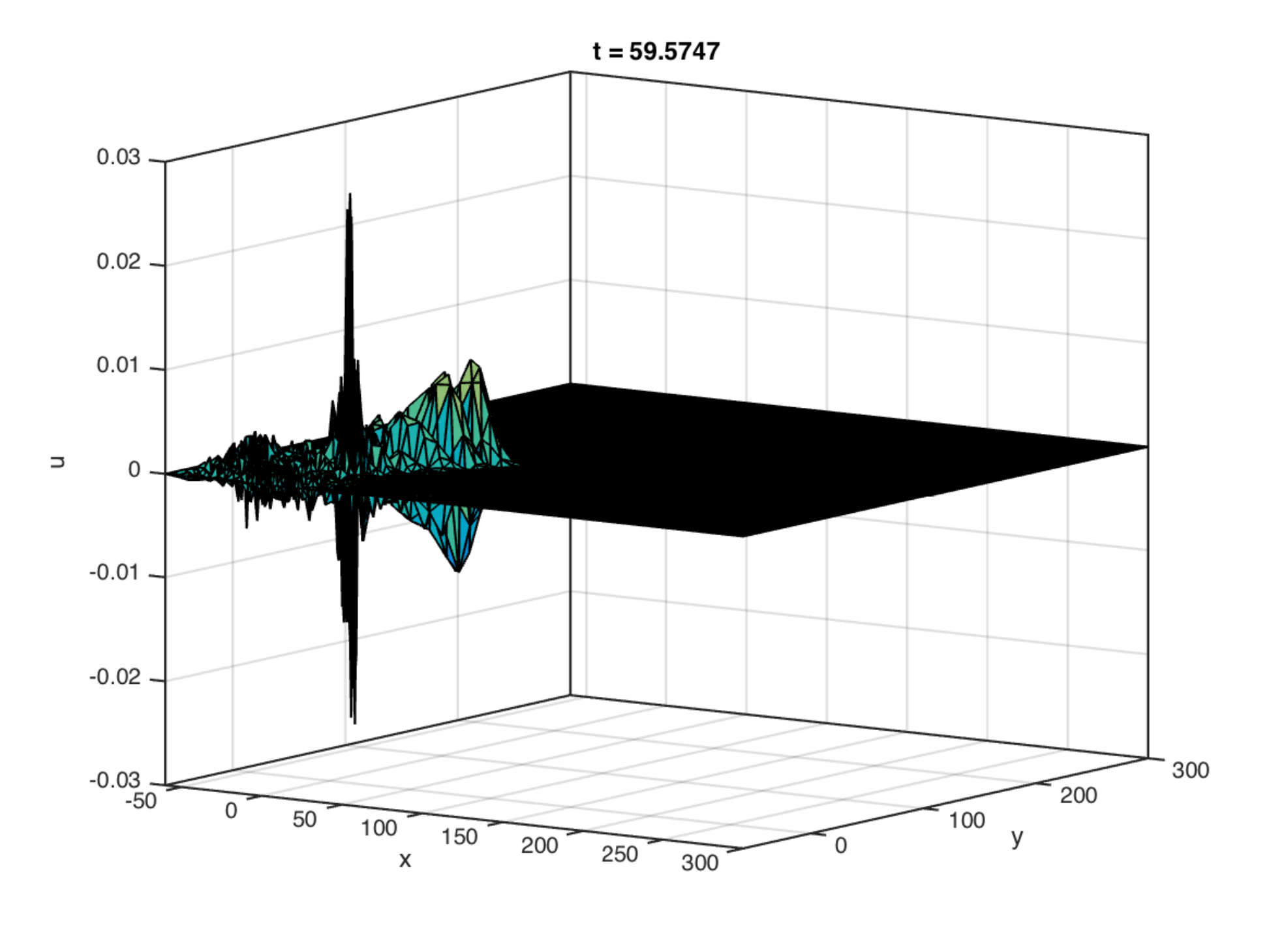}
\end{minipage}
\hspace{2mm}
\begin{minipage}[t]{2in}
\includegraphics[width=2.5in]{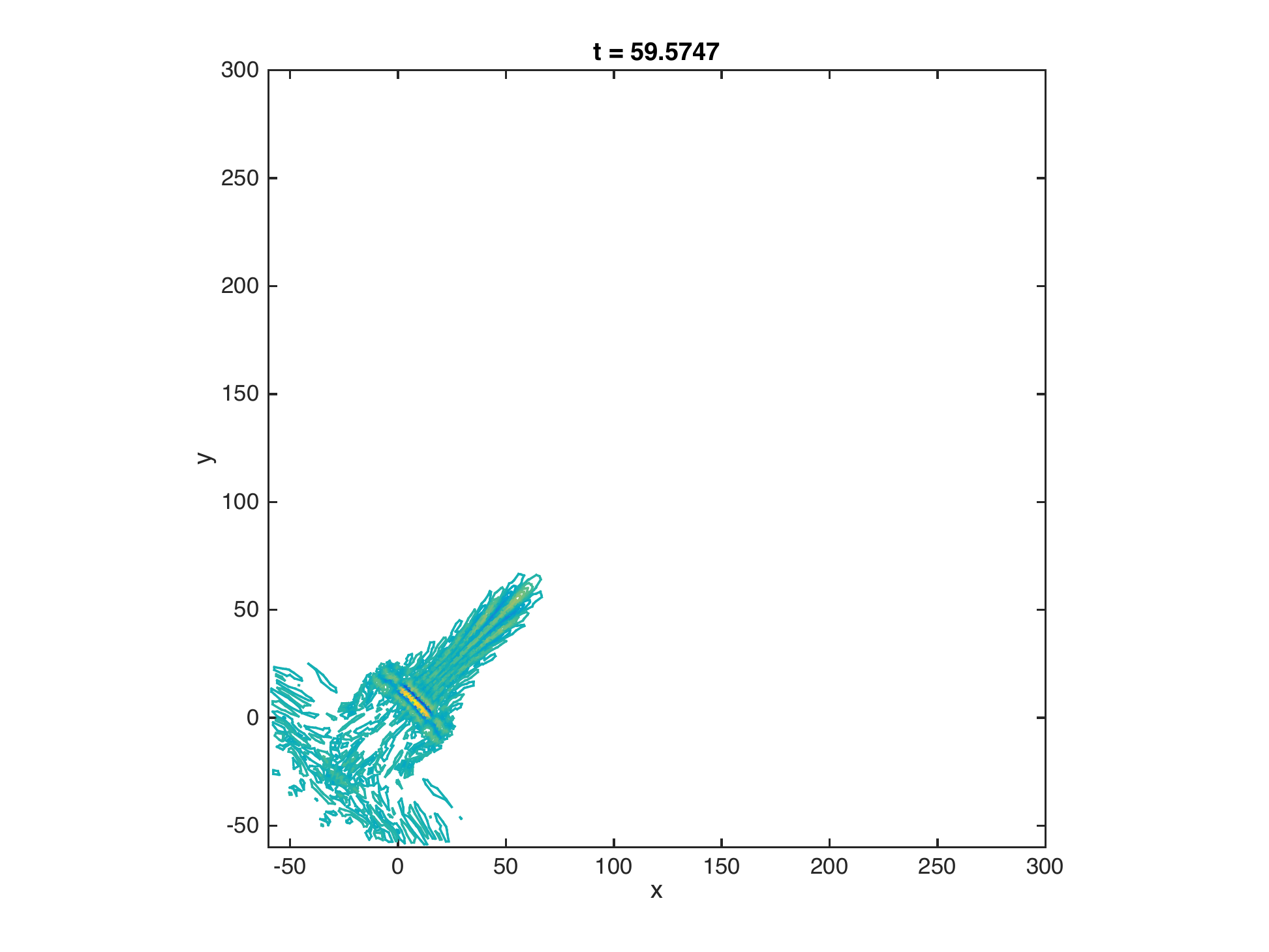}
\end{minipage}
\hspace{2mm}
\begin{minipage}[t]{2in}
\includegraphics[width=2.5in]{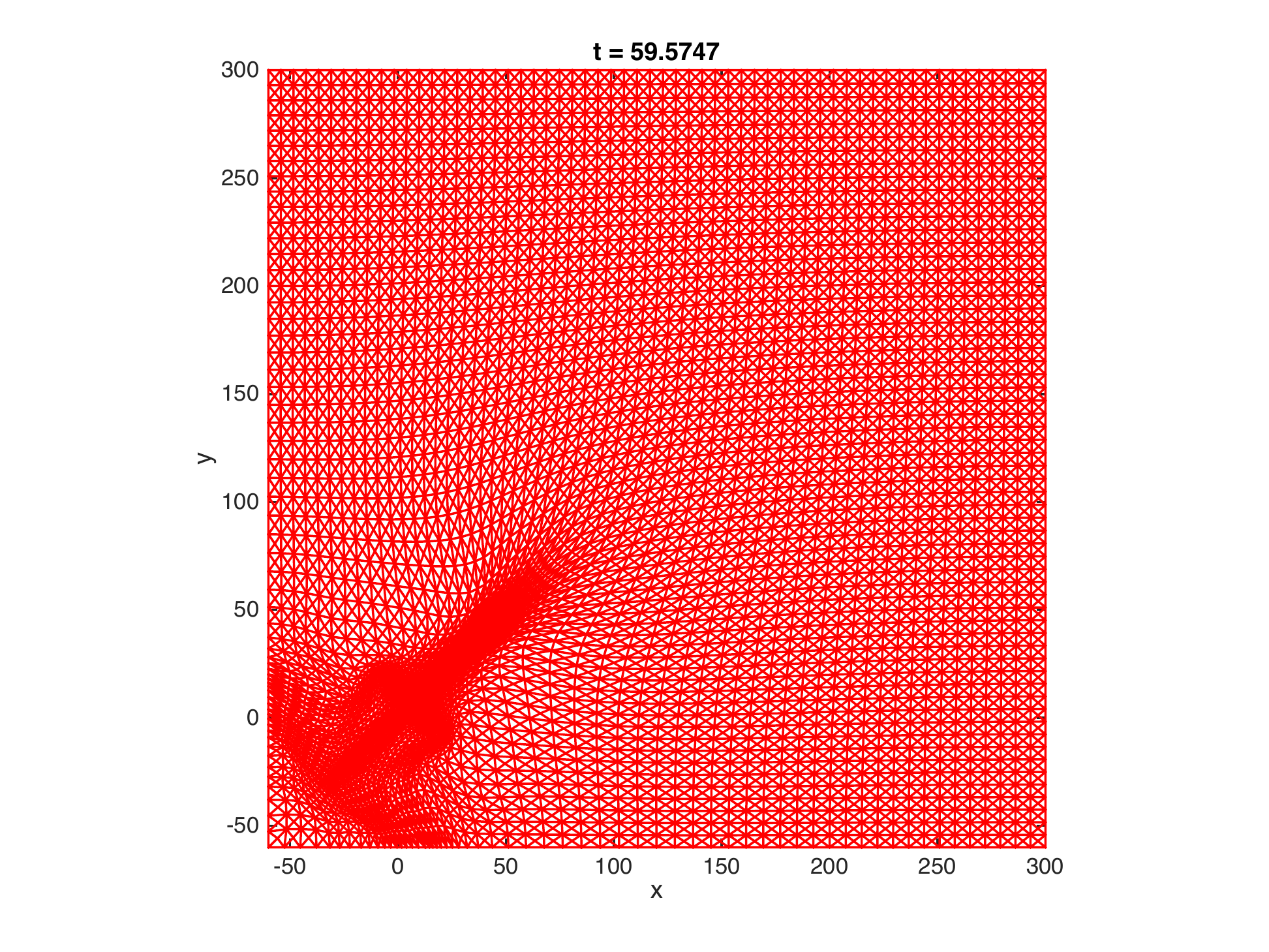}
\end{minipage}
\hspace{2mm}
}
\hbox{
\begin{minipage}[t]{2in}
\includegraphics[width=2.3in]{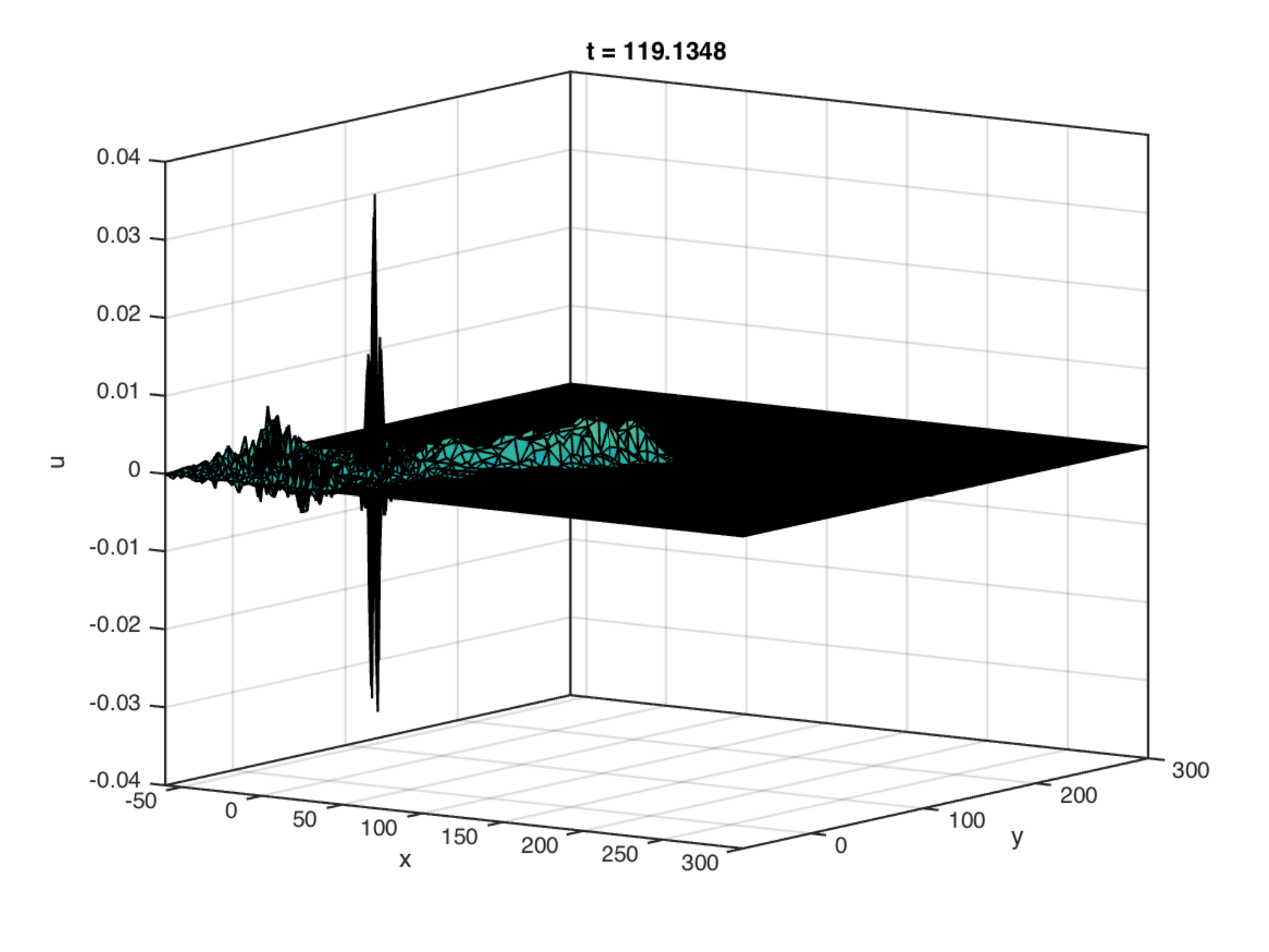}
\end{minipage}
\hspace{2mm}
\begin{minipage}[t]{2in}
\includegraphics[width=2.5in]{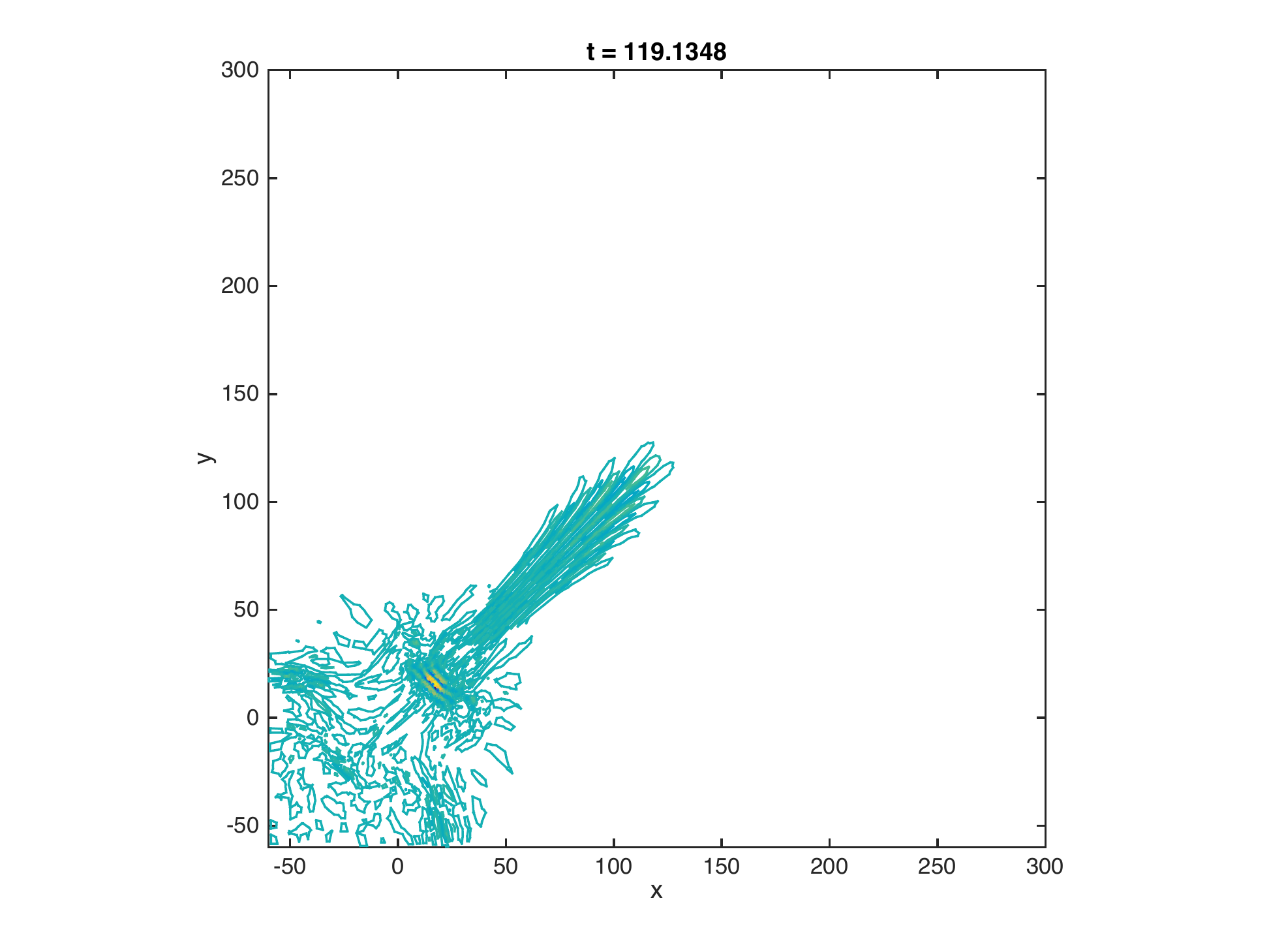}
\end{minipage}
\hspace{2mm}
\begin{minipage}[t]{2in}
\includegraphics[width=2.5in]{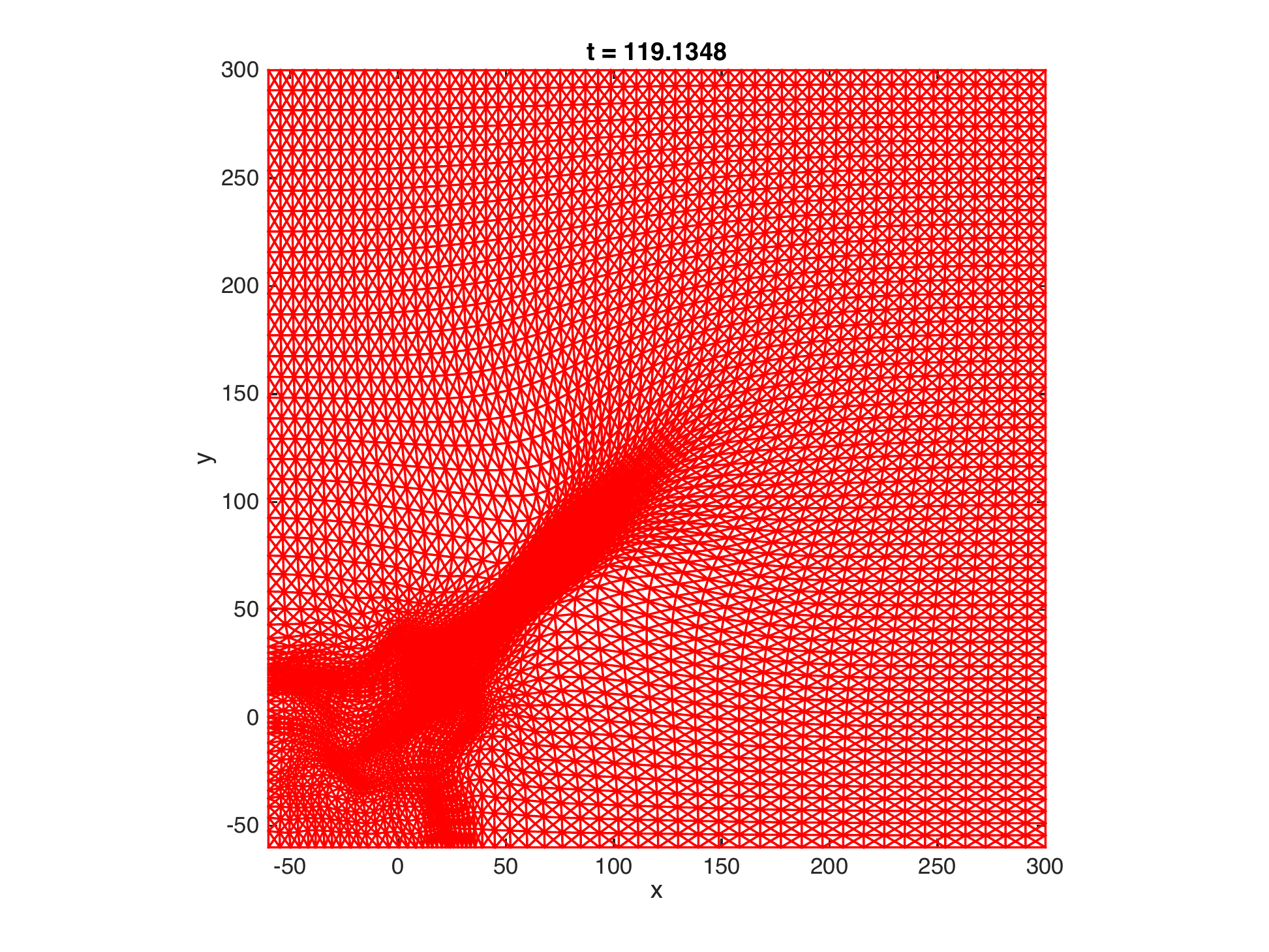}
\end{minipage}
\hspace{2mm}
}
\hbox{
\begin{minipage}[t]{2in}
\includegraphics[width=2.3in]{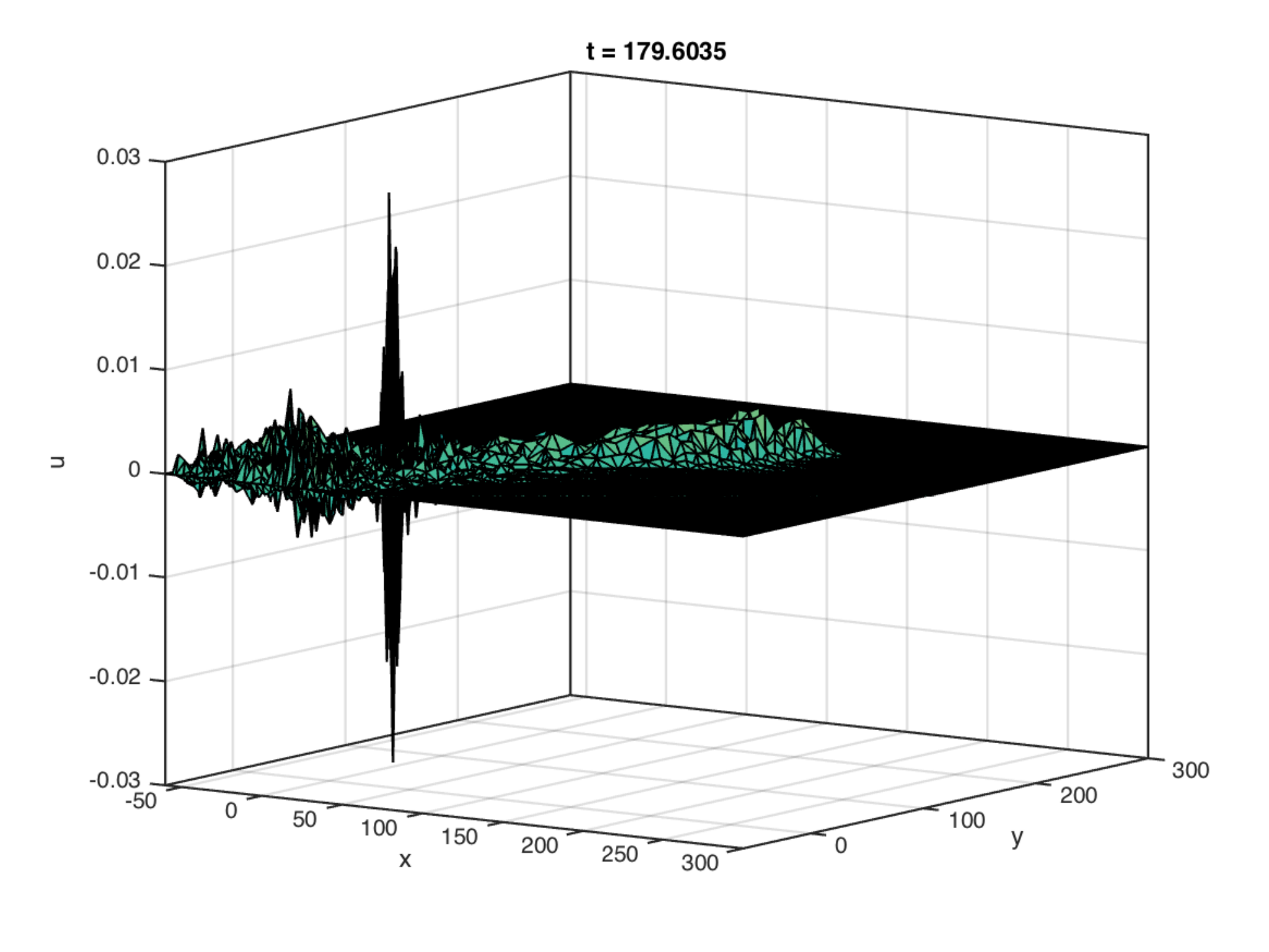}
\end{minipage}
\hspace{2mm}
\begin{minipage}[t]{2in}
\includegraphics[width=2.5in]{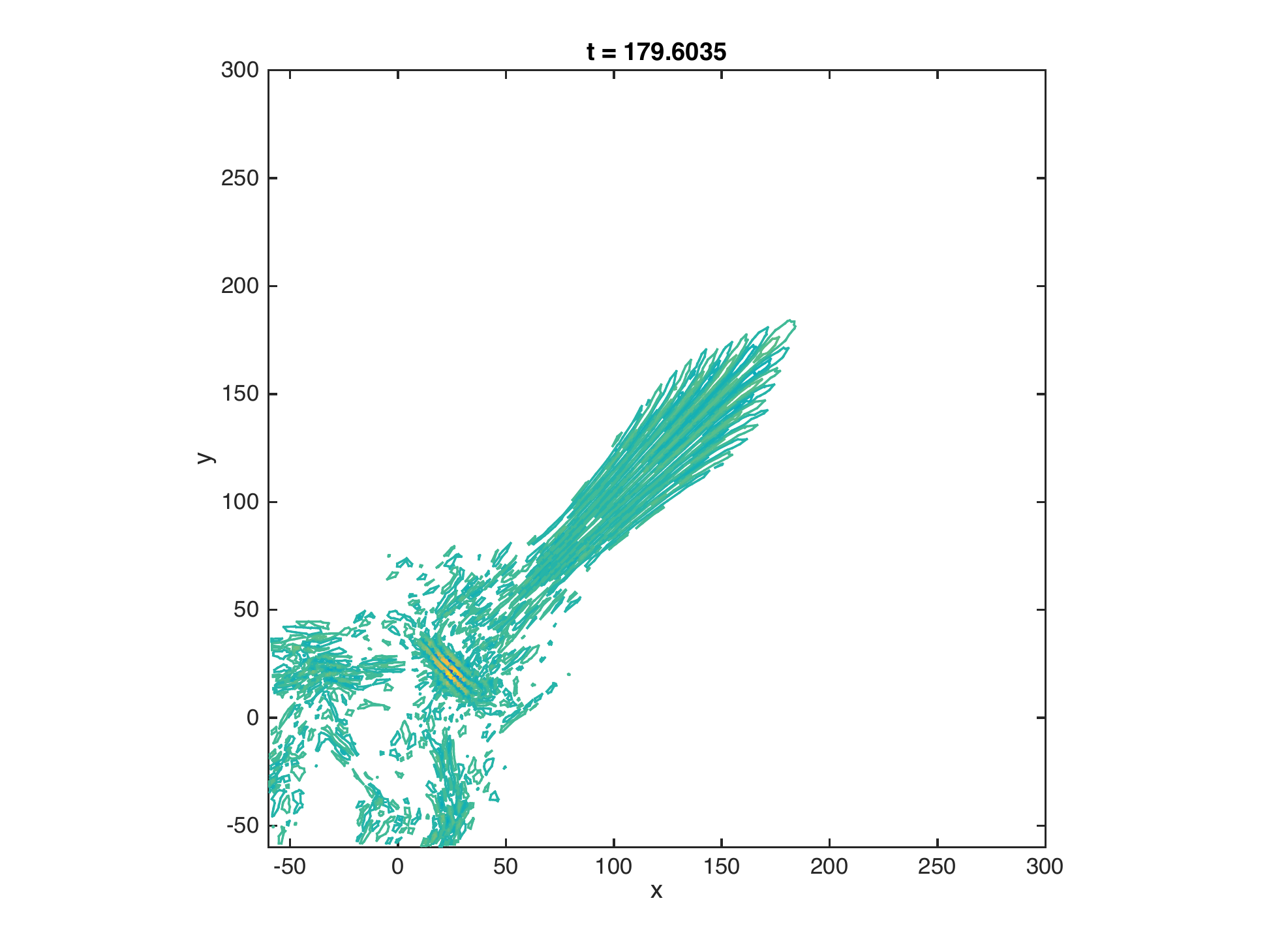}
\end{minipage}
\hspace{2mm}
\begin{minipage}[t]{2in}
\includegraphics[width=2.5in]{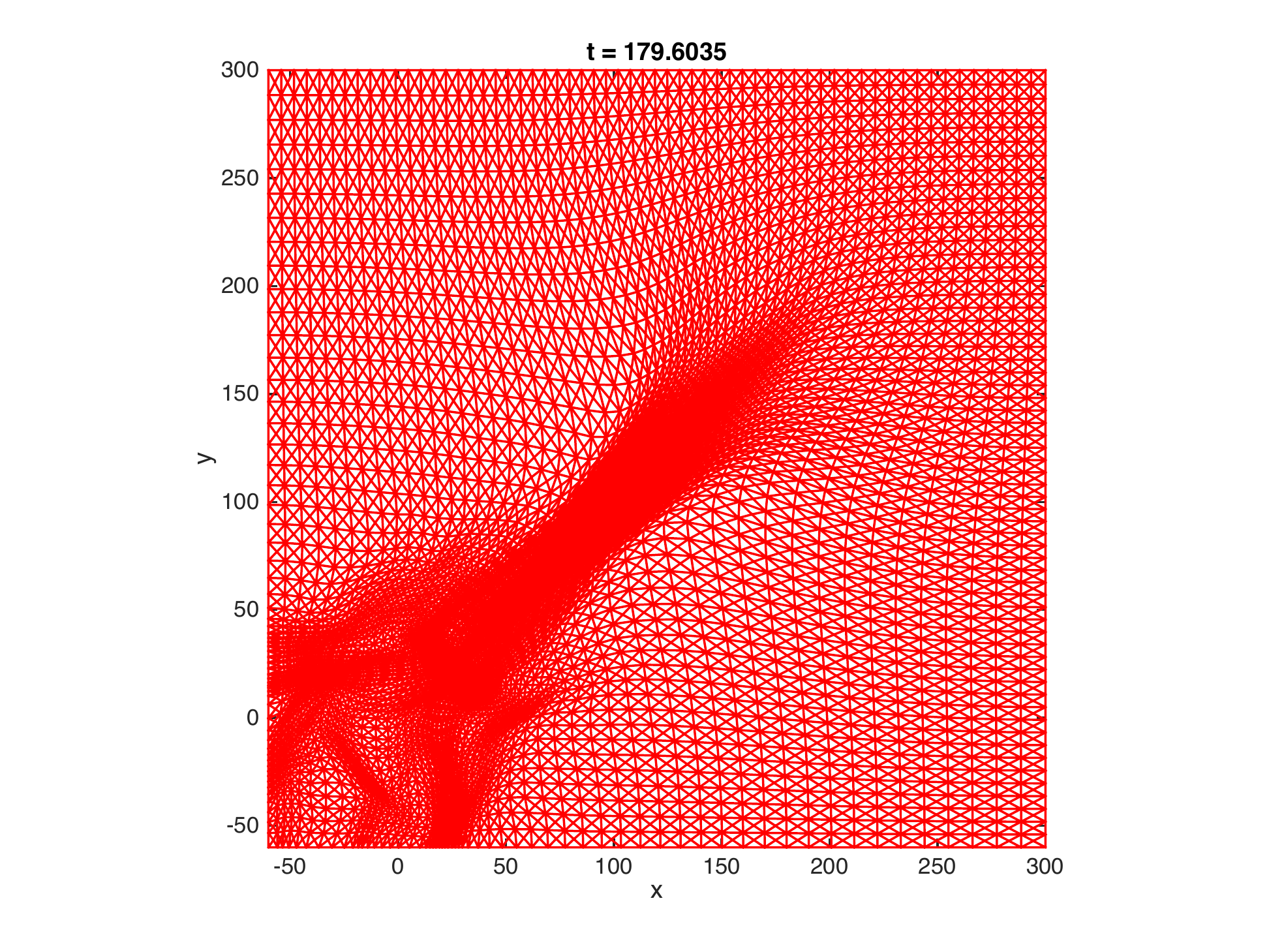}
\end{minipage}
\hspace{2mm}
}
\hbox{
\begin{minipage}[t]{2in}
\includegraphics[width=2.3in]{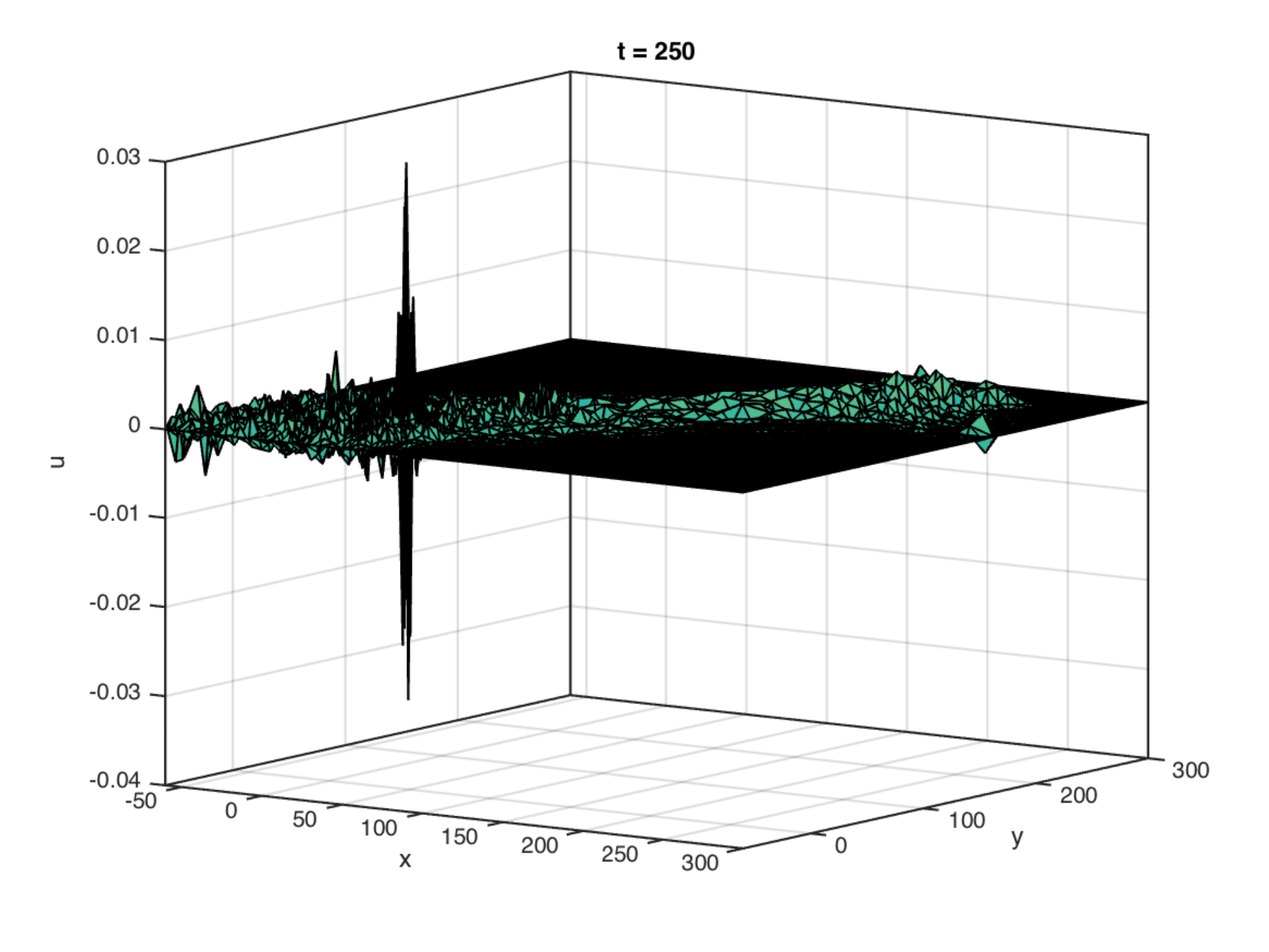}
\end{minipage}
\hspace{2mm}
\begin{minipage}[t]{2in}
\includegraphics[width=2.5in]{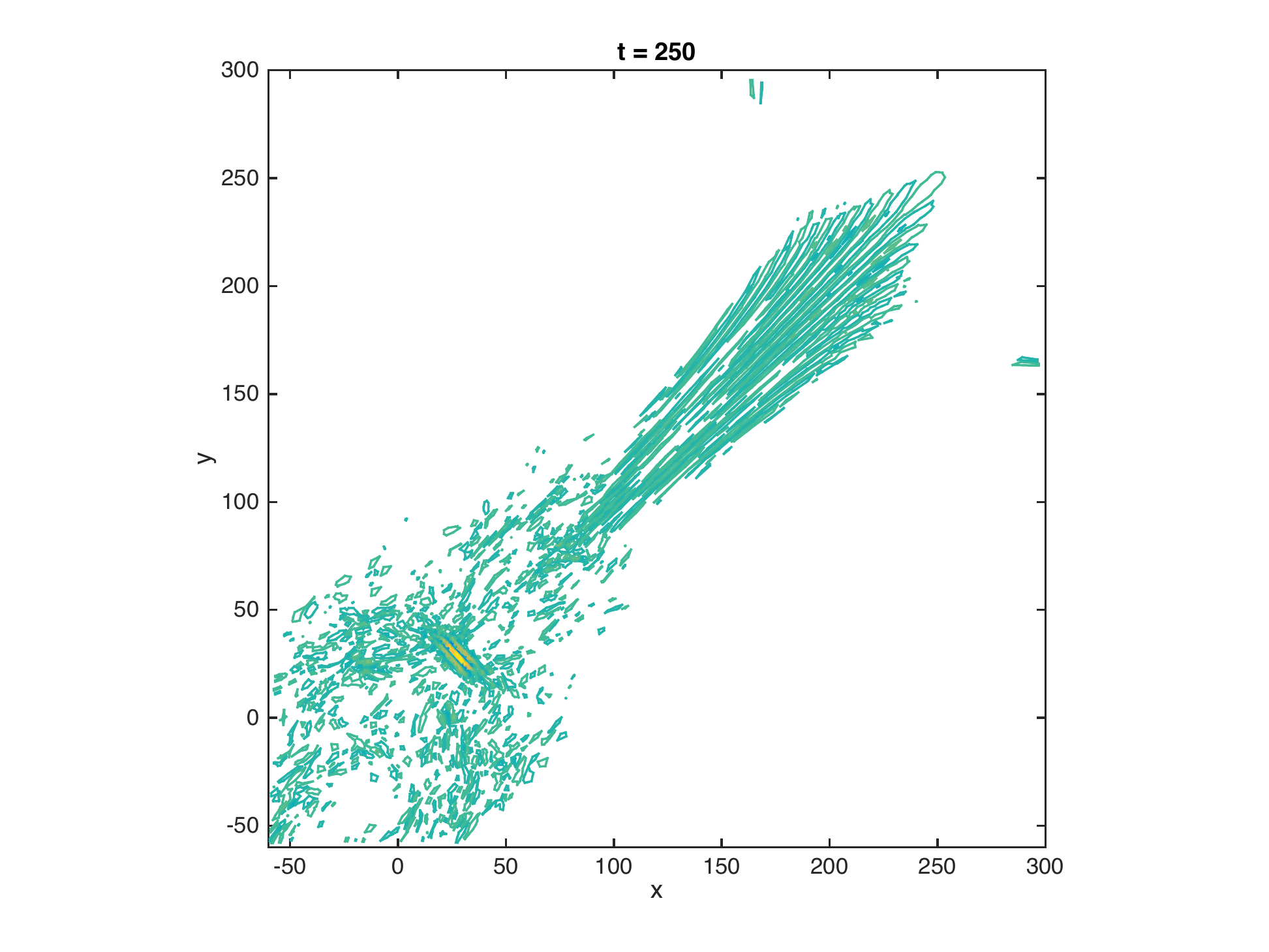}
\end{minipage}
\hspace{2mm}
\begin{minipage}[t]{2in}
\includegraphics[width=2.5in]{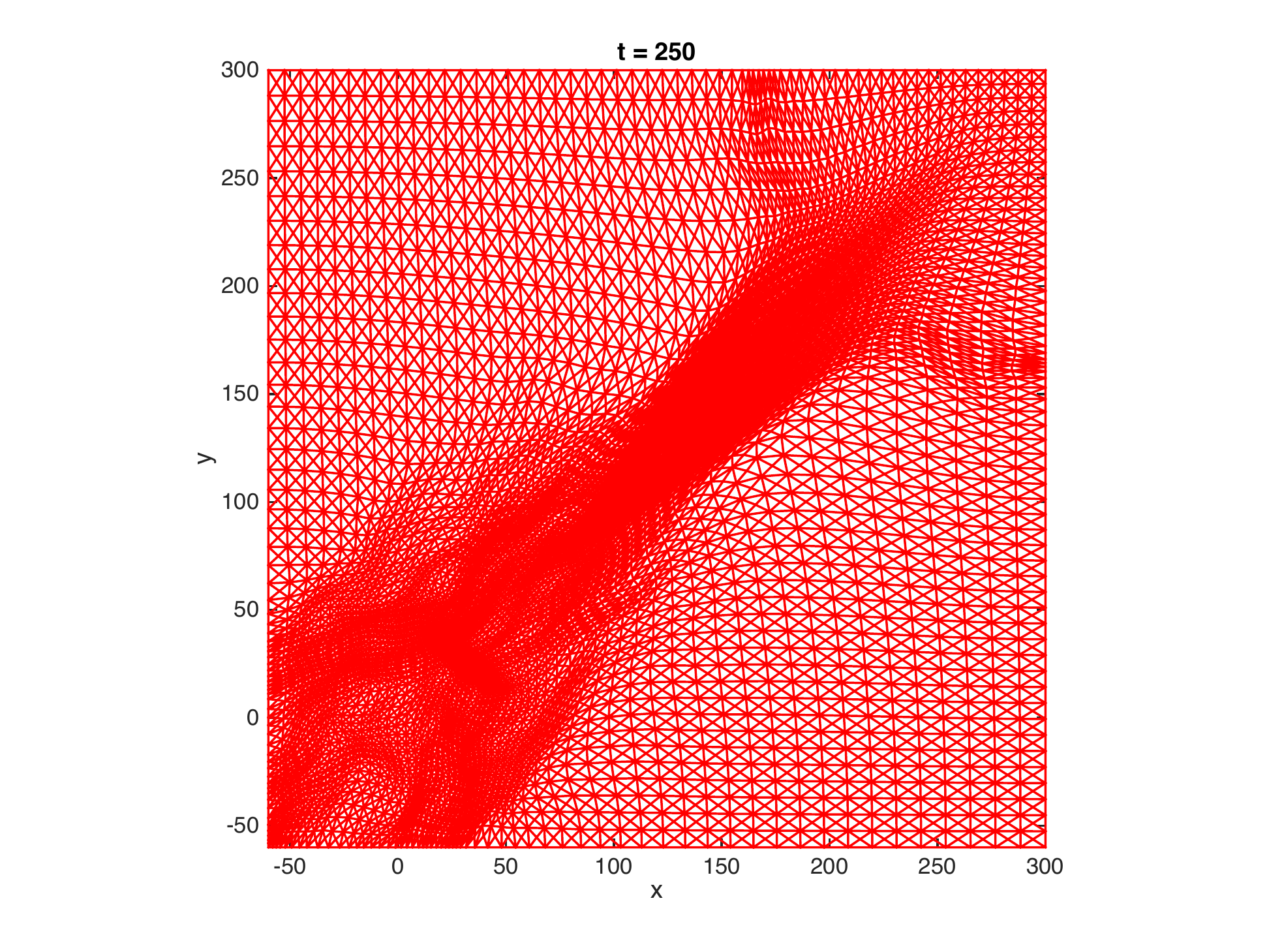}
\end{minipage}
\hspace{2mm}
}
\caption{Example~\ref{exam3.6}. Development of the 2D undular bore obtained with a moving mesh of $N=14400$.
The left column is for the numerical solution, the middle column is for the contours of the numerical solution,
and the right column is for the mesh.}
\label{exam-3.6-1}
\end{figure}

\begin{exam}
\label{exam3.6}
(2D RLW with undular bore)
This example is a two-dimensional generalization of Example~\ref{exam3.3} (the 1D undular bore).
The equation (\ref{RLW-1}) is subject to a homogeneous Dirichlet boundary condition and the initial condition
\[
u(x,y,0) = \frac{u_0}{2} \left (  1-\tanh\left ( (x-x_0)^2+(y-y_0)^2 - d^2\right ) \right ),
\]
where $\alpha = \beta = 1$, $\gamma = \delta = 1.5$, $\mu = 1/6$, $u_0 = 0.1$,
$x_0 = y_0 = 0$, and $d = 2$.
The computation is performed on $\Omega = (-60,300)\times(-60,300)$ with $T = 250$.

Fig.~\ref{exam-3.6-1} shows the development and expansion of the 2D undular bore
which propagates in a northeast direction. Compared to the 1D situation,
the propagation is slightly slower and the amplitude is smaller.
The mesh concentration correctly reflects the development of the undular bore.
\end{exam}

\begin{figure}[thb]
\centering
\hbox{
\begin{minipage}[t]{2in}
\includegraphics[width=2.3in]{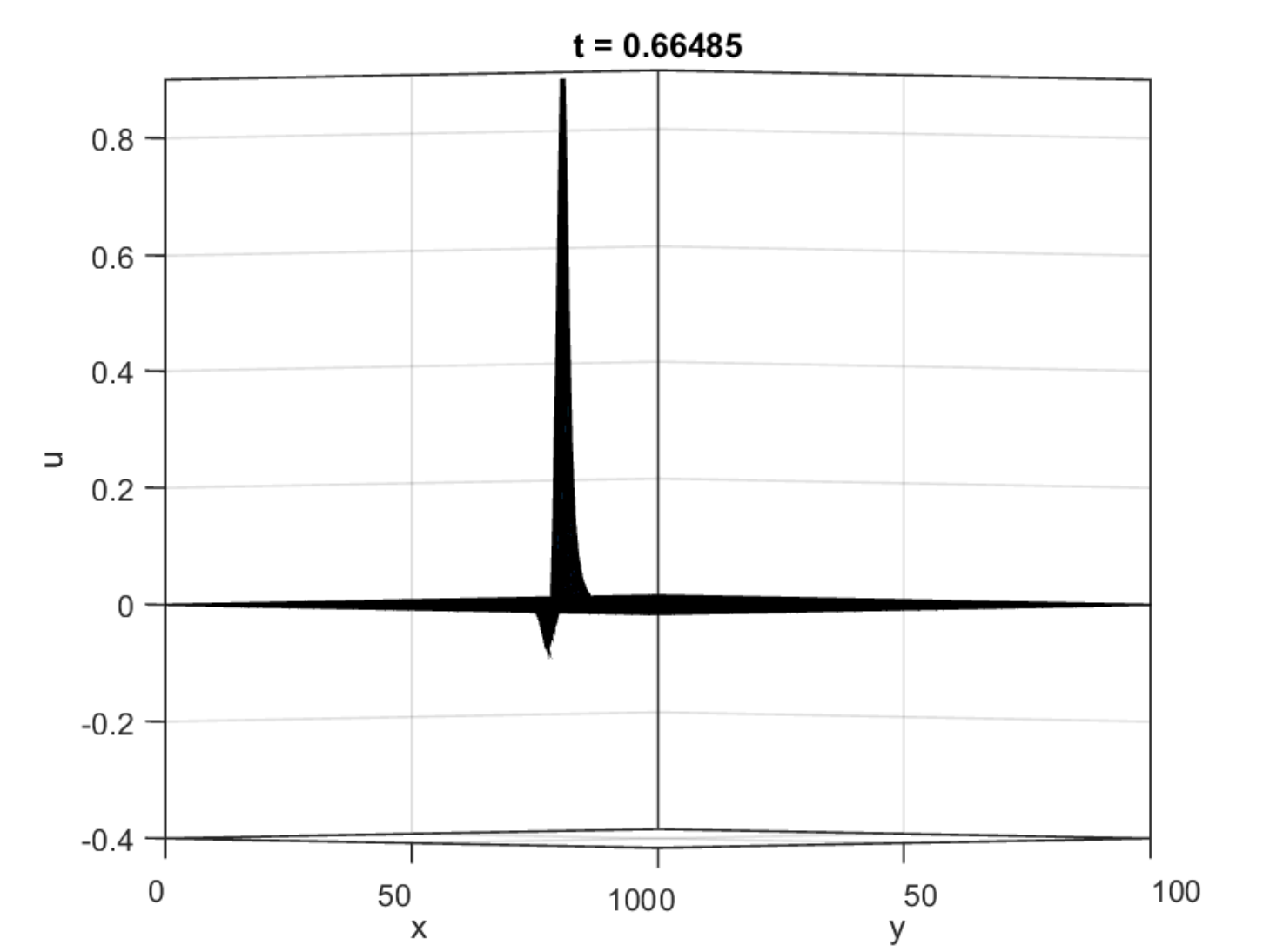}
\end{minipage}
\hspace{2mm}
\begin{minipage}[t]{2in}
\includegraphics[width=2.5in]{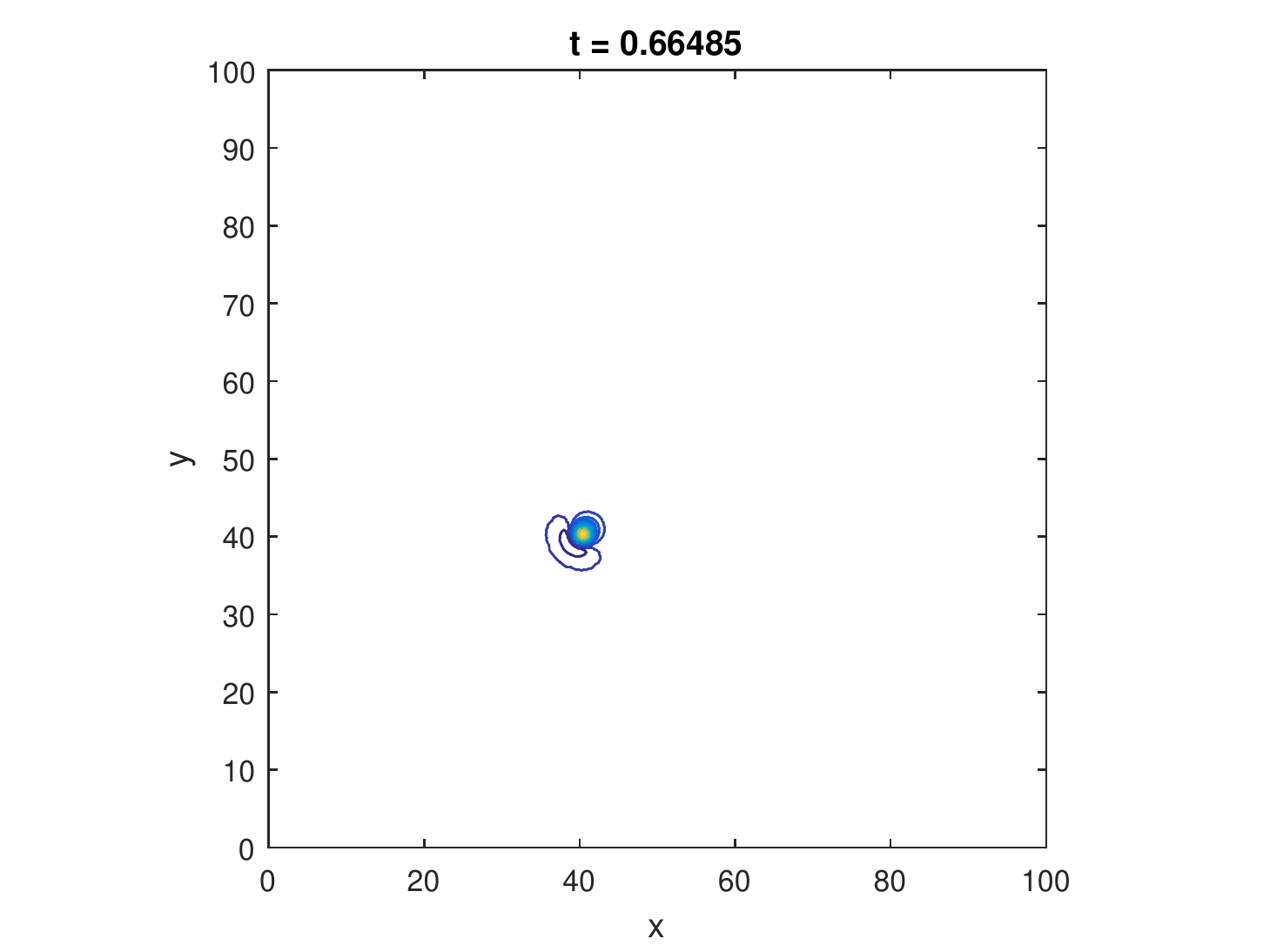}
\end{minipage}
\hspace{2mm}
\begin{minipage}[t]{2in}
\includegraphics[width=2.5in]{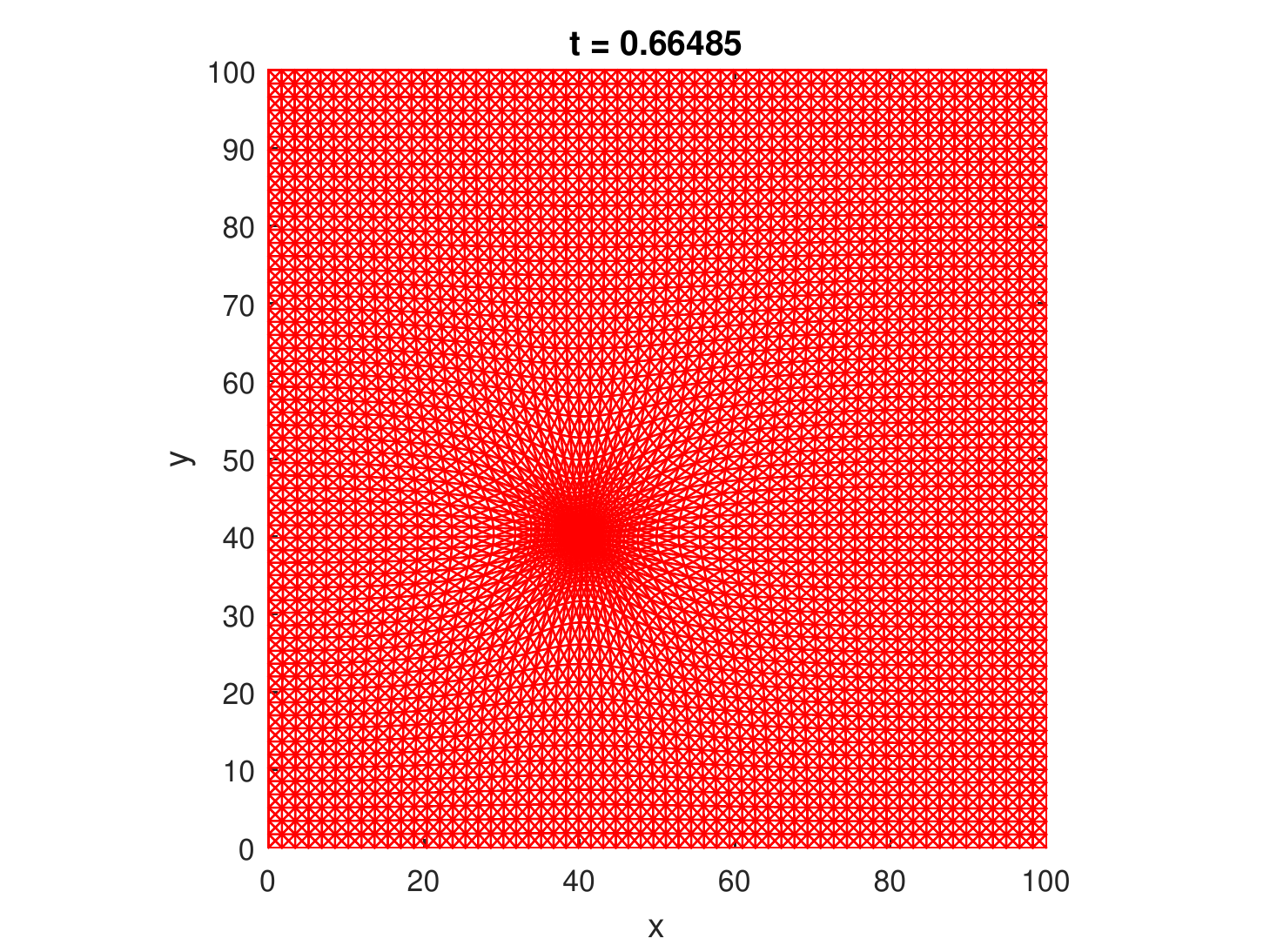}
\end{minipage}
\hspace{2mm}
}
\hbox{
\begin{minipage}[t]{2in}
\includegraphics[width=2.3in]{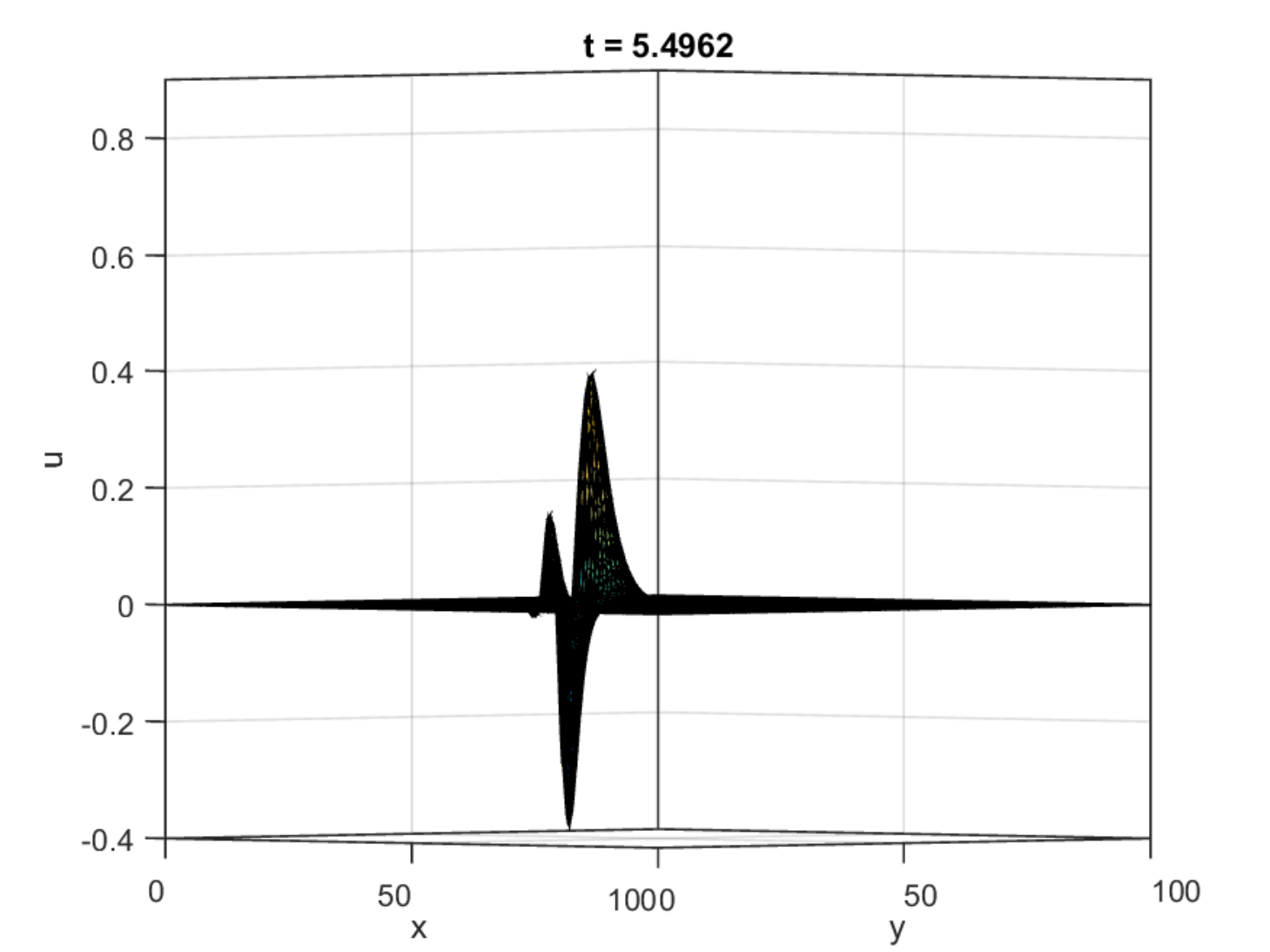}
\end{minipage}
\hspace{2mm}
\begin{minipage}[t]{2in}
\includegraphics[width=2.5in]{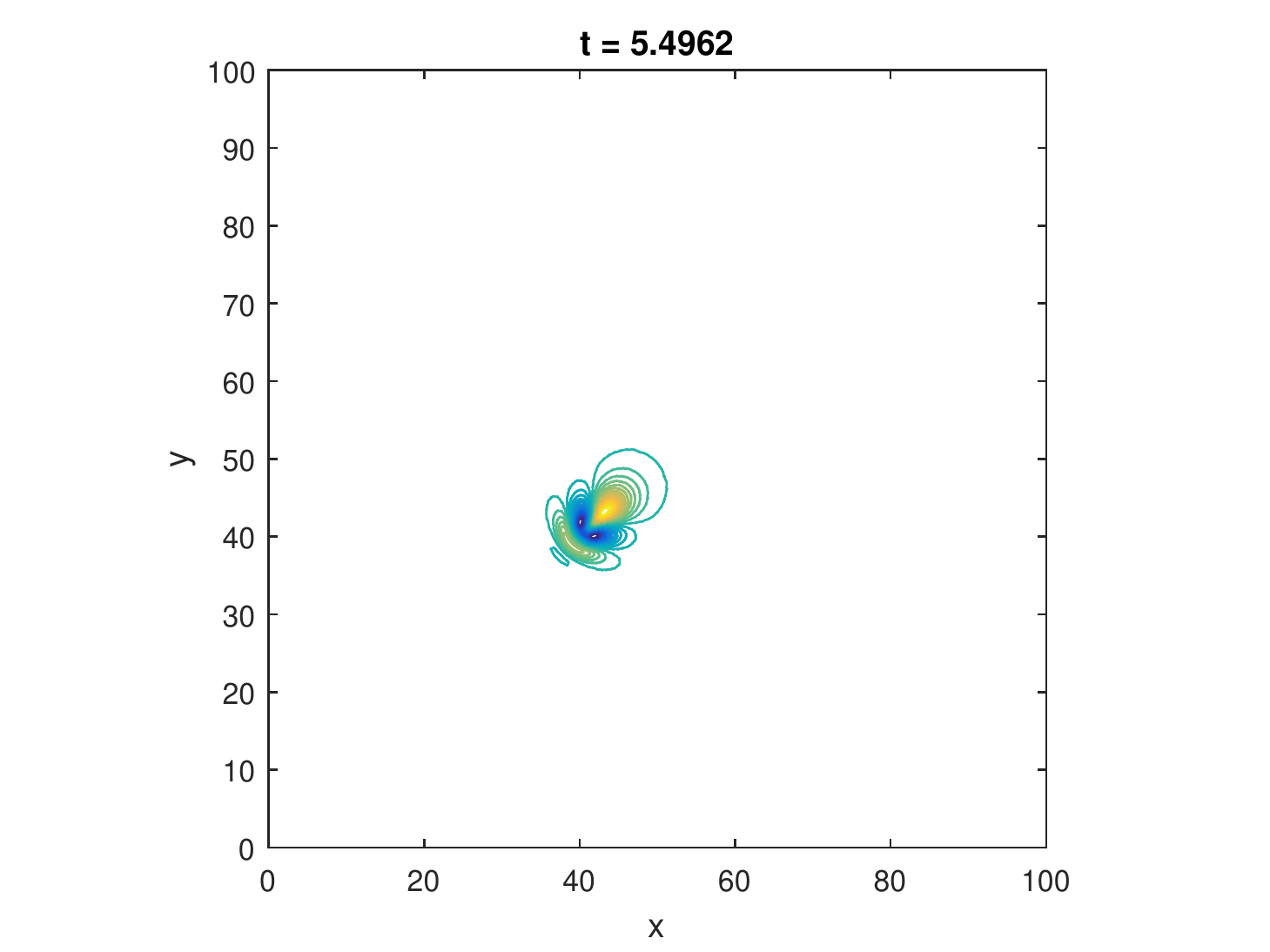}
\end{minipage}
\hspace{2mm}
\begin{minipage}[t]{2in}
\includegraphics[width=2.5in]{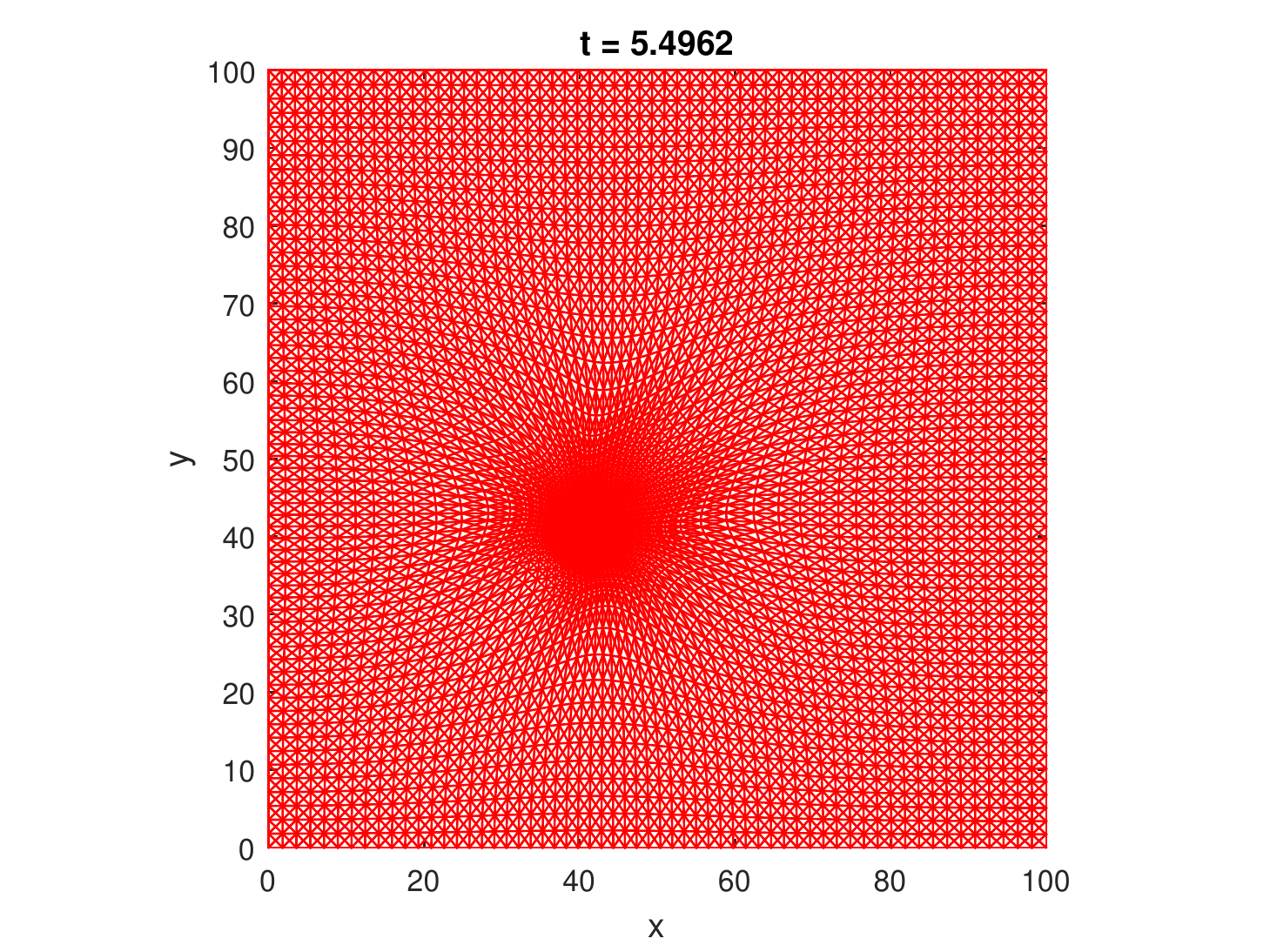}
\end{minipage}
\hspace{2mm}
}
\hbox{
\begin{minipage}[t]{2in}
\includegraphics[width=2.3in]{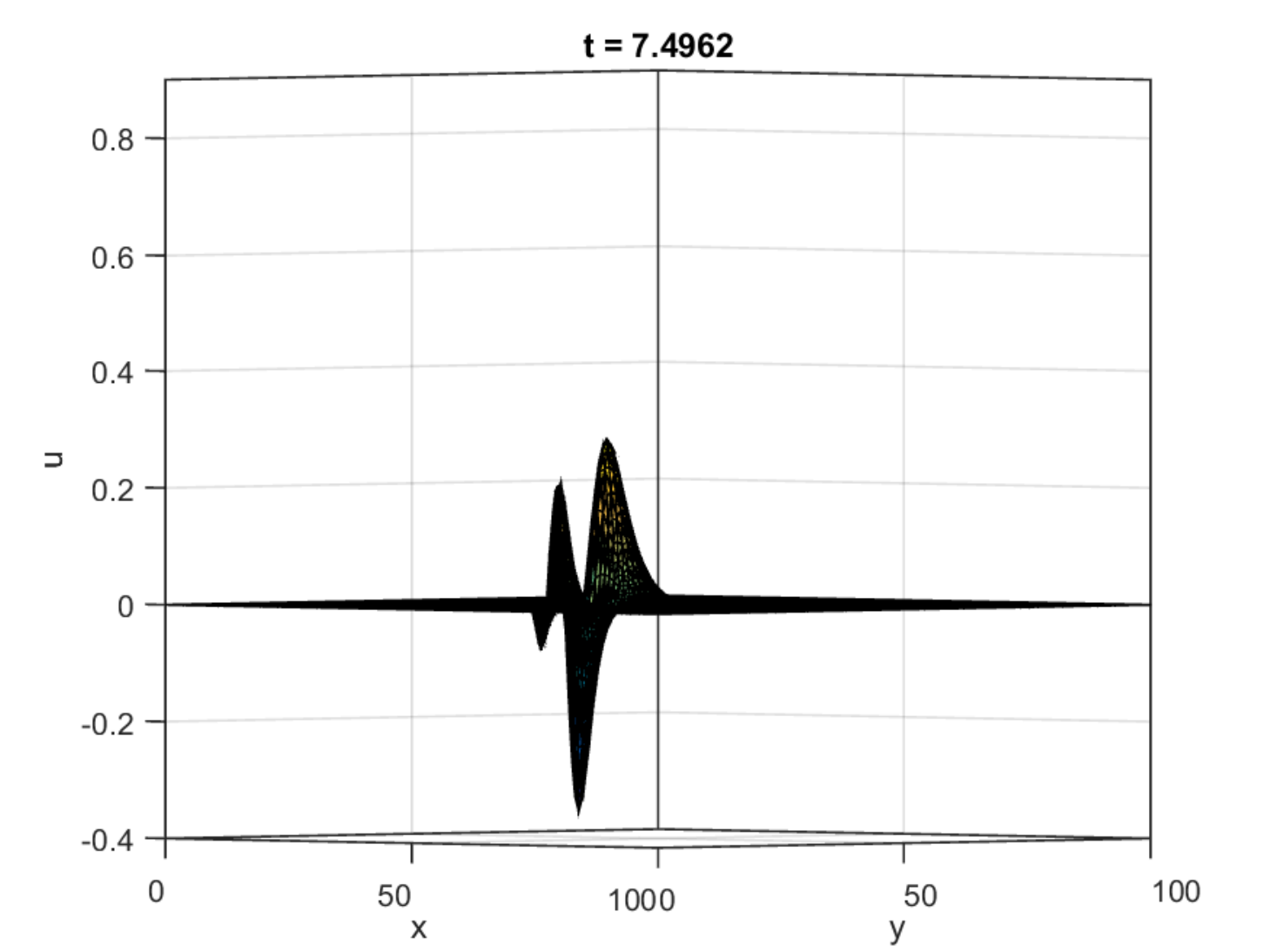}
\end{minipage}
\hspace{2mm}
\begin{minipage}[t]{2in}
\includegraphics[width=2.5in]{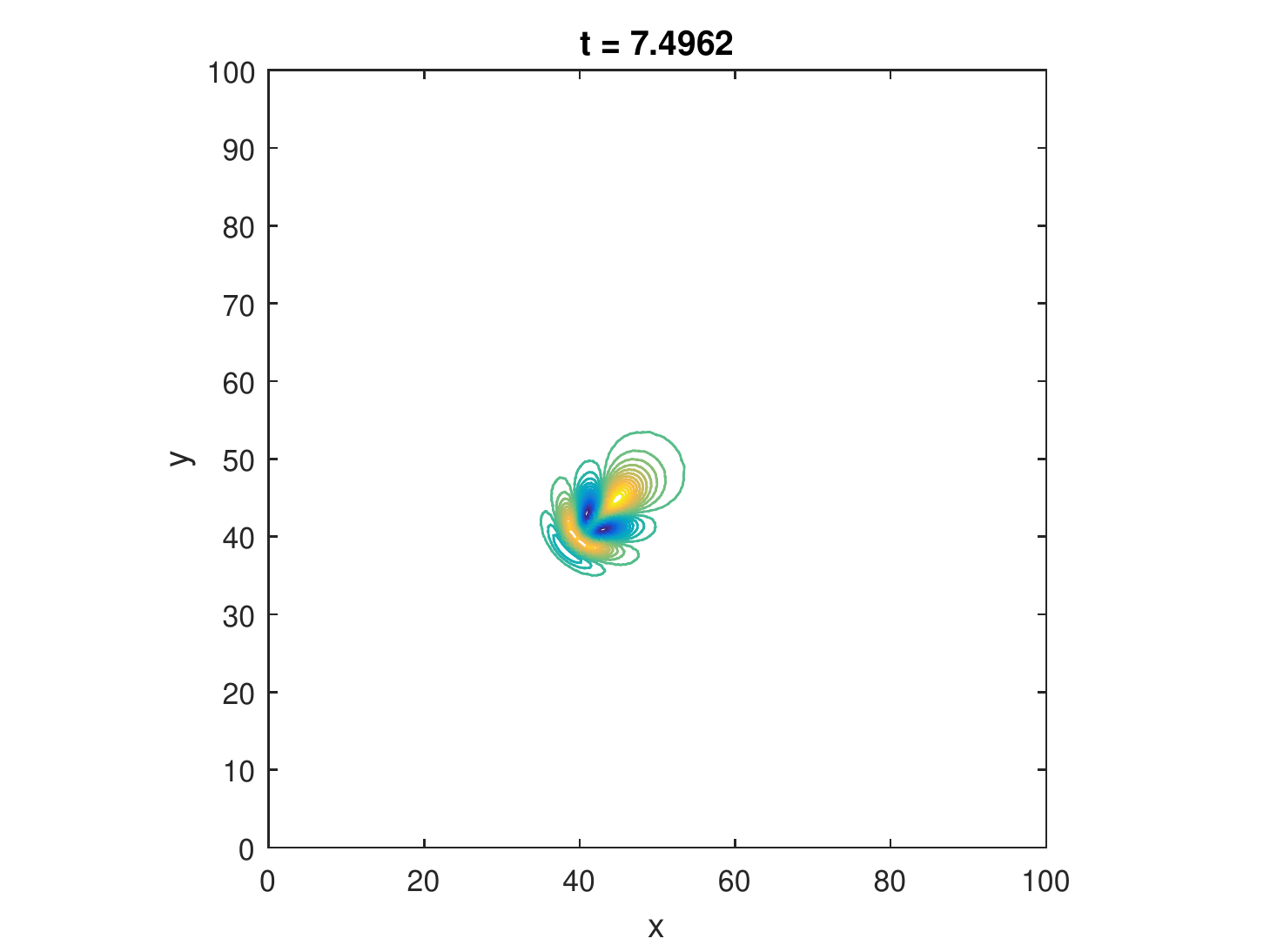}
\end{minipage}
\hspace{2mm}
\begin{minipage}[t]{2in}
\includegraphics[width=2.5in]{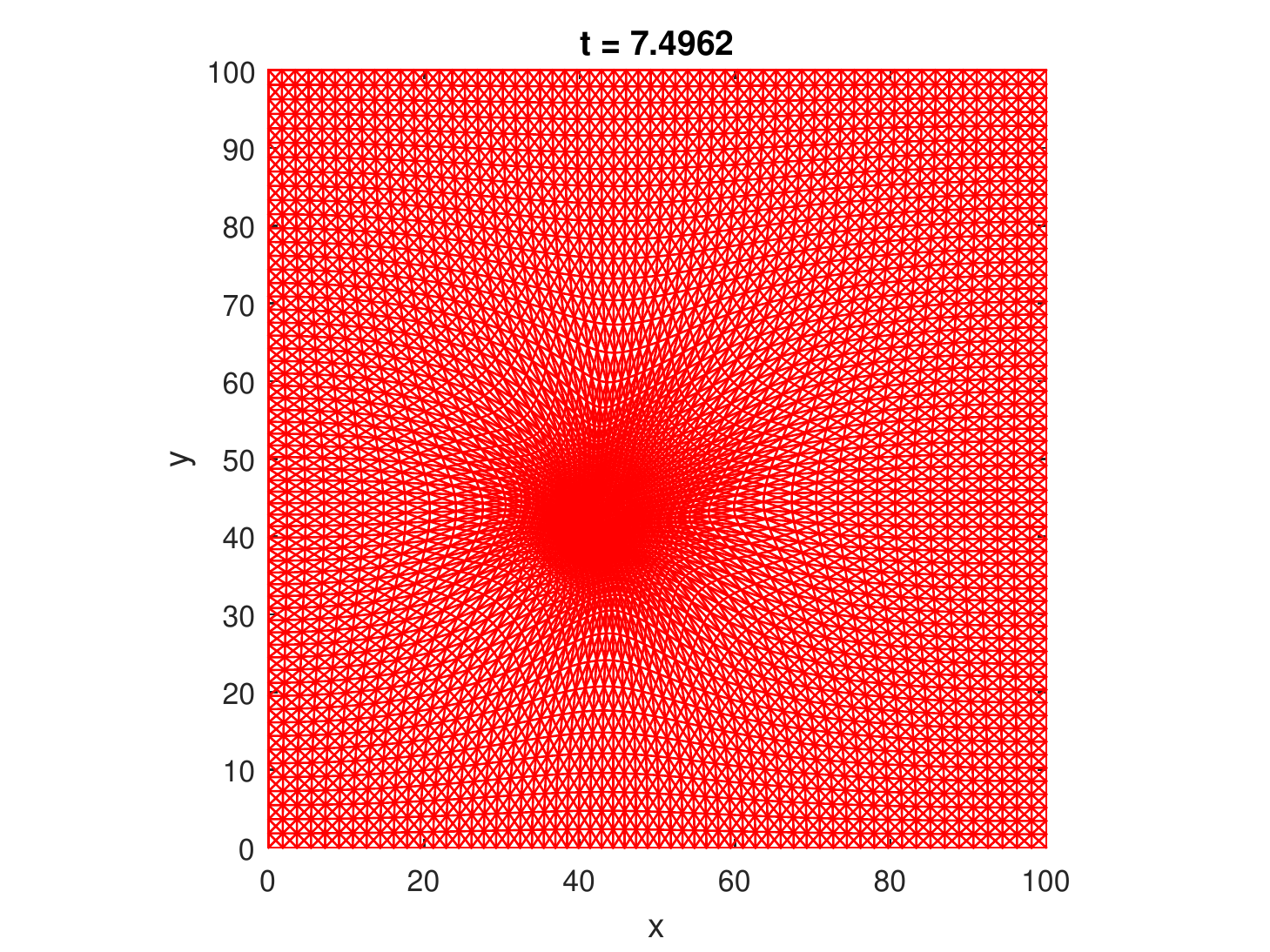}
\end{minipage}
\hspace{2mm}
}
\hbox{
\begin{minipage}[t]{2in}
\includegraphics[width=2.3in]{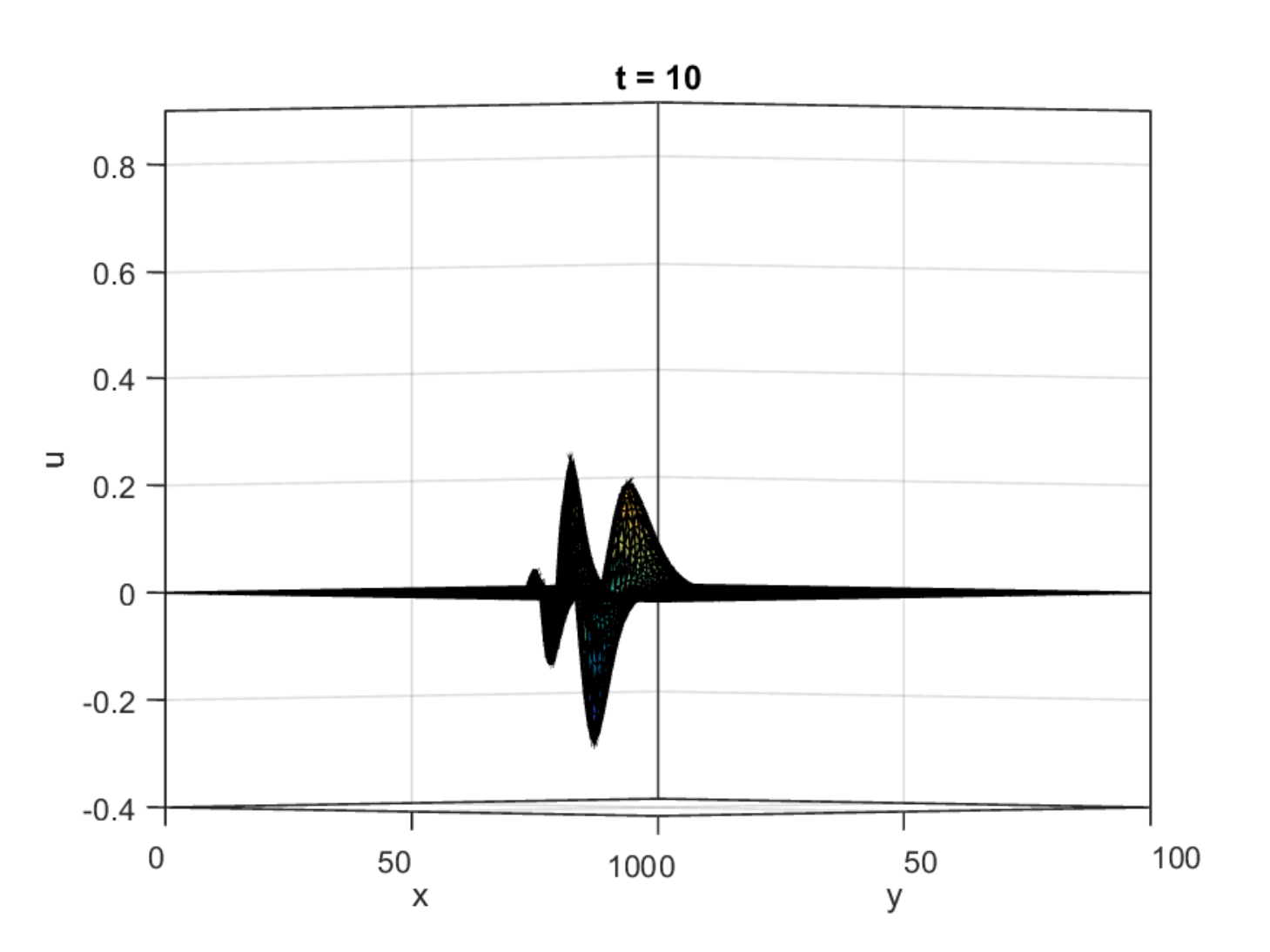}
\end{minipage}
\hspace{2mm}
\begin{minipage}[t]{2in}
\includegraphics[width=2.5in]{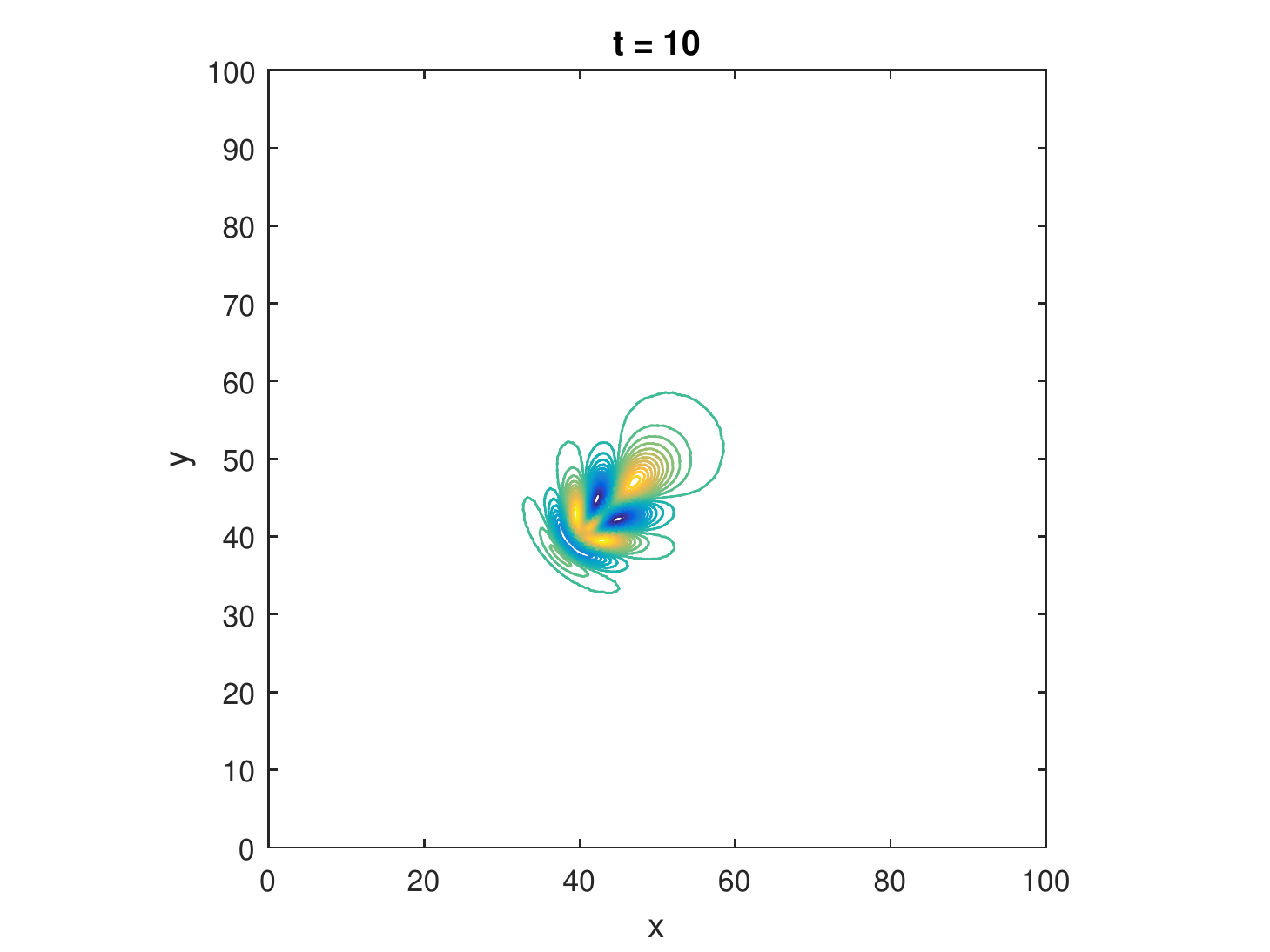}
\end{minipage}
\hspace{2mm}
\begin{minipage}[t]{2in}
\includegraphics[width=2.5in]{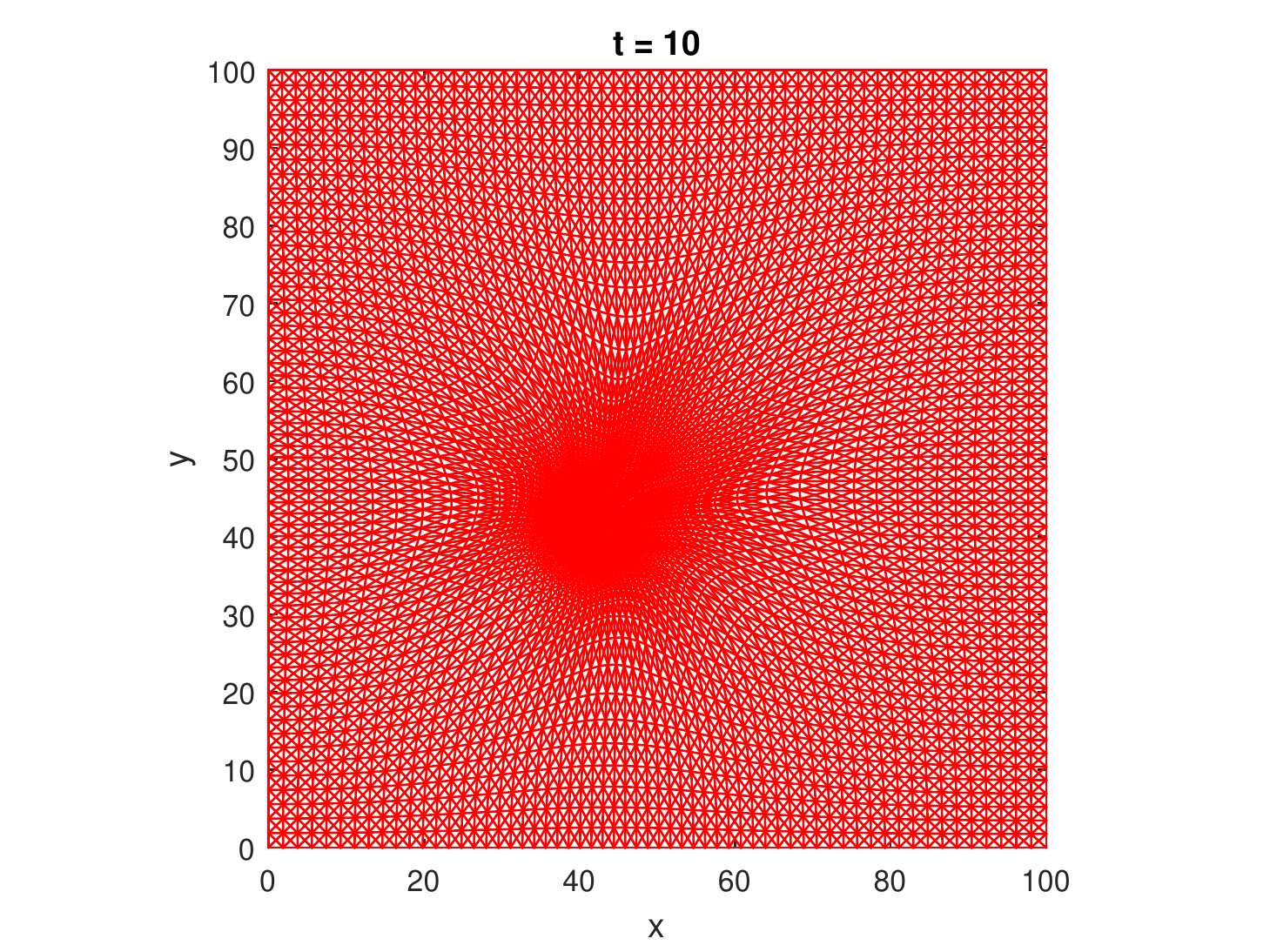}
\end{minipage}
\hspace{2mm}
}
\caption{Example~\ref{exam3.7}. The numerical solution, its contours, and the mesh
are shown at various time instants for the 2D Maxwellian initial condition case with $\mu = 1$.
A moving mesh of  $N = 14400$ is used.
}
\label{exam-3.7-1}
\end{figure}

\begin{figure}[thb]
\centering
\hbox{
\begin{minipage}[t]{2in}
\includegraphics[width=2.3in]{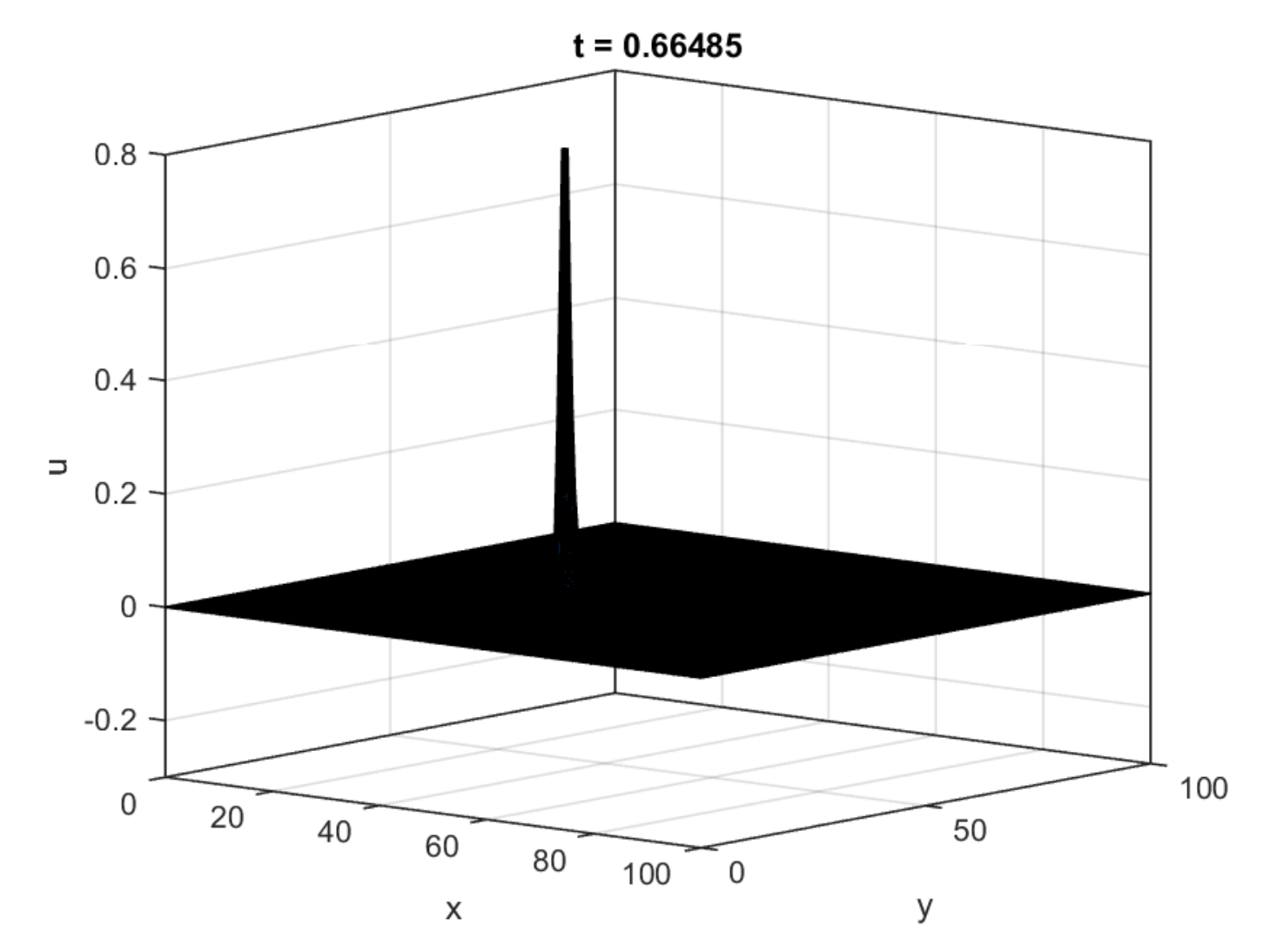}
\end{minipage}
\hspace{2mm}
\begin{minipage}[t]{2in}
\includegraphics[width=2.5in]{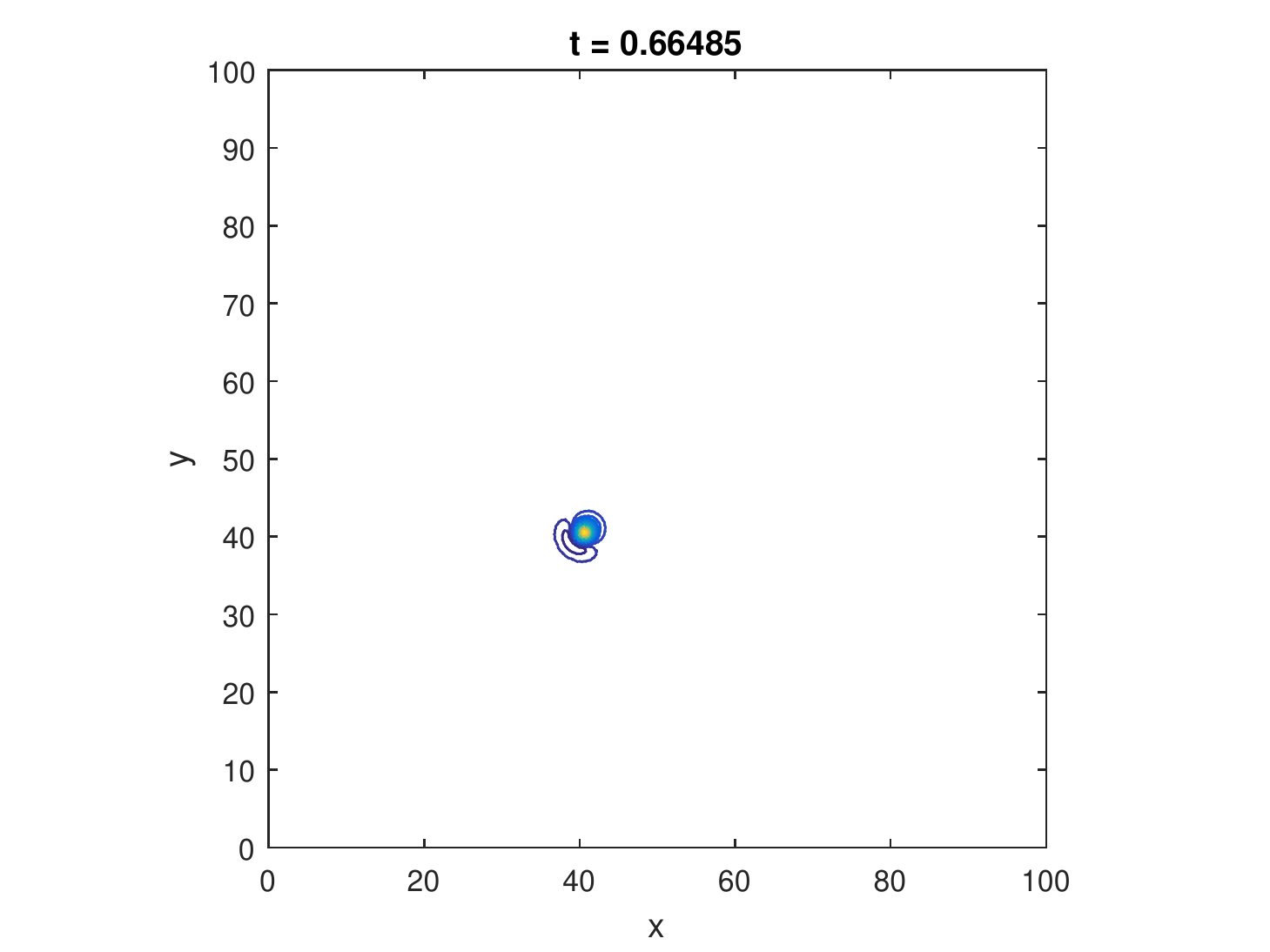}
\end{minipage}
\hspace{2mm}
\begin{minipage}[t]{2in}
\includegraphics[width=2.5in]{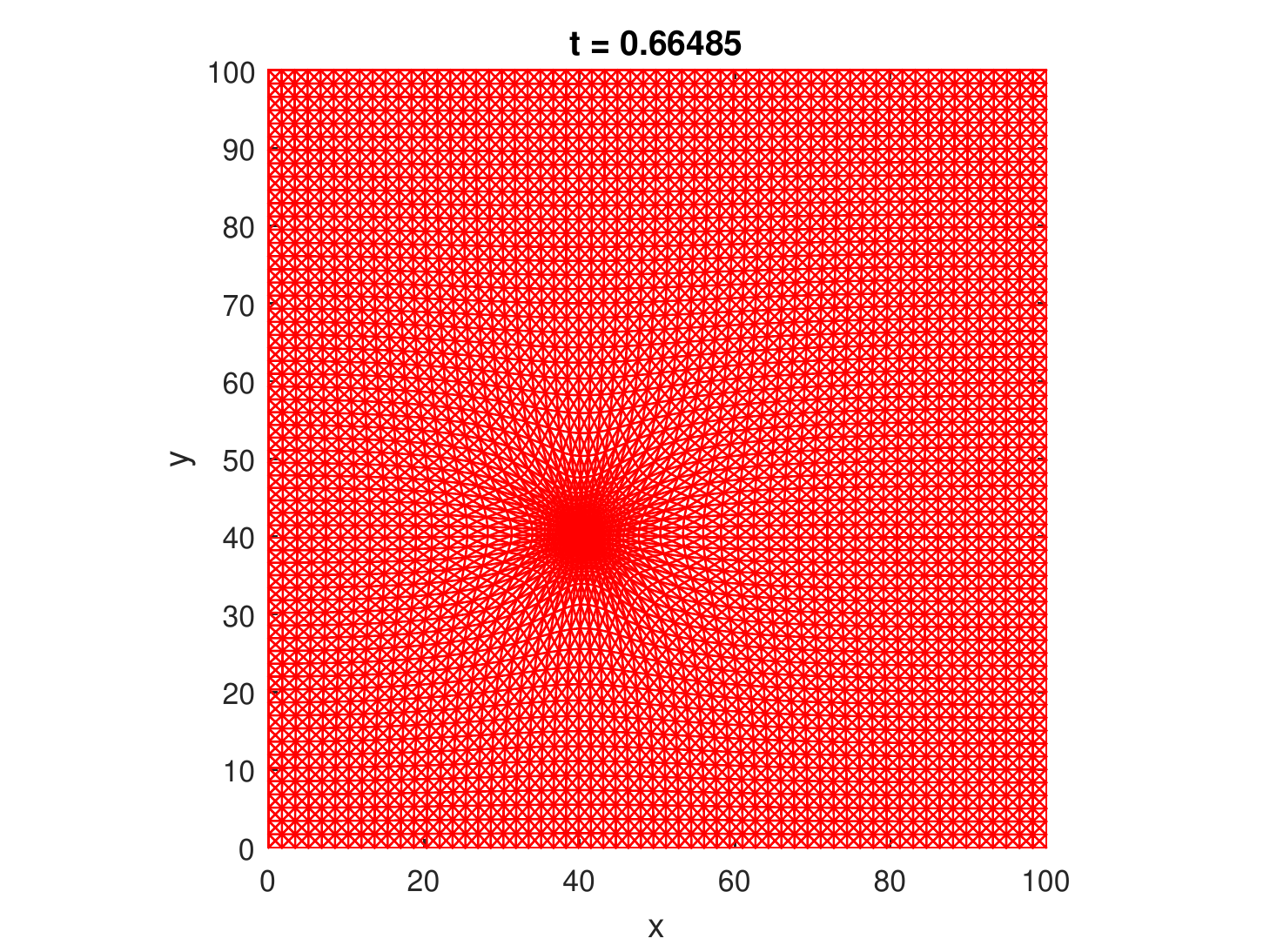}
\end{minipage}
\hspace{2mm}
}
\hbox{
\begin{minipage}[t]{2in}
\includegraphics[width=2.3in]{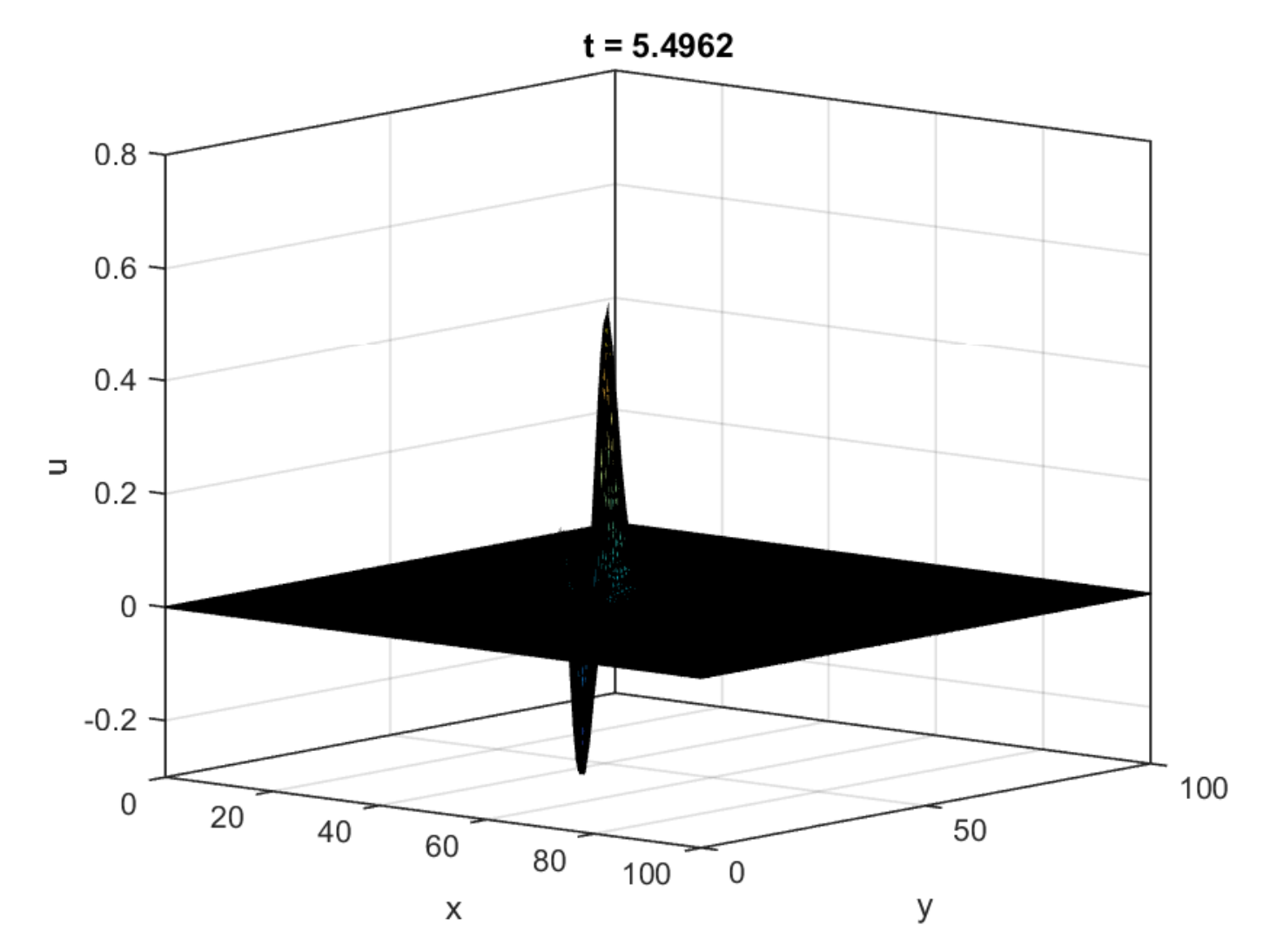}
\end{minipage}
\hspace{2mm}
\begin{minipage}[t]{2in}
\includegraphics[width=2.5in]{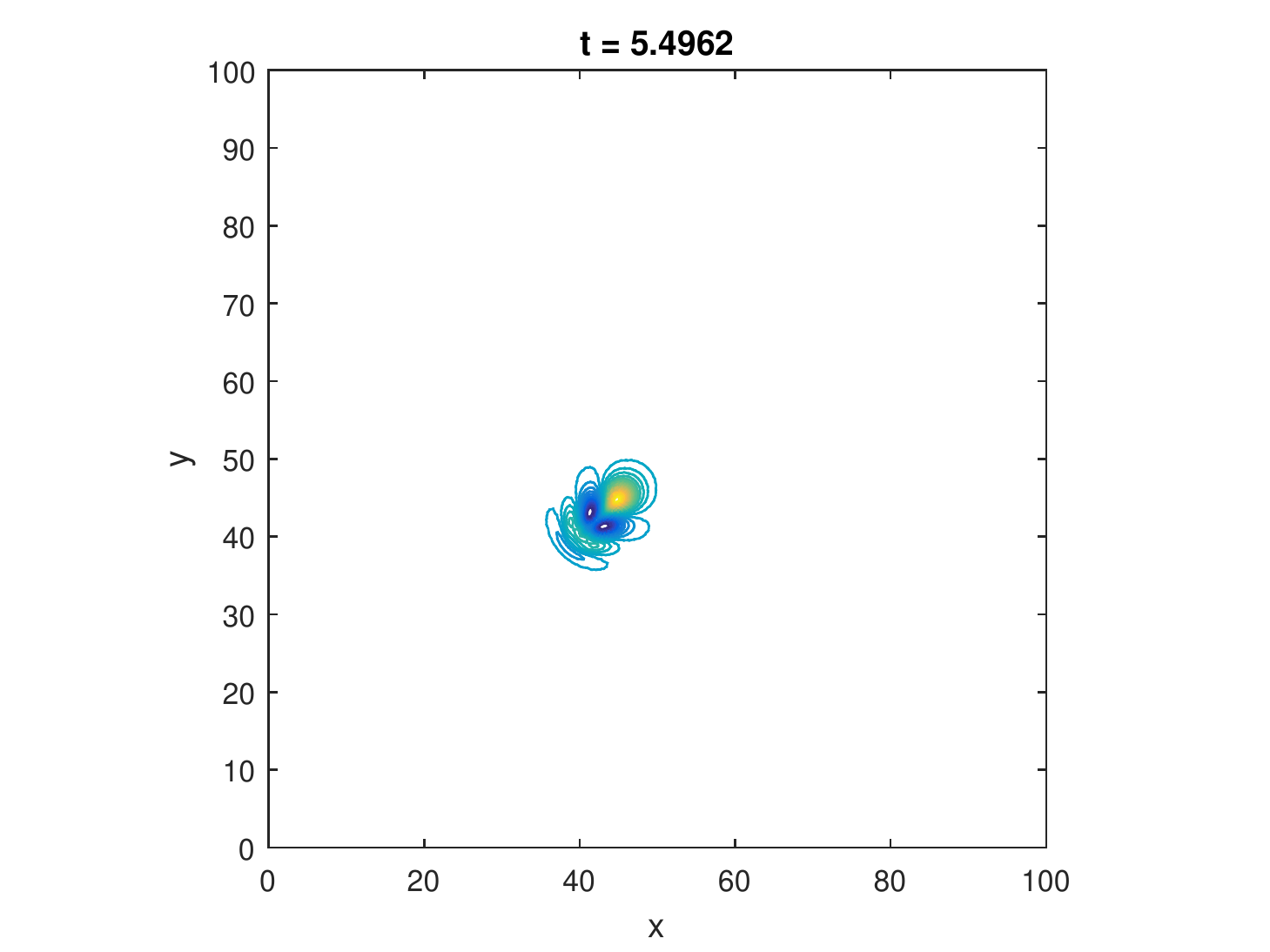}
\end{minipage}
\hspace{2mm}
\begin{minipage}[t]{2in}
\includegraphics[width=2.5in]{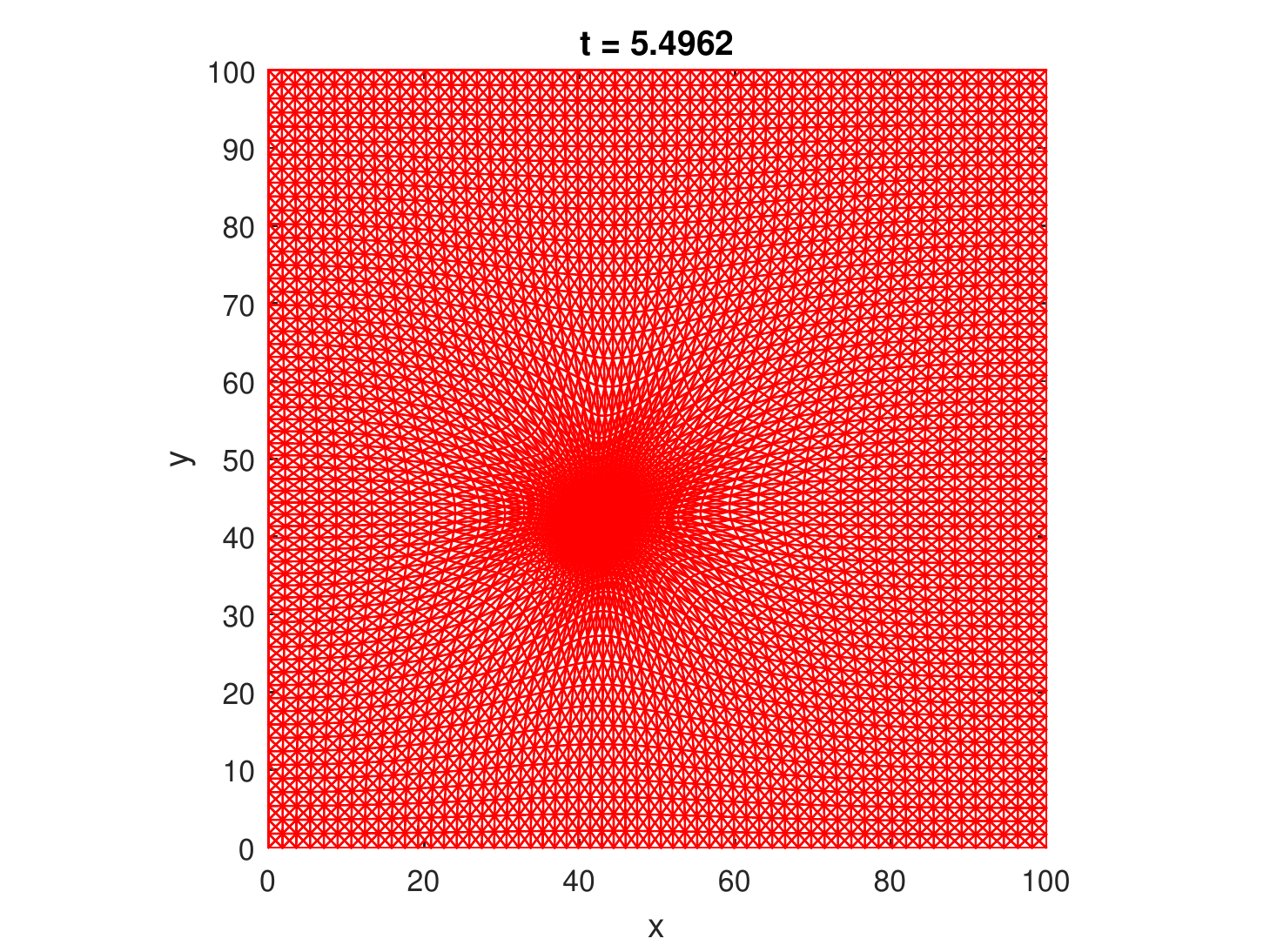}
\end{minipage}
\hspace{2mm}
}
\hbox{
\begin{minipage}[t]{2in}
\includegraphics[width=2.3in]{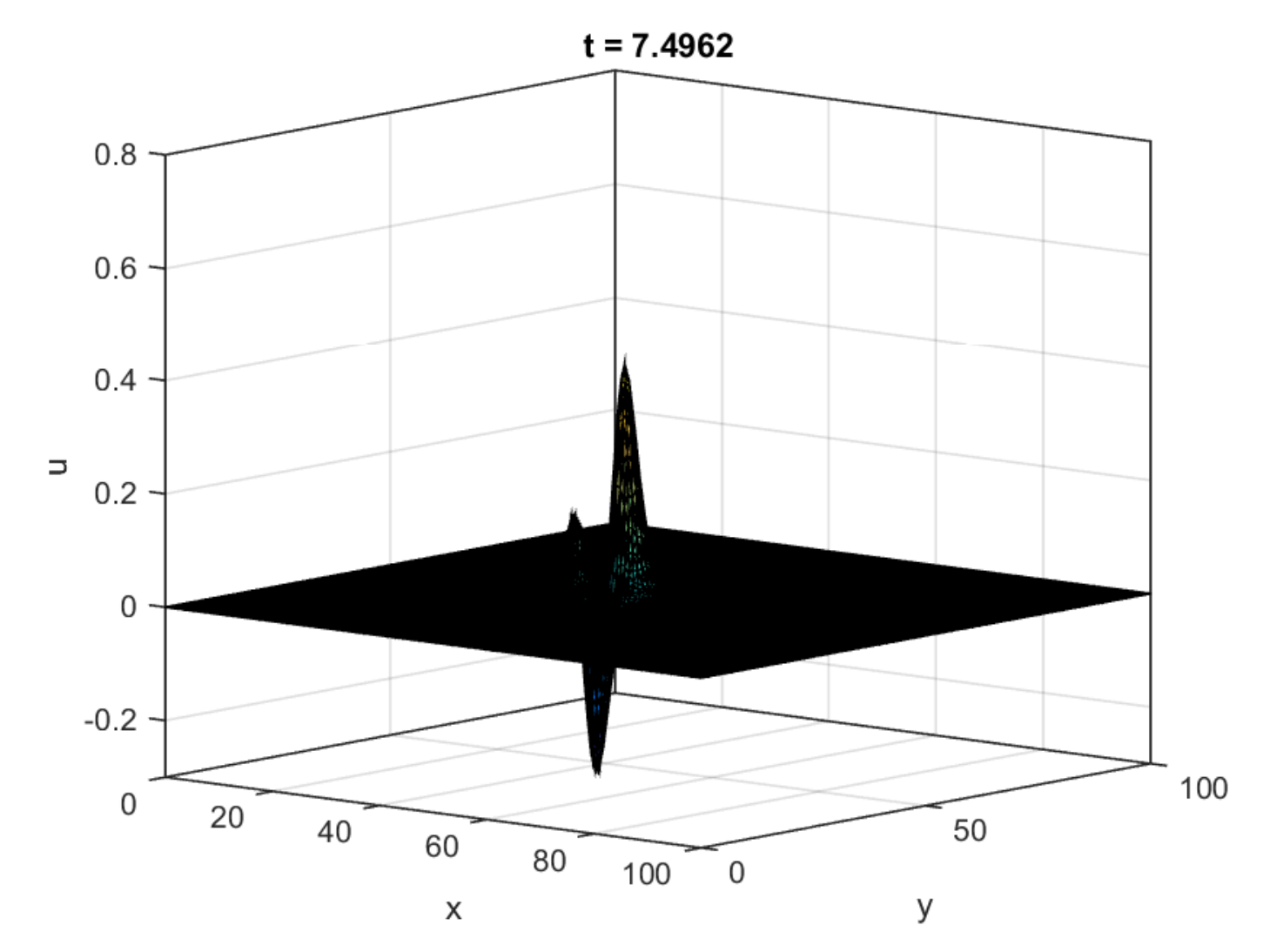}
\end{minipage}
\hspace{2mm}
\begin{minipage}[t]{2in}
\includegraphics[width=2.5in]{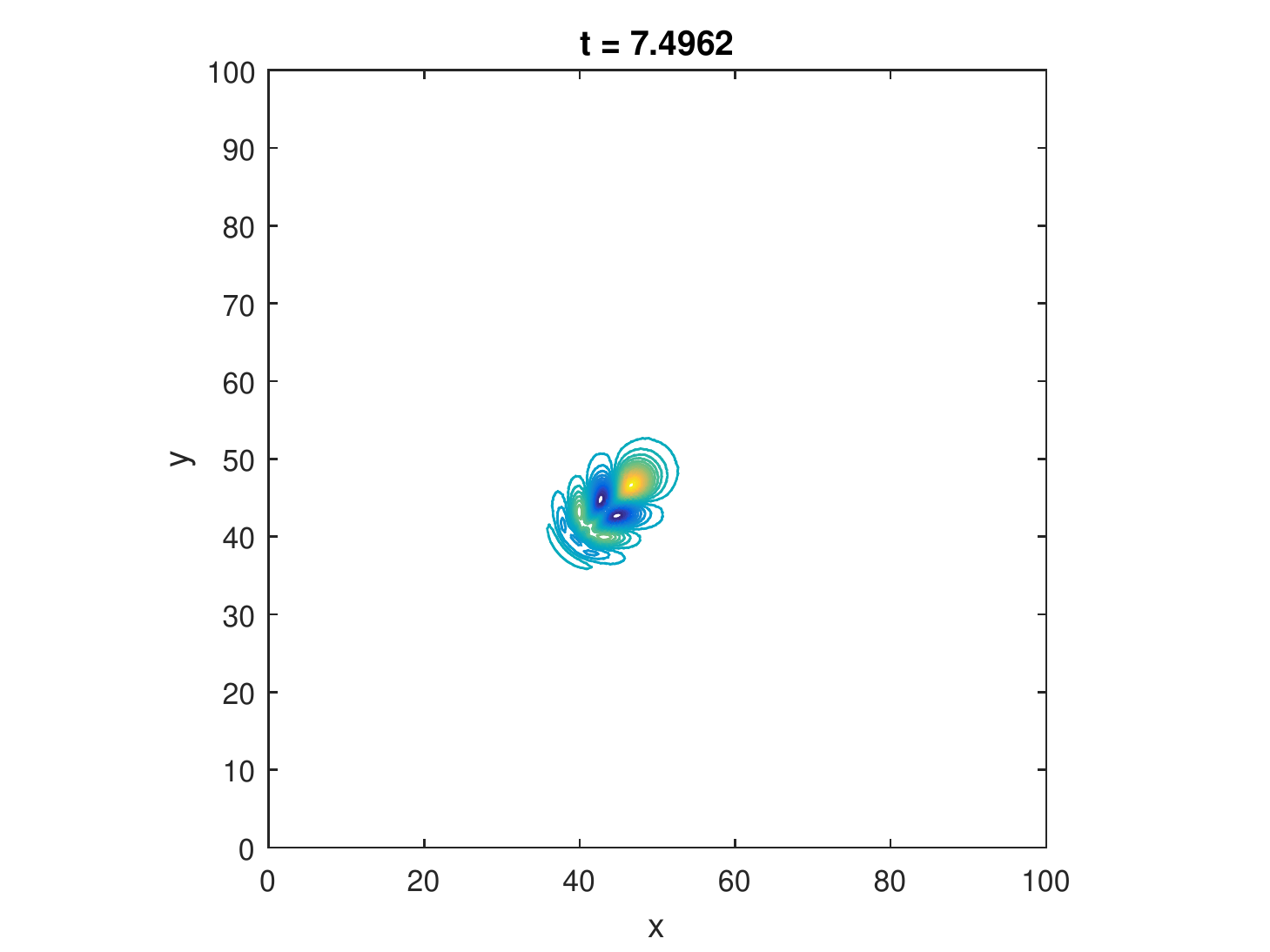}
\end{minipage}
\hspace{2mm}
\begin{minipage}[t]{2in}
\includegraphics[width=2.5in]{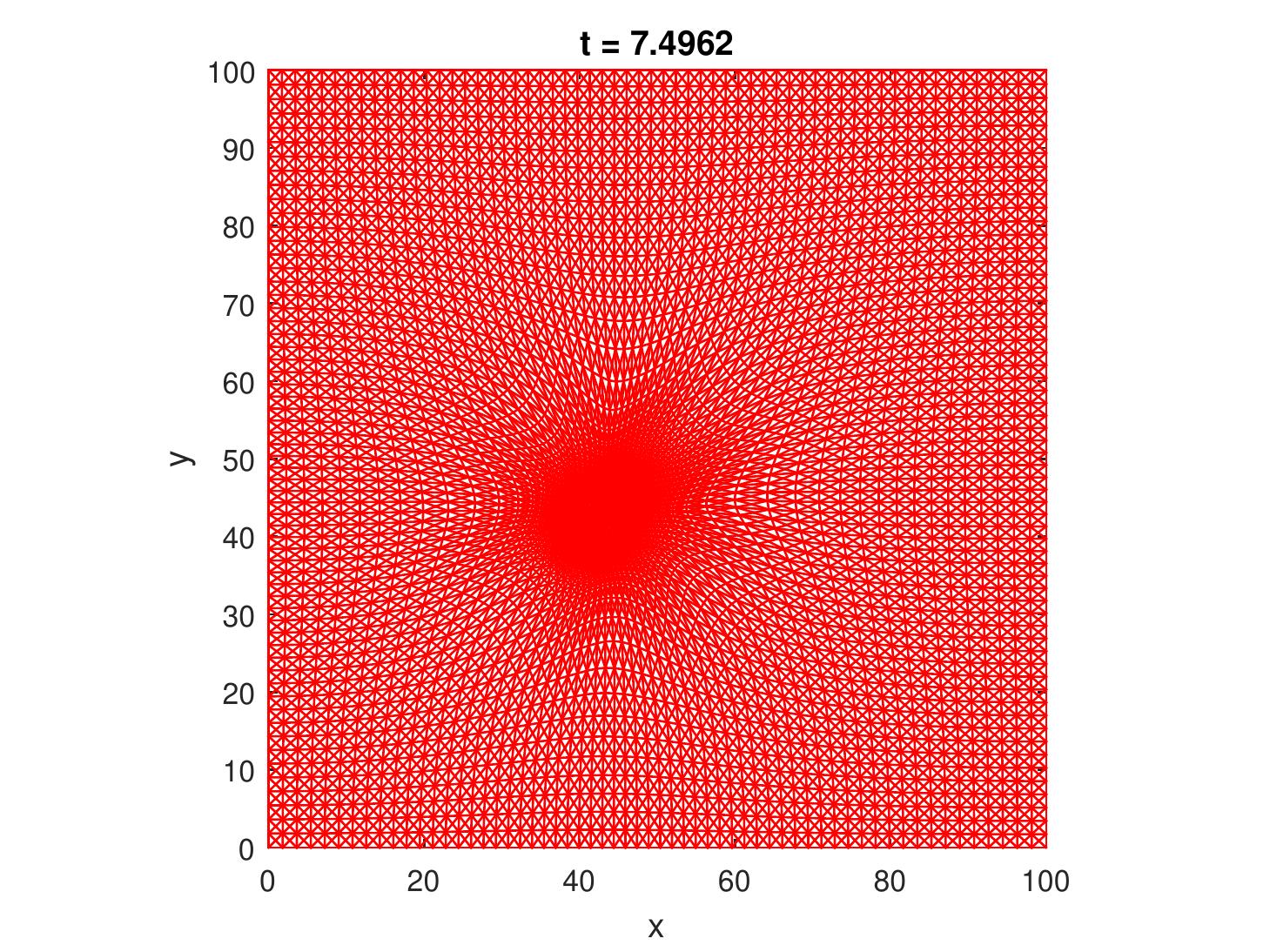}
\end{minipage}
\hspace{2mm}
}
\hbox{
\begin{minipage}[t]{2in}
\includegraphics[width=2.5in]{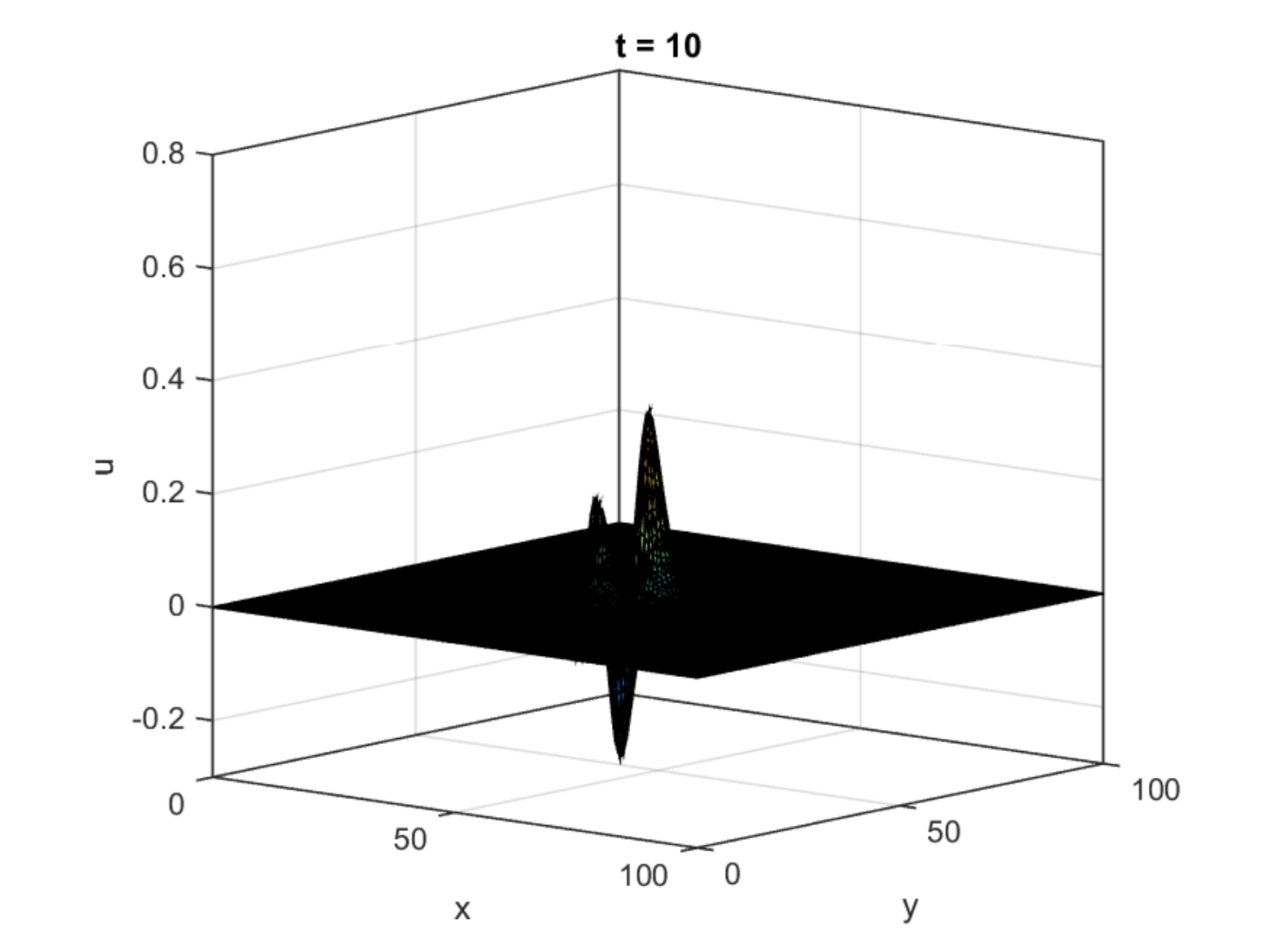}
\end{minipage}
\hspace{2mm}
\begin{minipage}[t]{2in}
\includegraphics[width=2.5in]{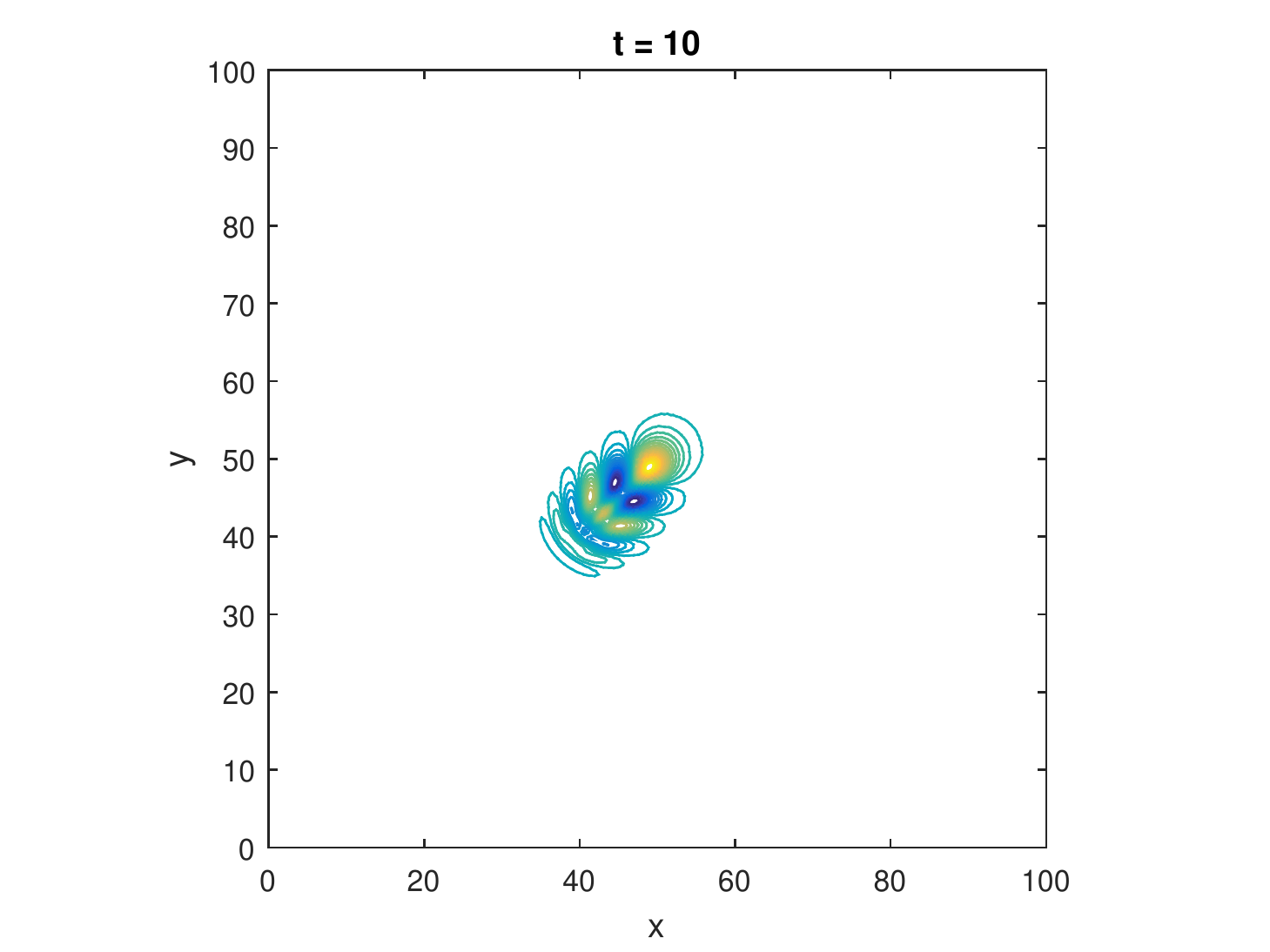}
\end{minipage}
\hspace{2mm}
\begin{minipage}[t]{2in}
\includegraphics[width=2.5in]{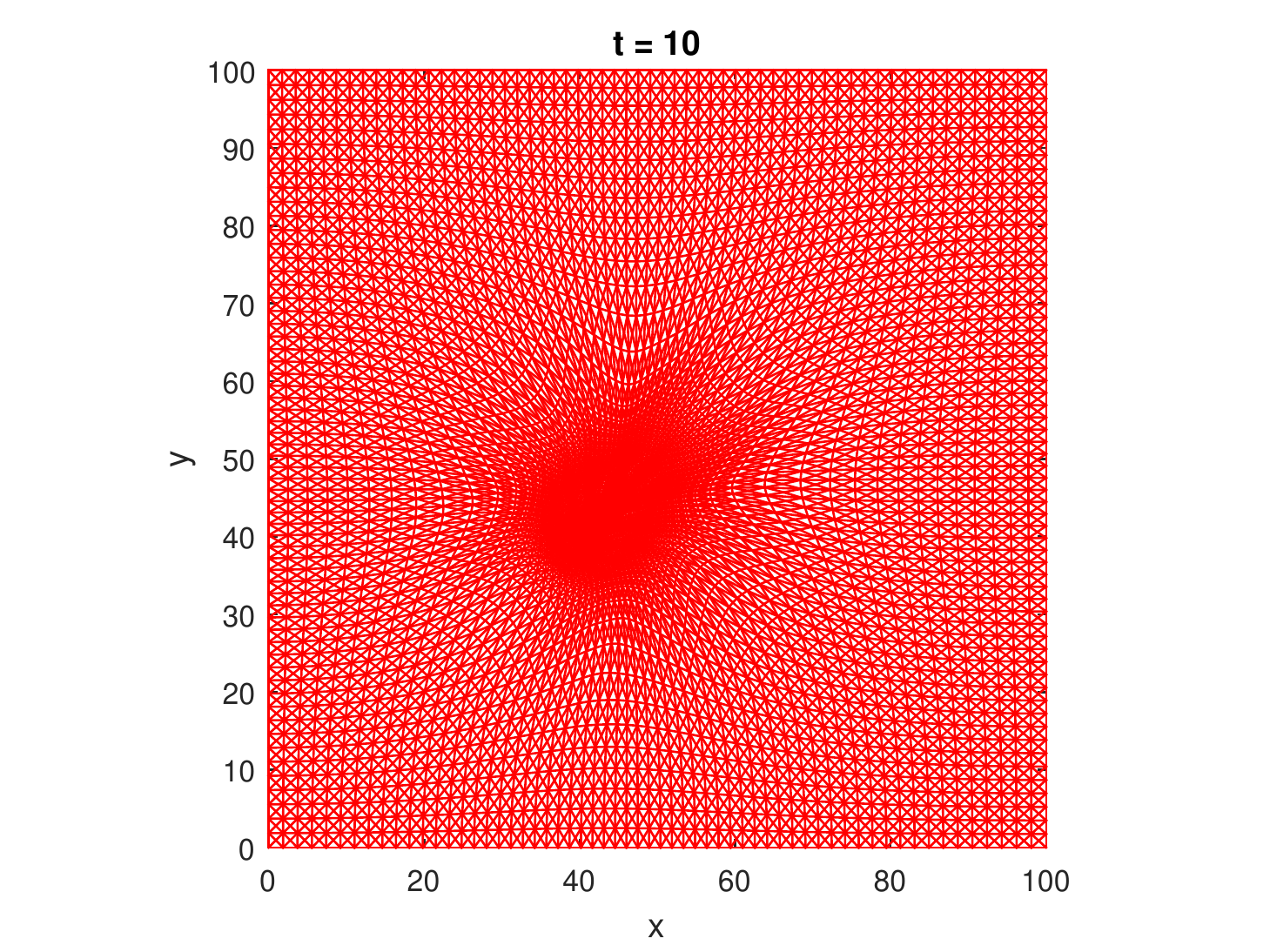}
\end{minipage}
\hspace{2mm}
}
\caption{Example~\ref{exam3.7}. The numerical solution, its contours, and the mesh
are shown at various time instants for the 2D Maxwellian initial condition case with $\mu = 0.5$.
A moving mesh of  $N = 14400$ is used.
}
\label{exam-3.7-2}
\end{figure}

\begin{exam}
\label{exam3.7}
(2D RLW with the Maxwellian initial condition)
In this final example, we consider the initial Maxwellian initial condition
\[
u(x,y,0) = e^{-((x-40)^2+(y-40)^2)}
\]
for the 2D MRLW equation 
\[
u_t + u_x + u_y + \gamma u^2 u_x + \delta u^2 u_y - \mu u_{xxt} - \mu u_{yyt}  = 0,
\]
where 
$\gamma = \delta = 6$, and $\mu = 0.5$ or $\mu = 1$.
A homogeneous Dirichlet boundary condition is used.
The computation is performed on $\Omega = (0,100) \times (0,100)$ with $T = 10$.

The numerical results are shown in Fig.~\ref{exam-3.7-1} for $\mu = 1$ and
Fig.~\ref{exam-3.7-2} for $\mu = 0.5$. It can be seen that the train
of solitary waves is developed mainly along the northeast direction. Moreover,
it is obvious that the mesh elements are concentrated in the peak region of the solitary waves.
\end{exam}

\section{Conclusions and further comments}
\label{SEC:conclusion}

In the previous sections we have studied an adaptive moving mesh finite element method
for the numerical solution of the RLW equation. The RLW equation represents a class of PDEs
containing spatial-time mixed derivatives. For the numerical solution of those PDEs, a $C^0$
finite element method cannot be applied on a moving mesh since the mixed derivatives of the finite
element approximation may not be defined. To avoid this difficulty, a new variable (\ref{v-1}) was
introduced and the RLW equation was rewritten into a system of two coupled PDEs.
The system was then discretized in space using linear finite elements on a moving mesh
which is generated with a new implementation of the moving mesh PDE method.
The ODE system was integrated in time using the fifth-order Radau IIA scheme.

A range of numerical examples in one and two dimensions were presented. They include
the RLW equation with one or two solitary waves and special initial conditions that lead to
the undular bore and solitary train solutions. Numerical results have demonstrated that
the moving mesh finite element method has a second order convergence
as the mesh is being refined and is able to move and adapt
the mesh to the evolving features in the solution of the RLW equation. Moreover,
the method produces an error an order of magnitude smaller than that with a fixed mesh of the same
number of elements. 

It should be mentioned that the finite element approximation with both fixed and moving meshes
does not preserve $E_1$ (the mass) but the error quickly decreases to the level of roundoff error
as the mesh is refined. On the other hand, the moving mesh finite element method does a worse job to
conserve $E_2$ (the energy) than the fixed mesh finite element method although the former is more accurate.
It would be interesting to know what advantages the conservation of this quantity may give
the scheme for the RLW equation.
A major difficulty for the moving mesh method to conserve $E_2$ comes from the mesh movement,
which makes the mass matrix time dependent and introduces an extra convection term (see (\ref{vt-1})).
How to design a moving mesh method
that conserves this quantity will also be an interesting research topic.

\vspace{20pt}

{\bf Acknowledgment.} 
The work was partially supported by NSFC through grants
91230110, 11571290, and 41375115.
The authors are grateful to the anonymous referee for the valuable comments in improving the quality of the paper.


\end{document}